\documentclass[11pt,leqno,a4paper,psamsfonts]{amsart}
\usepackage[french]{babel}
\usepackage{amssymb,amsmath,amscd,enumerate,amsthm}
\usepackage[psamsfonts]{amsfonts}
\usepackage{latexsym}
\usepackage{amsbsy}
\usepackage[mathscr]{eucal}
\usepackage{graphicx}
\usepackage{enumerate}
\usepackage{verbatim}
\newcommand{\C}{{\mathbb{C}}}
\newcommand{\D}{{\mathbb{D}}}
\newcommand{\DD}{{\overline{\D}}}
\newcommand{\DDD}{{\DD\times\DD}}
\newcommand{\ssD}{{D^{\sharp}}}

\newcommand{\F}{{\mathcal{F}}}
\newcommand{\N}{\mathbb{N}}
\newcommand{\R}{{\mathbb{R}}}
\newcommand{\Z}{{\mathbb{Z}}}
\newcommand{\corjf}[1]{#1}
\newcommand{\correction}[1]{#1}
\newcommand{\II}{{|\!|\!|}}
\newcommand{\be}{{\mathbf{e}}}
\newcommand{\mf}[1]{{\mathfrak{#1}}}
\newcommand{\mc}[1]{{\mathcal{#1}}}
\newcommand{\mb}[1]{{\mathbb{#1}}}
\newcommand{\wt}[1]{{\widetilde{#1}}}
\let\wF\F
\newcommand{\FN}{{\F_{p_0 /q_0,\,k,\, \alpha }}}
\newcommand{\FNU}{{\F_{1,\,1,\, \alpha }}}
\newcommand{\wG}{\widehat{\mathcal{G}}}
\newcommand{\wh}[1]{{\widehat{#1}}}

\newcommand{\ucoF}{{\overset{\scriptscriptstyle 1}{\underset{\scriptscriptstyle\F}\looparrowright}}}

\newcommand{\uco}[1]{{\overset{\scriptscriptstyle
1}{\underset{\scriptscriptstyle{#1}}\looparrowright}}}
\let\ucon\ucoF

\let\loarF\ucoF

\newcommand{\zco}[1]{{\overset{\scriptscriptstyle 0}{\underset{\scriptscriptstyle{#1}}\looparrowright}}}%%%%%%%%%%%%%%%%%%%%%%%%%%%%%%%%%%%%%%%%%%%%%%%%%
\newcommand{\iso}{{\overset{\sim}{\longrightarrow}}}
\newcommand{\inte}[1]{\overset{\circ}{{#1}}}
\let\CAL=\mathcal%
\def\mathcal#1{{\CAL#1}}%
\newcounter{thmm}

\newtheorem*{thm}{\textsc{Th\'{e}or\`{e}me Principal}}
\newtheorem{teo}{Th\'{e}or\`{e}me}[subsection]
\newtheorem{lema}[teo]{Lemme}
\newtheorem{sublema}[teo]{Sous-Lemme}
\newtheorem{prop}[teo]{Proposition}
\newtheorem{defin}[teo]{D\'{e}finition}
\newtheorem{cor}[teo]{Corollaire}
\newtheorem{obs2}[teo]{Remarque}
\newenvironment{obs}{\begin{obs2}\rm}{\end{obs2}}
\newenvironment{dem}{\begin{proof}[Preuve]}{\end{proof}}
\newenvironment{dem2}[1]{\begin{proof}[Preuve #1]}{\end{proof}}
\def\bibartp#1#2#3#4#5#6#7#8
{\bibitem{#1} {\sc #2}, {\it #3}, {#4}, {\bf t. #5}, pages { #6} \`a
{ #7}, {({#8})}}
\def\bibart#1#2#3#4#5#6
{\bibitem{#1} {\sc #2}, {\it #3}, {#4}, {\bf #5}, {({#6})}}
\def\bibliv#1#2#3#4#5
{\bibitem{#1} {\sc #2}, {\it #3}, {#4}, {({#5})}}
\def\bibaart#1#2#3#4
{\bibitem{#1} {\sc #2}, {\it #3}, {#4}}
\begin{document}
\title[Incompressibilit\'{e} des feuilles  de feuilletages holomorphes]
{L'incompressibilit\'{e} des feuilles de
germes de feuilletages holomorphes singuliers}
\date{\today}
\author{David Mar\'{\i}n et Jean-Fran\c{c}ois Mattei}
\address{Departament de Matem\`{a}tiques \\
Universitat Aut\`{o}noma de
Barcelona \\
E-08193 Bellaterra (Barcelona)\\
Spain} \email{davidmp@mat.uab.es}
%\author{Jean-Fran\c{c}ois Mattei }
\address{Institut de Math\'{e}matiques\\ Laboratoire Emile Picard \\
Universit\'{e} Paul Sabatier\\
118, Route de Narbonne\\
F-31062 Toulouse Cedex 4, France} \email{mattei@picard.ups-tlse.fr}
\begin{abstract}
    We consider a non-dicritic germ of foliations $\mathcal{F}$  defined in
some ball ${\mathbb{B}\subset\mathbb{C}^{2}}$ with separatrix set
$S$, satisfying some additional but generic hypothesis. We prove
that there exists an open subset $U\supset S$ of $\mathbb{B}$  such
that for every leaf $L$ of $\mathcal{F}_{|(U\setminus S)}$ the
natural inclusion $\imath: L\hookrightarrow U\setminus S$ induces a
monomorphism $\imath_*:\pi_1( L)\hookrightarrow\pi_1(U\setminus S)$
at the fundamental group level.
\end{abstract}
\maketitle

\tableofcontents

\section*{Introduction}
Soit $\F_\omega $ un feuilletage holo\-morphe singulier d\'{e}fini par
une 1-forme diff\'{e}\-ren\-tielle $\omega$ \`{a} coefficients holomorphes
sur la boule ouverte $\mb B_{\varepsilon}\subset\C^2$ de centre $0$
et de rayon $\varepsilon>0$. Nous supposons $\omega$ \`{a} singularit\'{e}
isol\'{e}e en $0$ et {\textbf{non-dicritique}} \cite{MM}, i.e. les
germes \`{a} l'origine de courbes analytiques irr\'{e}ductibles $S_j$ telles
que $\omega_{\mid S_j} \equiv 0$, appel\'{e}es \textbf{s\'{e}paratrices de
$\F_\omega $}, sont en nombre fini : $ j = 1, \ldots,
\correction{\varrho}$ -et non-nul d'apr\`{e}s \cite{Camacho.Sad}. Nous
choisissons $\varepsilon_0$ assez petit pour que, dans la boule
ferm\'{e}e $\overline{\mb B}_{\varepsilon_0}$, \textbf{la s\'{e}paratrice
totale} $S:=  \corjf{\bigcup}_{j=1}^{{\varrho}} S_j$ soit analytique
ferm\'{e}e, \`{a} singularit\'{e} isol\'{e}e $0$ et transverse \`{a} chaque sph\`{e}re
$\partial\overline{\mb B}_r$, $0 <r \leq \varepsilon_0$. Fixons
aussi une fonction holomorphe r\'{e}duite $f$ \`{a} valeur dans le disque
$\D_{\eta'} := \left\{ z \in \C \; ; \; |z|<\eta' \right\}$ qui
d\'{e}finit $S$ sur $\overline{\mb
B_{\varepsilon_0}}$. \\

La restriction de $f$ \`{a} l'ouvert $${T}_\eta := f^{-1}(\D_\eta)\cap
\mb B_{\varepsilon_0}\,, \qquad 0<\eta\ll\varepsilon_0\,,$$ que nous
appellerons ici \textit{tube de Milnor}, est une fibration
diff\'{e}rentiable \cite{Milnor} localement triviale au dessus du disque
\'{e}point\'{e} $\D^\ast_\eta := \D_\eta \setminus \{0\}$. Notons
\begin{equation}\label{tubeMilnorCreux}
  {   T_\eta^*:= T_\eta \setminus
S = f^{-1}(\D_\eta^*)\cap \mb B_{\varepsilon}}{\,.}
\end{equation}
Lorsque $f$ est une int\'{e}grale premi\`{e}re de $\F_\omega $, pour chaque
feuille $L$ de la restriction ${\F_{\omega \,| T_\eta^*}}$ de
$\F_\omega $ \`{a} ${T_\eta^*}$ la suite exacte d'{homotopie} donne:
$${1\to\pi_1(L)\stackrel{\imath_*}{\to}\pi_1( T^*_\eta)\stackrel{f_*}{\to}\pi_1(\D_{\eta}^*)\to 1\,,}$$
o\`u ${\imath:L\hookrightarrow  T_\eta^*}$ d\'{e}signe l'inclusion
naturelle. En particulier, les feuilles de ${\F_{\omega \,|
T_\eta^*}}$ sont incompressibles dans {$T_\eta^*$}. L'objet de ce
travail est d\'{e}monter un
r\'{e}sultat analogue dans un cadre g\'{e}n\'{e}ral.\\

Consid\'{e}rons {$E: \mc T_\eta \to T_\eta$} le morphisme de r\'{e}duction
de {$\F_{\omega}$}, cf. \cite{Seidenberg} ou \cite{MM}. Le
\textbf{transform\'{e} total ${\mc D} := E^{-1}(S)$ de $S$}, que nous
appelons ici \textbf{diviseur total}, est \`{a} croisements normaux. Ses
composantes irr\'{e}ductibles sont : les composantes irr\'{e}ductibles $\mc
E_j$, $j= 1,\ldots, \kappa$ du \textbf{diviseur exceptionnel} $\mc E
:= E^{-1}(0)$ et les \textbf{transform\'{e}es strictes} $\mc S_j :=
\overline{E^{-1}(S_j) - \mc{D}_\omega }$ des $S_j\,$, $j=1,\ldots,
\varrho$. L'image r\'{e}ciproque $E^\ast\omega$ permet de d\'{e}finir sur
$\wt{B}_{\varepsilon_0,\, \eta_0}$ un feuilletage $\F$ \`{a}
singularit\'{e}es isol\'{e}es, dont le \textbf{lieu singulier} $Sing(\F)$
est contenu dans $\mc E$. En chaque point $c\in Sing(\F)$, le germe
$\F_c$ de $\F$ peut \^{e}tre d\'{e}crit par un germe de 1-forme $\wt \omega
_c$ qui s'\'{e}crit, dans des coordonn\'{e}es $z_1$, $z_2$ appropri\'{e}es :
\begin{equation}\label{formered}
    \wt \omega _c = (\lambda _cz_1 + \cdots)dz_2 + (\mu _cz_2 +
    \cdots)dz_1
\,,\quad \mathrm{avec}\quad \lambda _c\not=0\,, \mu _c / \mu _c
\notin \mb Q_{<0}{\,.}
\end{equation}
Nous dirons ici que $\F$ est \textbf{de type g\'{e}n\'{e}ral} si, pour
chaque $c\in Sing(\F)$ on a :
\begin{enumerate}[(i)]
\item$\lambda _c\, \mu _c \neq 0$
  \item si $\lambda _c/\mu _c$ est un r\'{e}el  irrationnel,
  alors le germe $\F_c$ est
  lin\'{e}arisable.
\end{enumerate}
Rappelons que si $\lambda _c/\mu _c$ est un irrationnel $<0$, ou
bien si appartient \`{a} un ensemble $\mf B\subset \R_+$ de mesure
pleine, appel\'{e} ensemble de Brjuno {\cite{Yoccoz}}, alors le germe
$\F_c$ est toujours lin\'{e}arisable.

\begin{thm} Soit $\F_\omega $ un germe \`{a} l'origine de $\mb C^2$ de
feuilletage holomorphe singulier non dicritique, de type g\'{e}n\'{e}ral et
soit {$T_{\eta_0}$} un tube de Milnor pour la s\'{e}paratrice totale
$S$. Alors il existe un voisinage $U$ de $S$ dans {$T_{\eta_0}$ tel}
que :
\begin{enumerate}
  \item l'inclusion $(U\setminus S) \hookrightarrow {(T_{\eta_0}\setminus S)}$
  induit un isomorphisme $\pi _1(U\setminus S, \mathbf{\cdot})\; \iso\; \pi _1
  {(T_{\eta_0} \setminus S, \mathbf{\cdot})}$,
  \item toute feuille $L$ de la
restriction $\F_{\omega \,|(U \setminus S)}$ est incompressible dans
$(U\setminus S)$, i.e. l'inclusion naturelle $\imath: L
\hookrightarrow (U \setminus S)$ induit un monomorphisme des groupes
fondamentaux $\imath_*:\pi_1(L \,, \,\mathbf{\cdot})
\hookrightarrow\pi_1(U \setminus S \,, \mathbf{\cdot})\,$,
\item $U$ contient {$T_{\eta}$} pour $\eta>0$
  assez petit.
\end{enumerate}
\end{thm}

Il est bien connu \cite{Milnor} que l'application d'inclusion de {$(
T_{\eta_0}\setminus S)$} dans $(\mb B_{\varepsilon_0}\setminus S)$
induit un isomorphisme au niveau du  groupe fondamental.
Ainsi le th\'{e}or\`{e}me pr\'{e}c\'{e}dent permet de construire un syst\`{e}me
fondamental $(U_n)_n$ de voisinage de $S $ dans la boule $\mb
B_{\varepsilon_0}$ tel que, pour tout $n$, chaque feuille de
$\F_\omega {}_{|(U_n\setminus S)}$ est incompressible dans $(\mb
B_{\varepsilon_0}\setminus S)$.\\

Remarquons que g\'{e}n\'{e}riquement les feuilles de $\F$ dont le groupe
fondamental est non-nul forment un ensemble dense. C'est  le cas
lorsque le groupe d'holonomie d'une composante du diviseur
exceptionnel est non r\'{e}soluble. Il existe alors \cite{lor.liou.bel}
un ensemble dense de points fixes attractifs d'\'{e}l\'{e}ments du
pseudo-groupe d'holonomie. Ces points correspondent n\'{e}cessairement \`{a}
des lacets trac\'{e}s dans une feuille homotopiquement non-trivaux. La
densit\'{e} pour la topologie de Krull {ce}
type de feuilletage est montr\'{e}e dans \cite{LeFloch}.\\

Les hypoth\`{e}ses que nous donnons ici peuvent \^{e}tre affaiblies; en
particulier {il est possible d'adapter l'\'{e}nonc\'{e} et la preuve de ce
r\'{e}sultat au cas \corjf{des feuilletages dicritiques. Afin d'all\'{e}ger
ce texte nous avons pr\'{e}f\'{e}r\'{e} ne consid\'{e}rer ici que le cas
non-dicritique et traiter du cas dicritique dans un article
ult\'{e}rieur \cite{david.jf.2} o\`{u} l'analyse est faite dans le cadre
plus g\'{e}n\'{e}ral des {feuilletages holomorphes} au voisinage d'un
diviseur compact quelconque.}}\\

La structure de l'article est la suivante : Au chapitre 1 nous
intro\-duisons une notion tr\`{e}s g\'{e}n\'{e}rale, la $1$-connexit\'{e}
feuillet\'{e}e, qui permet de ``localiser'' le probl\`{e}me, gr\^{a}ce \`{a} un
th\'{e}or\`{e}me de type Van Kampen (\ref{teo4}). On obtient en suite au
chapitre 2 une {technique} d'assemblage feuillet\'{e}, qui est le
pendant feuillet\'{e} des techniques de plombage \cite{EisenbudNeumann}.
Cependant un ``algorithme de plombage'' ne sera possible qu'en
contr\^{o}lant la rugosit\'{e} des bords des blocs, notion introduite en
(\ref{subsecrugos}). Dans ce chapitre nous ramenons la preuve du
{th\'{e}or\`{e}me} principal \`{a} celle d'un th\'{e}or\`{e}me d'existence de {blocs}
``adaptables \`{a} rugosit\'{e} contr\^{o}l\'{e}e'' (\ref{teoexistence}). La preuve
de celui-ci occupera les deux autres chapitres.

\section{Notions de connexit\'{e}
feuillet\'{e}e}\label{connfeuillete}
\subsection{Notions de $0$- et
$1$-connexit\'{e} feuillet\'{e}e}\label{ssectconfeuillete} De mani\`{e}re
g\'{e}n\'{e}rale nous consid\'{e}rons une vari\'{e}t\'{e} diff\'{e}rentiable $M$ munie d'un
feuilletage r\'{e}gulier $\F$ de classe $C^{1}$. On notera aussi $\F$
l'ensemble des feuilles du feuilletage. Pour tout sous-ensemble $A$
de $M$, $\F_{|A}$ d\'{e}signera la collection des composantes connexes
des intersections  $L\cap A$, $L\in\F$, des feuilles de $\F$ avec
$A$. Pour $A\subset B$ nous d\'{e}signons  par $Sat_\F( A, B)$ et
appelons \textbf{satur\'{e} de $A$ par $\F_{|B}$}, le sous-ensemble
\begin{equation}\label{satuteAB}
    Sat_\F( A, B) := \bigcup_{L\in \mathfrak{A}} L \subset B\,,\qquad
    \mathfrak{A}:= \{{L\in \F_{|B}\;/ \;L\cap A\neq
    \emptyset}\}\,.
\end{equation}

Soient $A,B$ deux sous-ensembles de $M$, avec $A\subset B$.
\begin{defin}
Nous dirons que $A$ est \textbf{$0$-$\F$-connexe dans $B$} et nous
noterons $A \zco{\F} B$, si pour tout $L\in\F_{|B}$ l'application
$\pi_{0}(L\cap A)\to\pi_{0}(A)$ induite par l'inclusion de $L\cap A$
dans $A$ est injective, {i.e. pour tout chemin $\alpha :[0,1]\to A$
d'extr\'{e}mit\'{e}s dans $L\cap A$, il existe un chemin $\beta :[0,1]\to
L\cap A$ de m\^{e}mes extr\'{e}mit\'{e}s que $\alpha $.} Nous dirons aussi que
$A$ est \textbf{strictement $0$-$\F$-connexe dans $B$}, si $A$ est \`{a}
la fois $0$-$\F$-connexe dans $B$ et \textbf{incompressible dans}
$B$, i.e. l'inclusion naturelle de $A$ dans $B$ induit un morphisme
injective des groupes fondamentaux
$\pi_{1}(A,p)\hookrightarrow\pi_{1}(B,p)$ pour chaque $p\in A$.
\end{defin}
Pour $K\subset B$, notons $\mathcal{H}(K)$ l'ensemble des classes
d'homotopie des chemins trac\'{e}s dans $K$,  i.e. $\mathcal{H}(K)$ est
le quotient de $C^{0}([0,1],K)$ par la relation d'\'{e}quivalence
$a\sim_{K}b$ $:\Leftrightarrow$ il existe $H:[0,1]^2 \to K$ continue
telle que : $H(0,t)=a(t)$, $H(1,t)=b(t)$, $H(s,0)=a(0)=b(0)$ et
$H(s,1)=a(1)=b(1)$   pour tout $s,t\in[0,1]$.
\begin{defin}
Nous dirons que $A$ est \textbf{$1$-$\F$-connexe dans $B$} et nous
noterons $A\uco{\F} B$, si pour toute feuille $L\in\F$, $L\cap
A\neq\emptyset$, la suite suivante est exacte :
$$\mathcal{H}(L\cap A)\stackrel{\alpha}{\to}\mathcal{H}(A)\times\mathcal{H}(L\cap B)\begin{array}{c}\stackrel{\beta_{1}}{\to}\\ \stackrel{\textstyle\to}{{\scriptstyle\beta_{2}}}%{\to}
\end{array}\mathcal{H}(B)\,,$$
avec $\alpha([c]_{L\cap A}):=([c]_{A}, \,[c]_{L\cap B})\,$,  $
\beta_{1}([a]_{A},[b]_{L\cap B}):=[a]_{B}\,$ et $
\beta_{2}([a]_{A},[b]_{L\cap B}):=[b]_{B}$. Nous dirons aussi que
$A$ est \textbf{strictement $1$-$\F$-connexe dans $B$} si $A\uco{\F}
B$ et $A$ est incompressible dans $B$.
\end{defin}
\noindent Explicitement cela signifie que pour tout chemins
 $a:[0,1]\to A$ et $b:[0,1]\to L\cap B$ tels que $a\sim_{B}b$, il existe
un chemin $c:[0,1]\to L\cap A$ tel que $c\sim_{A}a$ et $c\sim_{L\cap
B}b$.
\begin{obs} Soit $(M,\F)$ une vari\'{e}t\'{e} feuillet\'{e} et $A\subset B$ des
sous-ensembles de $M$. Les assertions suivantes sont \'{e}videntes :
\begin{enumerate}[(i)]
\item La r\'{e}lation $A\ucoF B$ entre parties de $M$ est transitive :
$$A\loarF B\quad\textrm{et}\quad B\loarF C\quad\Longrightarrow\quad
A\loarF C.$$ Il en est de m\^{e}me des relations de $1$-$\F$-connexit\'{e}
stricte, de $0$-$\F$-connexit\'{e} et de $0$-$\F$-connexit\'{e} stricte.
\item Si $A$ est r\'{e}duit \`{a} un point, $A=\{p\}$, alors $A\loarF B$ si et
seulement si la feuille $L_{p}$ de $\F_{|B}$ qui passe par $p$ est
incompressible dans $B$ :
$\imath_\ast:\pi_{1}(L_{p},p)\hookrightarrow\pi_{1}(B,p)$, o\`{u}
{$\imath_{*}$} d\'{e}signe le morphisme induit par l'application
d'inclusion {$\imath$} de $L_{p}$ dans $B$.
\item  $A$ est $1$-$\F$-connexe dans $B$ si et seulement si
chacune de ses composantes connexes est $1$-$\F$-connexe dans $B$,
\item Si $A$ est $\F$-satur\'{e} dans $B$ alors $A\uco{\F}B$.

\end{enumerate}
\end{obs}
\begin{prop}\label{prop3}
Soit $\F$ un feuilletage sur $M$ de dimension r\'{e}elle $1$, d\'{e}finit
par un champ de vecteurs ne poss\'{e}dant pas d'orbite p\'{e}riodique, alors
tout ouvert $A$  strictement $0$-$\F$-connexe dans $M$ est
$1$-$\F$-connexe dans $M$.
\end{prop}
\begin{dem}
Soient $a:[0,1]\to A$ et $b:[0,1]\to L\cap B$ tels que $a\sim_{B}b$.
Gr\^{a}ce a la $0$-$\F$-connexit\'{e} de $A$ dans $B$, il existe $c:[0,1]\to
L\cap A$ de m\^{e}mes extr\'{e}mit\'{e}s que $a$ (et que $b$). Comme $\F$ ne
poss\`{e}de pas d'orbites p\'{e}riodiques, chaque composante connexe de
$L\cap B$ est simplement connexe et l'on a : $c\sim_{L\cap B}b$.
D'o\`{u} $c\sim_{B}a$. Les images de $a$ et de $b$ sont contenues dans
$A$ et $A$ est incompressible dans $B$. Ainsi $c\sim_{A}a$.
\end{dem}

\vspace{1em}

\subsection{Un th\'{e}or\`{e}me de type Van Kampen} Soit $T$ une sous-vari\'{e}t\'{e} ferm\'{e}e de codimension r\'{e}elle
$1$, non n\'{e}ces\-sai\-rement connexe, transversalement orientable,
d'une vari\'{e}t\'{e} diff\'{e}rentiable $M$ et soit $\F$ un feuilletage
(r\'{e}gulier) de classe $C^{1}$ sur $M$. Consid\'{e}rons deux
sous-ensembles connexes $V_{1}, V_{2}$ de $M$ dont l'intersection
est $T$ et dont l'union $V$ est un voisinage de $T$ :
$$T=V_{1}\cap V_{2}\subset W\subset V:=V_{1}\cup V_{2},\quad W\subset M \textrm{ ouvert}.$$
\begin{teo}\label{teo4}
Supposons que $T$ est transverse \`{a} $\F$ et est strictement
$1$-$\F$-connexe dans $V_{1}$ ainsi que dans $V_{2}$ :
$$T  \pitchfork\F,\quad T\loarF V_{i},\quad \pi_{1}(T,\ast)\hookrightarrow\pi_{1}(V_{i},\ast),\quad i=1,2.$$
Alors chaque $V_{1}$ et $V_{2}$ sont tous deux strictement
$1$-$\F$-connexes dans $V$ :
$$V_{i}\loarF V,\quad \pi_{1}(V_{i},\ast)\hookrightarrow\pi_{1}(V,\ast),\quad i=1,2.$$
\end{teo}
\noindent La preuve de ce th\'{e}or\`{e}me consiste en une ``chirurgie''
d'homotopie. Pr\'{e}cisons cette technique.\\

Soit $H: [0,1]\times[0,1]\to V$ une homotopie entre deux chemins
$\gamma _0 := H(0, \cdot)$ et $\gamma _1 := H(1, \cdot)$.

\begin{defin}\label{homotbonne}
Nous dirons ici que $H$ est \textbf{$T$-admissible} si  $H^{-1}(T)$
est une union finie disjointe de courbes de classe $C^1$ ferm\'{e}es et
proprement plong\'{e}s dans $[0,1]\times [0,1]$.
\end{defin}

\begin{obs}\label{genericite}
Si $H$ est $C^1$ et transverse \`{a} $T$ et \`{a} $\partial T$, alors $H$
est $T$-admissible. Ainsi lorsque $\gamma _0$ et $\gamma _1$ sont
$C^1$ et transverses \`{a} $T$ et \`{a} $\partial T$, on peut approcher $H$
par une homotopie diff\'{e}rentiable qui est g\'{e}n\'{e}riquement
$T$-admissible.
\end{obs}

Supposons que $H$ est $T$-admissible. D\'{e}signons par $\mc I(H; T)$
l'ensemble des composantes connexes de $H^{-1}(T)$ hom\'{e}omorphes \`{a} un
segment et par $\mc J(H;T)$ l'ensemble des composantes connexes de
$H^{-1}(T)$ qui sont des courbes de Jordan.  Munissons ces courbes
de leur param\'{e}trisation par abscisse curviligne. Si $\delta \in \mc
J(H;T)$ ne coupe pas le \textbf{bord horizontal} $[0,1]\times \{0,
1\}$ de $[0,1]\times [0,1]$, nous noterons $\Delta_\delta
^\mathrm{ext} $, resp. $\Delta_\delta ^\mathrm{int}$ la
\textbf{composante connexe ext\'{e}rieure}, resp. \textbf{int\'{e}rieure} de
$([0,1]\times [0,1]) \setminus |\delta |$, c'est \`{a} dire celle qui
contient, resp. qui ne contient pas $[0,1]\times \{0, 1\}$. D'autre
part nous dirons que $\theta \in \mc I(H;T)$ est un \textbf{\'{e}l\'{e}ment
extr\'{e}\-mal de} $\mathcal{I}(H, T)$, si ses extr\'{e}mit\'{e}s sont situ\'{e}es
sur $\{0\}\times [0,1]$; nous les noterons $(0,s_{0}(\theta))$,
$(0,s_{1}(\theta))$, et nous orientons $\theta$ pour avoir
$s_{0}(\theta)\le s_{1}(\theta)$. Le lacet simple
$\delta_\theta:=\theta\vee\alpha_{\theta}$ borde un disque conforme,
$\alpha_\theta$ d\'{e}signant la param\'{e}trisation naturelle de
$I_\theta:=[s_{0}(\theta),s_{1}(\theta)]$. Nous notons $\Delta
_\theta^{\mathrm{int}}$ l'int\'{e}rieur de ce disque et $\Delta
_\theta^{\mathrm{ext}} := ([0,1]\times [0,1])\setminus
\overline{\Delta} _\theta^{\mathrm{int}}$.

\begin{defin}\label{chirurgie}
Nous dirons qu'une homotopie $H'$ entre les chemins $a$ et $b$ est
\textbf{obtenue \`{a} partir de $H$ par chirurgie le long d'un lacet
simple $\delta' $} qui ne coupe pas le bord horizontal de
$[0,1]\times [0,1]$, si l'on a l'\'{e}galit\'{e} des restrictions :
$$H_{|\Delta_{\delta'} ^\mathrm{ext}} = H'_{|\Delta_{\delta'}
^\mathrm{ext}}.$$
\end{defin}

\begin{center}
\begin{figure}[htp!]
\includegraphics[width=12cm]{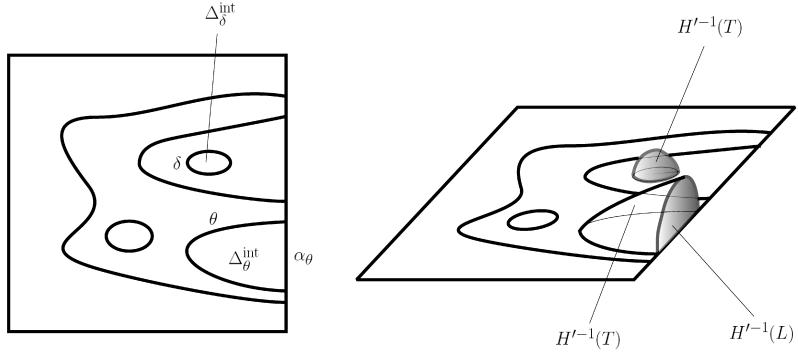}
\vspace{-12cm} \caption{Chirurgie d'homotopies.}
\end{figure}
\end{center}

\noindent On montre tr\`{e}s facilement les assertions suivantes :

\begin{lema}\label{chirgenerale}
Soit $\Omega$ une sous-vari\'{e}t\'{e} (non n\'{e}cessairement ferm\'{e}e) contenue
dans l'adh\'{e}rence d'une composante connexe de $V\setminus T$.
\begin{enumerate}[(i)]
\item Si  $\delta \in \mathcal{J}(H, T)$ v\'{e}rifie $H(|\delta
|)\subset \Omega $ et est homotope dans $\Omega $ \`{a} un point, alors
il existe une homotopie $H'$ obtenue par chirurgie le long de
$\delta $ telle que
$H'(\Delta_\delta ^\mathrm{int})\subset \Omega$.\\
\item  Soit $\theta $ un chemin extr\'{e}mal v\'{e}rifiant $H(|\theta |)
\subset \Omega $ et homotope dans $\Omega $ a un chemin $\eta $.
Consid\'{e}rons le chemin $a'$ \'{e}gal \`{a} $\eta $ en restriction \`{a} $I_\theta
$ et \'{e}gal \`{a} $a$ en restriction \`{a} $[0,1]\setminus I_\theta $. Alors
il existe une homotopie $H'$ entre $a$ et $a'$, obtenue par
chirurgie le long de $\theta  $, telle que $H'(\Delta
_\theta^{\mathrm{ext}})\subset \Omega $.
\end{enumerate}
\end{lema}

\begin{dem2}{du th\'{e}or\`{e}me \ref{teo4}}
Montrons d'abord les relations $V_{i}\loarF V$, $i=1,2$. Consid\'{e}rons
dans $V$ une homotopie $H$ entre un chemin $a$ trac\'{e}e dans $V_{i}$
dont les extr\'{e}mit\'{e}s sont dans une m\^{e}me feuille $L$ de $\F_{|V}$ et
un chemin $b$ trac\'{e}e dans $L$. La feuille $L$ est transverse \`{a} $T$
et quitte \`{a} effectuer une petite homotopie dans $L$ nous supposons
$b$ diff\'{e}rentiable et transverse \`{a} $T$. Gr\^{a}ce \`{a} la remarque
(\ref{genericite}) nous supposons aussi que $H$ est $T$-admissible.
Raisonnons par r\'{e}curence sur l'entier $N(H, T) := \# \mathcal{J}(H,
T) + \# \mathcal{I}(H,T)$. Lorsque $N(H, T) =0$, le chemin $b$ est
contenu dans $V_i$ et le r\'{e}sultat est trivial. Supposons $N(H, T) =
N+1$ et le th\'{e}or\`{e}me v\'{e}rifi\'{e} pour $N(H, T) \leq N$. Distinguons deux
cas suivant que  $\mathcal{J}(H,T)$ est vide ou non.\\

\indent \textit{\underline{- Premier cas : $\mathcal{J}(H, T)\neq
\emptyset$}. } Il existe visiblement une courbe de Jordan $\delta
\in \mathcal{J}(H, T)$, telle que ${\Delta}_{\delta}^\mathrm{int}$
n'intersecte aucun \'{e}l\'{e}ment de $\mathcal{I}(H, T)\cup \mathcal{J}(H,
T)$. Visiblement on a :
$$H(\delta)\subset T,\quad H(\Delta^{\mathrm{int}}_{\delta})\subset V_{j}\quad\textrm{et}\quad  H(\Delta_{\delta}^{\mathrm{ext}})\subset V_{k},\quad\textrm{avec}\quad \{j,k\}=\{1,2\}.$$
Comme $T$ est incompressible dans $V_{j}$, le lacet $\delta $ est
homotope \`{a} un point dans $T$. Le lemme \ref{chirgenerale} (i) donne
une homotopie $H'$ telle que
$H'(\overline{\Delta}_{\delta}^\mathrm{int})\subset T$. De plus
$\{\overline{\Delta}_{\delta}^\mathrm{int}\}\cup\mathcal{I}(H,
T)\cup \mathcal{J}(H, T)\setminus\{\delta \}$ est la collection des
composantes connexes de ${H'}^{-1}(T)$. Comme $T$ est transverse \`{a}
$\F$ et transversalement orientable, on peut se donner au voisinage
de $H'({\Delta}_{\delta}^\mathrm{int})$  un champ de vecteurs
diff\'{e}rentiable $Z$ \`{a} support compact tel que :
\begin{enumerate}
    \item[(i)] $H'({\overline{\Delta}}_{\delta}^\mathrm{int})\subset
    \mathrm{supp}(Z)\subset W$,
    \item[(ii)] $Z(z)$ est transverse \`{a} $T$ et pointe vers $V_{k}$,
    pour $z\in H'({\overline{\Delta}}_{\delta}^\mathrm{int})$,
    \item[(iii)] $\mathrm{supp}(Z)\cap H(\lambda) = \emptyset$ pour tout $\lambda
\in \mathcal{I}(H, T)\cup \mathcal{J}(H, T)\setminus\{\delta \}$.
\end{enumerate}
Le compos\'{e} $H'_\epsilon $ de $H'$ avec le flot de $Z$ est une
homotopie entre $a$ et $b $ qui v\'{e}rifie : $\mathcal{I}(H_\epsilon ,
T) = \mathcal{I}(H, T)$ et $\mathcal{J}(H_\epsilon , T) =
\mathcal{J}(H, T) \setminus \{\delta \}$ pour un temps $\epsilon >0$
assez petit. Ceci ach\`{e}ve la r\'{e}curence, dans ce cas.\\

\indent \textit{\underline{- Deuxi\`{e}me cas : $\mathcal{J}(H, T)=
\emptyset$ et $\mathcal{I}(H, T)\neq \emptyset$ . }} Il existe alors
un \'{e}l\'{e}ment extr\'{e}mal $\theta $ de $\mathcal{I}(H, T)$ tel que $\Delta
_\theta^{\mathrm{int}}$ n'intersecte aucun \'{e}l\'{e}ment de
$\mathcal{I}(H, T)$. Notons encore $V_{j}$ celui des deux ensembles
$V_1$ ou $V_2$ qui contient $H(\Delta _\theta^{\mathrm{int}})$ et
notons $V_{k}$ l'autre. Comme $T\loarF V_{j}$ le chemin $\alpha
_\theta $ est homotope dans $L\cap V_j$ \`{a} un lacet $\eta $ contenu
dans $T$. L'injection de $\pi_1(T)$ dans $\pi_1 (V_{j})$ donne une
homotopie dans $T$ entre $\eta $ et $\alpha_\theta $. On conclu
gr\^{a}ce au lemme \ref{chirgenerale} (ii) et \`{a} l'int\'{e}gration d'un champ
de vecteurs, comme dans le premier cas.\\

Il reste a montrer les injections des groupes fondamentaux de $\pi
_1(V_i)$ dans $\pi _1(V)$, $i=1$, $2$. Elles r\'{e}sultent en fait d'un
th\'{e}or\`{e}me classique de combinatoire des groupes en consid\'{e}rant $\pi
_1(V)$ comme le produit de $\pi _1(V_1)$ et $\pi _1( V_2)$ amalgam\'{e}e
sur $\pi _1(T)$, cf.\cite{Magnus}, th\'{e}or\`{e}me 4.3. page 199.
\end{dem2}

\vspace{1em}

{\section{Assemblage bord \`{a} bord feuillet\'{e}}}
%
%Vari\'{e}t\'{e}s feuillet\'{e}es construites par assemblage bord \`{a}
%bord}
\subsection{{D\'{e}finition d'assemblage feuillet\'{e} et th\'{e}or\`{e}me de localisation}}
%\subsection{Assemblage feuillet\'{e}}
Soit $V$ une vari\'{e}t\'{e} \`{a} bord r\'{e}elle de dimension 4,
non-n\'{e}cessairement compacte ni connexe, munie d'un feuilletage
r\'{e}gulier $\F_V$ de classe $C^1$ et $\left( V_j\right)_{j\in I}$,
$I\subset \N$, une famille de sous-vari\'{e}t\'{e}s de $V$ de m\^{e}me dimension
que $V$, dont les bords $\partial V_j := V_j\setminus \inte{V_j}$
sont transversalement orientables.
\begin{defin}\label{adaptable} Nous dirons qu'un \'{e}l\'{e}ment $V_j$ est
un \textbf{bloc feuillet\'{e} $\F_V$-adaptable} s'il satisfait les
propri\'{e}t\'{e} suivantes :
\begin{enumerate}
  \item\label{bordincompr} chaque composante connexe  de $\partial V_j$ est
  incompressible dans $V_j$,
  \item\label{bordtrans} $\F_V$ est tranverse \`{a} $ \partial V_j$,
  \item\label{feuilleincompr} chaque feuille de $\F_{|V_j}$ est incompressible dans $V_j$,
  \item\label{feuilleucon} chaque composante connexe de $\partial V_j$ est
  $1$-$\F$-connexe dans $V_j$.
\end{enumerate}
\noindent Nous dirons que $V$ est un \textbf{assemblage bord \`{a} bord
des $V_j$ bien construit}, si pour tout $j\in I$ la condition
(\ref{bordincompr}) et la condition suppl\'{e}mentaire suivante sont
satisfaites :
\begin{enumerate}\setcounter{enumi}{4}
  \item\label{bordabord} pour tout $i$, $j \in I$ distincts l'une des
deux \'{e}ventualit\'{e}s suivante est r\'{e}alis\'{e}e : ou bien $V_i \cap V_j
=\emptyset\,$, ou bien $V_i \cap V_j$ est  une composante connexe de
$\partial V_i$ et est une composante connexe de $\partial V_j\,$.
\end{enumerate}
Si chaque $V_j$ est $\F_V$-adaptable et si $V$ est un assemblage
bord \`{a} bord des $V_j$, nous dirons que $V$ est un \textbf{assemblage
bord \`{a} bord feuillet\'{e} bien construit}.
\end{defin}

\begin{obs}\label{rVK}
Une cons\'{e}quence imm\'{e}diate du th\'{e}or\`{e}me classique de Seifert-Van
Kampen est la suivante. Soient $A\subset A'$, $B\subset B'$ et
$A\cap B\subset A'\cap B'$ deux sous-ensembles connexes d'un espace
topologique tels que chaque inclusion pr\'{e}c\'{e}dente induise un
isomorphismes au niveau du groupe fondamental. Alors il est de m\^{e}me
pour l'inclusion  $A\cup B\subset A'\cup B'$.
\end{obs}

Par r\'{e}curence sur le nombre de blocs, on d\'{e}duit de cette remarque la
proposition suivante :

\begin{prop}\label{assemblages}
Soient $V = \bigcup_{j\in I} V_j$ et  $V' =\bigcup_{j\in I} V'_j$,
$\# I<\infty$, deux assemblages bord \`{a} bord  bien construits.
Supposons que pour chaque $j\in I$, $V'_j$ est contenu dans $V_j$ et
que l'application d'inclusion induise un isomorphisme $\pi _1(V'_j,
\, \cdot)\iso \pi _1(V_j, \cdot)$. Alors l'application d'inclusion
$V'\hookrightarrow V$ induit aussi un isomorphisme $\pi _1(V', \,
\cdot)\iso \pi _1(V, \cdot)$.
\end{prop}

\begin{obs}\label{reassemblage}
%Notons $\mathfrak{C}_j$ l'ensemble des composantes connexes $W$ du
%bord de $V_j$, telles qu'il existe un indice $r\neq j$ avec $W\cap
%V_r \neq \emptyset$.
%Posons $$V_j^{\diamond} \; : =\; V_j \setminus
%(\bigcup_{W\in \mf C_j} \,W) $$
Supposons que les $V_j$ v\'{e}rifient
les propri\'{e}t\'{e}s (\ref{bordincompr}) \`{a} (\ref{feuilleucon}) et la
propri\'{e}t\'{e} suppl\'{e}\-men\-taire suivante :
%\begin{enumerate}
% \item [(5')] Pour tout $i, j \in I$ distincts $V_i\cap V_j$ est ou
% bien vide, ou bien {\'{e}gal \`{a}} $\partial V_i\cap \partial V_j$ et
%il existe alors une sous-vari\'{e}t\'{e} ouverte
%$V'_{ij}$ de $V_i\cap V_j$ qui est incompressible et
%$1$-$\F$-connexe dans $V_i\setminus$ et dans $V_j$.
%\end{enumerate}
{\begin{enumerate}
  \item [(5')] \it Pour tout $i, j \in I$ distincts $V_i\cap V_j$ est ou
  bien vide, ou bien {\'{e}gal \`{a}} $\partial V_i\cap \partial V_j$ et
  il existe alors une sous-vari\'{e}t\'{e} ouverte
$V'_{ij}$ de $V_i\cap V_j$ qui est incompressible et
$1$-$\F$-connexe dans $(V_i\setminus V_i\cap V_j)\cup V'_{ij}$ et
dans $(V_j\setminus V_i\cap V_j)\cup V'_{ij}$.
\end{enumerate}
Consid\'{e}rons les vari\'{e}t\'{e}s \`{a} bord
$$V_j^{\diamond} \; : =\; \left( V_j
\setminus \cup_{j\neq i} \,V_j\cap V_i \right)\cup  \left(\cup_{(i,
j)\in I\widehat{\times }I} V'_{ij}\right) \; \subset \; V_j\,,\quad
j\in I\,,
$$ o\`{u} $I\widehat{\times }I$ est l'ensemble des couples $(i, j)$ tels que
$V_i\cap V_j\neq \emptyset$.} Alors $V' := \bigcup_{j\in I}
V_j^\diamond$ est un assemblage bord \`{a} bord feuillet\'{e} bien construit
des $V_j^\diamond$.
\end{obs}

\begin{teo}[\bf de localisation]\label{teo.localisation}
Si $V$ est un assemblage bord \`{a} bord feuillet\'{e} bien construit, alors
:
\begin{enumerate}
  \item chaque feuille de $\F_V$ est incompressible dans $V$,
  \item chaque union de {blocs} $V' :=\bigcup_{k\in I'}
  V_k$, $I'\subset I$, est $1$-$\F_V$-connexe dans V,
  \item chaque union de {blocs} est incompressible dans V.
\end{enumerate}
\end{teo}
\begin{dem}Les deux premi\`{e}res assertions sont des cons\'{e}quences imm\'{e}diates
du th\'{e}or\`{e}me (\ref{teo4}). La troisi\`{e}me assertion r\'{e}sulte du th\'{e}or\`{e}me
de Van Kampen classique pr\'{e}c\'{e}demment cit\'{e} : \cite{Magnus}, th\'{e}or\`{e}me
4.3, page 199.
\end{dem}

Pour prouver le Th\'{e}or\`{e}me Principal \'{e}nonc\'{e} dans l'introduction, nous
construirons un voisinage ouvert $\wt U$ du diviseur total $\mc D$
tel que $\wt U^* := \wt U \setminus \wt D$ se d\'{e}compose comme un
assemblage bord \`{a} bord feuillet\'{e} bien construit, $\wt U^*=
\bigcup_\alpha \mc B_\alpha $. Le th\'{e}or\`{e}me de localisation permet
imm\'{e}diatement de conclure. Cette construction sera effectu\'{e}e \`{a}
partir d'une d\'{e}composition appropri\'{e}e du diviseur qui sera
{d\'{e}crite} %effectu\'{e}e
en (\ref{decompdiv}). Certains blocs, correspondant \`{a} des
singularit\'{e}s de $\mc D$, sont d\'{e}j\`{a} d\'{e}crits dans les exemples
ci-apr\`{e}s.

\subsection{Exemples utiles}\label{subsect.exemples}
Consid\'{e}rons un germe  $\mc G$ de feuilletage singulier r\'{e}duit \`{a}
l'origine de $\C^2$ donn\'{e} par un champ de vecteurs holomorphe $X :=x
\frac{\partial\phantom{x}}{\partial x} + y
 (\lambda  + \cdots)\frac{\partial\phantom{x}}{\partial y}$,
d\'{e}fini au voisinage {de}
%du sous-ensemble ferm\'{e}
$\mb K^{\zeta} :=
\{|x|\leq 1 ,\,|y|\leq 1,\,|xy|<\zeta \}$, avec $0<\zeta<1$.
D\'{e}signons par $\Phi^X_t$ le flot de $X$ et, pour tout sous-ensemble
$A$ de $\mb K^{\zeta}$, nous notons
\begin{equation}\label{notast}
    A^* := A\setminus \{xy=0\}\,,\quad \partial_1 A := A\cap \{|x| =
    1\}\,,\quad \partial_2 A := A\cap \{|y| =
    1\}\,.
\end{equation}
Supposons $\mc G$ \`{a} singularit\'{e} isol\'{e}e. Soit $\Sigma \subset
\{x_0\}\times \D(\zeta)$, avec $|x_0| = 1$, un disque conforme
ferm\'{e}, de bord analytique par morceaux. Consid\'{e}rons la
\textbf{suspension de $\Sigma $} au dessus du cercle $\{|x| = 1,\,
y=0\}$
$$
U_\Sigma  := \{\Phi^X_{it}(x_0, y)\;/\; 0\leq t\leq 2\pi ,\, (x_0,y)
\in \Sigma \}.
$$
Cet ensemble est bien d\'{e}fini et contenu dans le tore $T_{1} := \{|x|
=1, \, |y| < 1\}$ si $\Sigma $ est assez petit, ce que nous
supposerons. Nous supposons aussi que $\Sigma $ et son image
$h(\Sigma )$ par l'application d'holonomie $ h := \Phi^X_{2\pi i}$
sont \'{e}toil\'{e}s. Alors l'ensemble $U_\Sigma $ est un r\'{e}tract par
d\'{e}formation du tore $T_{\zeta}$. De plus il se r\'{e}tracte par
d\'{e}formation sur tout tore $T_\varepsilon$, avec $\varepsilon>0$
assez petit. En particulier $U_\Sigma^* $ est incompressible dans
$K^{\zeta *}$.\\

Nous d\'{e}sirons construire des blocs adaptables $\mc B $ contenus dans
$\mb K^{\zeta\,*}$ dont les adh\'{e}rences sont des voisinages des axes
dans $\mb K^\zeta$

\subsubsection{{\bf Exemple 1 :}} Supposons le champ $X$ lin\'{e}aire :
$X=x\frac{\partial}{\partial x}+\lambda y\frac{\partial}{\partial
y}$ avec $\lambda\in\C^*$. On a
\begin{equation*}
\partial_{i}\mb K^{\zeta\, *}
\uco{\mc G}\mb K^{\zeta\,
*},\qquad i=1,2.
\end{equation*}
En effet, consid\'{e}rons le sous-ensemble $\widetilde{\mb K}^{\zeta}$
des points $(z,w)$ dans $\C^{2}=\mathbb{R}^{4}$ qui v\'{e}rifient les
in\'{e}galit\'{e}s suivantes :
$$
\mathrm{Re}(z)\leq 0,\quad \mathrm{Re}(w)\leq 0\quad
\textrm{et}\quad \mathrm{Re}(z+w)<\log(\zeta).
$$
Consid\'{e}rons le rev\^{e}tement universel $\rho:\widetilde{\mb K}^*\to \mb
K^{\zeta\, *}$ d\'{e}finit par $\rho(z,w)=(e^{z},e^{w})$. Le feuilletage
$\widetilde{\mc G}:=\rho^*\mc G$ admet l'int\'{e}grale premi\`{e}re lin\'{e}aire
$w-\lambda z$.
Soient $a:[0,1]\to \partial_{i} \mb K^{\zeta\, *}$ et $b:[0,1]\to
L$, deux chemins homotopes dans $\mb K^{\zeta\, *}$ \`{a} extr\'{e}mit\'{e}s
fixes, o\`{u} $L$ est une feuille de $\mc G$ restreint \`{a} $\mb K^{\zeta\,
*}$. On peut choisir $\widetilde{a}:[0,1]\to
\partial_{i}\widetilde{\mb K}^{\zeta}
:= \rho^{-1}(\partial_{i}\mb K^{\zeta\, *})$ et
$\widetilde{b}:[0,1]\to\widetilde{L}$ deux $\rho$-rel\`{e}vements de $a$
et $b$ avec les m\^{e}mes extr\'{e}mit\'{e}s,  $\widetilde{L}$ \'{e}tant une feuille
de $\widetilde{\mc G}$ telle que $\rho(\widetilde{L})=L$. Comme
$\widetilde{L}$ est l'intersection d'un 2-plan avec ces trois
demi-espaces de $\mathbb{R}^{4}$, l'intersection $\widetilde{L}\cap
\partial_{i} \widetilde{\mb
K}^\zeta\neq\emptyset$ est soit une droite, soit une demi-droite. Il
est donc clair qu'il existe un chemin lin\'{e}aire
$\widetilde{c}:[0,1]\to \widetilde{L}\cap\partial \widetilde{\mb
K}^{\zeta}$ joignant les deux extr\'{e}mit\'{e}s communes de $\widetilde{a}$
et $\widetilde{b}$. Comme $\partial\widetilde{\mb K}^{\zeta}$ et
$\widetilde{L}$ sont simplement connexes, $\widetilde{c}$ est
homotope \`{a} $\widetilde{a}$ dans $\partial_{i} \widetilde{\mb
K}^\zeta$ et \`{a} $\widetilde{b}$ dans $\widetilde{L}$. Si on d\'{e}finit
$c:=\rho(\widetilde{c})$ on peut redescendre ces homotopies \`{a} $\mb
K^{\zeta\, *}$ et  conclure que $c\sim_{\partial_{i} \mb K^{\zeta\,
*}}a$ et $c\sim_{L}b$.

Soient maintenant $\Sigma_{1}\subset\{x_{0}\}\times\D(\zeta)$ et
$\Sigma_{2}\subset\D(\zeta)\times\{y_{0}\}$, $|x_{0}|=1$,
$|y_{0}|=1$, deux disques conformes sur lesquels les transformations
d'holonomie $h_{1}:\Sigma_{1}\to \{x_{0}\}\times \D(1)$ et
$h_{2}:\Sigma_{2}\to \D(1)\times \{y_{0}\}$,
$$h_{1}(x_{0},y)=(x_{0},y e^{2i\pi\lambda}), \qquad h_{2}(x,y_{0})=(x e^{2i\pi/\lambda},y_{0}),$$
sont bien d\'{e}finies. On peut distinguer trois cas :
\\
\indent \textbf{A)} $\mathbf{\lambda\notin\mathbb{Q}.}$   Si
$h_{i}(\Sigma_{i})\subset\Sigma_{i}$ ou si $\Sigma_{i} \subset
h_{i}(\Sigma_{i})$, alors la suspension $U_{\Sigma_{i}}$ v\'{e}rifie :
$U_{\Sigma_{i}}^*\uco{\mc G}\partial_{i} \mb K^{\zeta\,
*}$. Ceci est une cons\'{e}quence immediate de la proposition
(\ref{prop3}). D'autre part chaque feuille est contractile. Ainsi,
par transitivit\'{e} de la connexit\'{e} {feuillet\'{e}} l'ensemble
$$\mc B_\zeta :=\inte{\mb
K}{}^{\zeta\, *}\cup U_{\Sigma_{1}}^*\cup U_{\Sigma_{2}}^*$$ est un
bloc $\mc G$-adaptable.
\\
\indent \textbf{B)} $\mathbf{\lambda\in\mathbb{R}_{>0}.}$ Dans ce
cas les hypersurfaces r\'{e}elles $\{ |x|^{-\lambda}\,| y|  = \kappa
\}$, $\kappa \in \R_{>0}$, sont satur\'{e}es pour le feuilletage $\mc G$
et s\'{e}parent le polydisque en deux composantes connexes. Pour
$\kappa$ assez grand, $\Omega _\kappa  :=
\{|x|\,|y|^{-\frac{1}{\lambda}} > \kappa \}\cap \mb K$ est un
voisinage $\mc G$-satur\'{e} dans $\mb K$ de l'axe des $x$ \'{e}point\'{e} de
l'origine. De m\^{e}me pour $\tau >0$ assez grand $\Omega '_\tau :=
\{|x|^{-\lambda}\,|y| > \tau \}\cap \mb K^{\zeta}$ est un voisinage
$\mc G$-satur\'{e} dans $\mb K^{\zeta}$ de l'axe des $y$ \'{e}point\'{e}.
Visiblement ces ensembles se r\'{e}tractent par d\'{e}formation sur
$W_\kappa := \Omega _\kappa \cap \{|x|=1\}$ et $W_\tau' := \Omega
'_\tau \cap \{|y|=1\}$ respectivement. La r\'{e}traction est donn\'{e}e par
une homotopie laissant invariant chaque feuille de $\mc G$. On en
d\'{e}duit imm\'{e}diatement les relations :
$$W_\kappa^* \;\uco{\mc G}\; \Omega _\kappa ^*\; \uco{\mc G}\;
\mc B_{\zeta,\, \kappa ,\,\tau } \quad \mathrm{et} \quad W'_\tau{}^*
\;\uco{\mc G}\; \Omega '_\tau {}^*\; \uco{\mc G}\; \mc B_{\zeta,\,
\kappa ,\,\tau }\,,$$ avec : $$\mc B_{\zeta,\,\kappa ,\,\tau } :=
\inte{\mb K}{}^{\zeta\, *}\cup W_\kappa^* \cup W'_\tau{}^*\,.$$ On
en d\'{e}duit que   $\mc B_{\zeta,\, \kappa ,\,\tau }$ est un bloc $\mc
G$-adaptable. En effet il ne reste \`{a} v\'{e}rifier que
l'incompressibilit\'{e} des feuilles de $\mc G_{|\partial_i\mc
B_{\zeta,\,\kappa ,\,\tau }}$ dans $\partial_i \mc B_{\zeta,\,\kappa
,\,\tau }$, $i=1,2$. Lorsque $\lambda $ est irrationnel ceci est
trivial, car les feuilles sont contractiles. Lorsque $\lambda
 = p/q$, $p, q \in \mb N^*$ les feuilles sont les courbes d'\'{e}quation
$y^q = C x^p$, et leurs traces sur $\partial_i \mc B_{\zeta,\,\kappa
,\,\tau }$ sont les lacets $( e^{iq\theta}, c e^{ip\theta})$ ou bien
$(c e^{iq\theta},  e^{ip\theta})$, $\theta\in [0, 2\pi ]$, $c$
constante $>0$. Elles induisent dans les deux cas l'\'{e}l\'{e}ment non nul
$(q, p)$ de $\pi _1(\partial_i \mc B_{\zeta,\,\kappa ,\,\tau }) \iso
\mb Z\times \mb Z$. Remarquons aussi que  chaque composante
$W_\kappa^* $ et $ W'_\tau{}^*$ du bord
de $\mc B_{\zeta,\,\kappa ,\,\tau }$ est de type suspension.\\
\indent \textbf{C)}
$\mathbf{\lambda=-\frac{p}{q}\in\mathbb{Q}_{<0}^*.}$ Le feuilletage
poss\`{e}de l'int\'{e}grale premi\`{e}re $x^py^q$. Le sous-ensemble
$$\mc B_\zeta :=\{|x|\leq 1,\, |y|\leq 1,\,
0<|x^{p}y^{q}|< \zeta\}$$ est $1$-$\mc G$-connexe dans $\mb
K^{\zeta}$ pour tout $\zeta>0$ car il est $\mc G$-satur\'{e}. Il est
facile de construire des r\'{e}tractions par d\'{e}formation de $\mc
B_{\zeta }$ sur $\partial_1 \mc B_{\zeta}$ et sur $\partial_2 \mc
B_{\zeta }$, laissant invariant chaque feuille de $\mc G_{|\mc
B_{\zeta }}$. Comme pr\'{e}c\'{e}demment les feuilles des restrictions
$\partial_i \mc B_\zeta$ induisent l'\'{e}l\'{e}ment $(q, -p)$ de $\pi
_1(\partial_i \mc B_\zeta) \iso \mb Z\times \mb Z$, ce qui montre
leur incompressibilit\'{e} dans dans le bord de $\mc B_\zeta$ -et donc
aussi dans $\mc B_\zeta$. Ainsi $\mc B_\zeta$ est un bloc $\mc
G$-adaptable. Notons que les deux composantes du bord sont encore de
type suspension.

\subsubsection{{\bf Exemple 2 :}} Si $\lambda   \in
\mathbb{Q}_{<0}$ et $\mc G$ n'est pas lin\'{e}arisable on sait --et cela
sera pr\'{e}cis\'{e} en (\ref{retractcollier}) au paragraphe
(\ref{passgecols})-- qu'il existe encore un collier feuillet\'{e}
$\Omega $. Mais maintenant, pour tout disque conforme assez petit
$\Sigma \subset \{x= x_0\}\cap \Omega $, $U_\Sigma ^*$ n'est jamais
$0$-connexe dans $W := \{ |x| = 1\}\cap \Omega $; on en d\'{e}duit
facilement que $U_\Sigma ^*$ n'est jamais $1$-connexe dans $W$ !
Ainsi, dans ce cas, pour obtenir un bloc $\mc G$-adaptable on ne
peut pas  adjoindre \`{a} $\inte \Omega $ des bords trop petits.
Cette difficult\'{e} sera lev\'{e}e au paragraphe \ref{decompcollier}.\\

\vspace{1em}

\subsection{Notion de rugosit\'{e}}\label{subsecrugos}
Dans les exemples pr\'{e}c\'{e}dent nous avons construit des ouverts
adaptables dont les composantes  bord sont des ensembles $U_\Sigma $
de type suspension. Cette propri\'{e}t\'{e} du bord sera essentielle pour
construire par induction un  assemblage bord \`{a} bord feuillet\'{e}. Dans
tous les exemples $U_\Sigma^* $ a le type d'homotopie d'un {tore %vide
$\mb S^{1}\times \mb S^{1}$}, car $h(\Sigma )\cap \Sigma $ est
connexe. S'il n'en \'{e}tait pas ainsi, $U_{{\Sigma}}^*$ aurait le type
d'homotopie d'un tore auquel serait {attach\'{e}} un bouquet de cercles.
Alors $U_\Sigma ^*$ ne serait pas incompressible dans $\mc B$ et la
propri\'{e}t\'{e} (\ref{bordincompr}) de la d\'{e}finition (\ref{adaptable}) de
bloc adaptable ne serait pas satisfaite. La notion de rugosit\'{e} que
nous introduisons ici permet de contr\^{o}ler le caract\`{e}re {\'{e}toil\'{e}  des
domaines $\Sigma $ et $h(\Sigma )$}, ce qui garanti la connexit\'{e} de
{leur} intersection.\\

Notons, pour $\;z \in \mb{C}, \; |z| = 1$
\begin{equation}\label{doublecrocher}
[[ z ]] := \left\lbrace
\begin{array}{lllll}
|\theta|\,, &\mathrm{si} & z = |z| \,e^{i\theta}\,,& \mathrm{avec}
&\theta \in \; ]-\frac{\pi}{2}, \; +\frac{\pi}{2}[ \, ,\\
+ \infty, & \mathrm{si} & z =  |z|\, e^{i\theta}\,,& \mathrm{avec}
&\theta \in \; [\frac{\pi}{2}, \; \frac{3\pi}{2}].
\end{array}
 \right.
\end{equation}
\noindent Tout d'abord consid\'{e}rons un chemin $ \gamma: [0, 1]
\longrightarrow \mb C^{\ast}$ analytique et \textit{lisse}, i.e.
$\gamma'(s) \neq 0$. Pour tout $s\in [0, 1]$ nous notons :
$${\mathbf{e}}\,(\gamma; s) := \left[\left[\,\frac{\gamma'(s)}
{i \, \gamma(s)} \, \right]\right].$$ \noindent Cette quantit\'{e}
``mesure''  l'angle de la courbe orient\'{e}e $\mathrm{Im}(\gamma ):=
\gamma ([0, 1])$ avec la tangente au cercle de
centre $0$ passant par $\gamma (s)$.\\

\indent Appelons ici  \textbf{chemin analytique lisse par morceaux}
(\textbf{a. l. p. m.} en abr\'{e}g\'{e}) tout chemin  $\mu = \mu_{1} \vee
\cdots \vee \mu_{p}$ qui est  une concat\'{e}nation  de chemins
analytiques lisses  $\mu_{j}$ tels que  chaque intersection
$\mathrm{Im }(\mu _i)\cap \mathrm{Im} (\mu _j)$, $1\leq i < j \leq
p$, est : ou bien un ensemble discret, ou bien une courbe simple sur
laquelle les orientations de $\mu _i$ et de $\mu _j$ co\"{\i}ncident.
\begin{defin}\label{defin31}
Nous appelons {\bf rugosit\'{e} de $\,\mu\,$ au point $\mu (s)$} et
notons encore $\be (\mu ; \, s)$ le maximum des limites \`{a} gauche et
\`{a} droite en $s$ de ${\be}\,(\mu;\, t)$. Nous appelons {\bf rugosit\'{e}
de $\,\mu\,$} l'\'{e}l\'{e}ment  $\be( \mu ) := \max \{\be (\mu ;\,s)\;/\;
s\in [0,1] \}$ de $\overline{\R}_+{=\R_{+}\cup\{+\infty\}}$.
\end{defin}

\begin{obs}\label{depimage}Remarquons que $\be(\mu )$ ne d\'{e}pend que de la  courbe orient\'{e}e
$\mathrm{Im} (\mu )$. Notons aussi que si $\,\mathbf{e}(\gamma; s)\;
< \;+\infty\,$, alors $\, \mathbf{e}(\gamma^{-1}; s) =
+\infty\,$.\end{obs}

{}%\indent
Il est clair qu'un domaine simplement connexe $\Delta \ni 0
\,$ de bord une courbe a. l. p. m. est \'{e}toil\'{e} par rapport \`{a}
l'origine d\`{e}s que $\mathbf{e}(\partial\Delta)\,$ est fini. D'autre
part, pour deux domaines $\Delta$, $\Delta'$ contenant l'origine on
a :
\begin{equation}\label{rugsom}
\be ( \partial (\Delta \cup \Delta'))\,, \quad \be ( \partial
(\Delta \cap \Delta'))\quad \leq \quad\max \left\{ \be (\partial
\Delta) , \be (\partial \Delta')\right\} \,. \end{equation}

\indent Plus g\'{e}n\'{e}ralement supposons  $\mu \,$ trac\'{e}e dans une
vari\'{e}t\'{e} analytique r\'{e}elle $M$, non-n\'{e}cessairement connexe. Soit
$\xi$ une fonction $\mathbb{R}$-analytique, d\'{e}finie sur un voisinage
ouvert de $\mathrm{Im} (\mu)$ et \`{a} {valeurs} dans $\mathbb{C}$ telle
que  $\xi \circ \mu$ n'est pas constant. Visiblement la courbe
$\xi\circ \mu $ est aussi a.l.p.m. Nous appelons
\textbf{$\xi$-rugosit\'{e} de $\mu$} l'\'{e}l\'{e}ment $${\mathbf{e}}_\xi\,(\mu)
:= {\mathbf{e}}\,(\xi \circ \mu)\in \overline{\R}_{+}\,.$$ \noindent
Nous adoptons dans tout le texte les notations suivantes :
\begin{equation}\label{diametre}
\|\mu\|_{\xi} := \max_{s\in[0,1]}\,|\xi \circ \mu(s)|\quad
\{\!\!\{\,  \theta \, \}\!\!\} := \left\{\begin{array}{ccc}
{|\theta|}& si & \theta\in ]-\pi /2 , \, + \pi /2[ \\
+\infty & si & \theta\notin ]-\pi /2 , \, + \pi /2[ \\
\end{array}\,.\right.
\end{equation}
\noindent Maintenant supposons que $\mu $ est trac\'{e} sur un
\textit{secteur ferm\'{e}}
$$\overline{S}_{R, \alpha, \beta} := \{r\, e^{i\theta} \, / \, 0 \leq r < R, \,
\alpha \leq \theta \leq \beta \,\} \, \subset \, \widetilde{\mb
C}\,,\qquad 0\leq \beta -\alpha \leq 2\pi\ $$ de la surface de
Riemann $\wt \C$ de $\C^*$. Consid\'{e}rons une application holomorphe
$g$ d\'{e}finie sur un secteur ouvert $S_{R', \alpha', \beta'}$
contenant $ \overline{S}_{R, \alpha, \beta} \,$ et qui poss\`{e}de sur
ce secteur un d\'{e}veloppement asymptotique s'\'{e}crivant $g(z) := \sum_{j
= \nu}^{\infty}c_{j}\,z^{j/p}$, avec $\nu\geq 1$, $p\in \N^*$ cf.
\cite{wasow}.

\begin{prop}\label{prop31}
Il existe une  constante $C_g> 0$ d\'{e}pendant continument de $g$ telle
que:
$${\mathbf{e}}_{\,g }\,(\mu) \, \leq \,
\{\!\!\{\,{\mathbf{e}}_{\,z}(\mu) + C_g\| \mu\|_{z}\,
\}\!\!\}\,,\quad \|\mu \|_g \leq C_g\|\mu \|_z$$
\end{prop}

\begin{dem} La seconde in\'{e}galit\'{e} est facile. Pour prouver
la premi\`{e}re in\'{e}galit\'{e}, on peut supposer  $\mu \,$ lisse et
$c_\nu\neq 0$. On a :
\begin{equation}\label{com.rugos}
    \frac{(g  \, \circ \, \mu)'}{i\cdot
g \, \circ \, \mu } = \psi \,\circ \mu \, \cdot \, \frac{ \mu'}{i
\mu}, \qquad \mathrm{avec} \quad \psi(z) := \frac{z \, g'(z)}{g(z)}.
\end{equation}
\noindent Classiquement \cite{wasow} lorsque $z$ tend vers $0$ sur
$\,S_{R, \alpha, \beta} \,$ et la fonction $\psi(z)$ poss\`{e}de un
d\'{e}velop\-pement assymptotique du type {$\nu/p + o(1)$}. Ainsi
$\psi(z)$ tend vers $\nu/p \in {\mathbb{Q}_{>0}}$ et
$|\arg(\psi(z)/z)| \leq C' |z| $ pour une constante $C'>0$
appropri\'{e}e et pour $|z|$ assez petit, disons $|z|\leq C_1 $. Toute
constant $C_g>0$ telle que $C_g\geq C'$ et $C_g C_1 >\pi /2$
convient. On en d\'{e}duit sans peine la majoration d\'{e}sir\'{e}e. Nous
laissons au lecteur le soin de v\'{e}rifier que si $g$ d\'{e}pend
continument d'un param\`{e}tre, on peut choisir la constante $C_g$
d\'{e}pendant continument de ce param\`{e}tre.
\end{dem}
\noindent Quitte \`{a} restreindre $R$ l'application $g$ poss\`{e}de une
``r\'{e}ciproque'' $g^{-1}$ d\'{e}finie sur un secteur $S_{R'', \alpha'',
\beta''}\subset \wt \C$ et \`{a} valeur sur un ouvert de $S_{R',
\alpha', \beta'}$ contenant $\overline{S}_{R, \alpha, \beta}$. On en
d\'{e}duit le corollaire suivant.
\begin{cor}\label{cor.est.rug.}
Il existe une  constante $C_g'> 0$ d\'{e}pendant continument de $g$
telle que:
$${\mathbf{e}}_{\,z }\,(\mu) \, \leq \,
\{\!\!\{\,{\mathbf{e}}_{\,g}(\mu) + C_g'\| \mu \|_{g}\, \}\!\!\}\,,
\quad \|\mu \|_z \leq C'_g \| \mu \|_g.$$
\end{cor}

\vspace{1em}

{\section{R\'{e}duction de la preuve du th\'{e}or\`{e}me principal}}

\subsection{D\'{e}composition du diviseur et construction de bonnes
fibrations}\label{decompdiv} Consi\-d\'{e}\-rons le diviseur total:
$$\mc D :=E^{-1}(S)= \mc E \cup \mc S \,,\quad
\mc E :=E^{-1}(0)\,,\quad \mc S := \mc
 S_1 \cup \cdots \cup \mc S_{\varrho} \,.$$
\noindent Nous reprenons le vocabulaire habituellement utilis\'{e} : la
\textbf{valence} $\upsilon(D)$ d'une composante irr\'{e}ductible) $D$ de
$\mc E$ est le nombre de singularit\'{e} de ${\F}$ situ\'{e}es sur $D$, une
\textbf{cha\^{\i}ne} est une union connexe maximale de composantes de
$\mc E$ de valence $\leq 2$, une \textbf{branche morte} est une
cha\^{\i}ne de composantes poss\'{e}dant une composante de valence 1, appel\'{e}e
\textbf{composante d'extr\'{e}mit\'{e}}. Deux unions connexes de composantes
de $\mc D$ sont dites \textbf{adjacentes} si leur intersection est
non-vide et r\'{e}duite \`{a} un point, appel\'{e}e \textbf{point d'attache}.
Nous appellerons \textbf{composante simple de $\mc E$} toute
composante $D$ de $\mc E$ qui n'intersecte aucune branche morte, et
\textbf{bloc agr\'{e}g\'{e} de $\mc E$} toute union connexe maximale de
composante non-simples de $\mc E$. Un bloc agr\'{e}g\'{e} est toujours
constitu\'{e} d'une composante de $\mc E$ appel\'{e}e \textbf{composante
centrale du bloc} et de l'union (non-vide) des
branche mortes adjacentes \`{a} cette composante. \\

Visiblement les transform\'{e}s strictes $\mc S_j$, les composantes
simples et les blocs agr\'{e}g\'{e}s forment un recouvrement de $\mc D$ en
sous-ensembles qui, deux \`{a} deux, sont ou bien adjacent, ou bien
d'intersection vide.\\

En chaque point singulier $s\in \mc D $ de $\wF$ nous nous fixons
maintenant des coordonn\'{e}es $(x_s, y_s) : K(s)\, \iso\, \D(2)\times
\D(2)$ telles que $\mc D\cap K(s) \subset \{x_s y_s =0\}$ et que :
\begin{itemize}
  \item si le germe $\wF_s$ de $\wF$ en $s$ est
lin\'{e}arisable, alors  il est est d\'{e}finit par un champ de vecteur
lin\'{e}aire dans ces coordonn\'{e}es $X_s =
x_s\frac{\partial\phantom{x_s}}{\partial x_s} + \lambda
y_s\frac{\partial\phantom{y_s}}{\partial y_s}$
  \item si $\wF_s$ est r\'{e}sonnant non-lin\'{e}arisable, alors il est  d\'{e}fini par
  une $1$-forme diff\'{e}ren\-tielle du type forme normale de Dulac:
  $x_s dy_s - y_s (\lambda +x_sy_s(\cdots))dx_s$, $\lambda \in \mb Q_{<0}$,
  cf. \cite{MM}.
\end{itemize}
 Si $D$ est une composante irr\'{e}ductible de $\mc D$ qui
contient $s$, alors nous notons
$$D_s := D\cap \left\{{|x_s|<
1,\,|y_s|<1}\right\}\,.$$ \noindent Visiblement $D_s$ un disque
conforme ouvert. Quitte \`{a} composer les coordonn\'{e}es construites par
des homoth\'{e}ties appropri\'{e}es, nous supposons que deux disques ferm\'{e}s
$\overline{D}_s$ et $\overline{D}_{s'}$ contenus dans une m\^{e}me
composante $D$ ne
s'intersectent jamais.\\

Consid\'{e}rons maintenant une composante $D$ de $\mc D$ et d\'{e}signons
par
 $s_j$, $j=1 , \ldots , \upsilon(D)$ les points singuliers de $\wF$ situ\'{e}s
 sur $D$. Supposons pour simplifier l'\'{e}criture  qu'en chaque point $y_{s_j}=0$ est l'\'{e}quation locale de $D$ en
 $s_j$.
%Si $D$ n'est pas une composante simple mais poss\`{e}de $\upsilon(D)-k$  branches mortes adjacentes, nous choisissons toujours la num\'{e}rotation des points singuliers pour que $s_{k+1},\ldots,s_{\upsilon(D)}$ soient les points d'attache de ces branches mortes.
{Nous noterons}
\begin{equation}\label{disquetroue}
    D^\ast \; := \; D\; \setminus\; \bigcup_{j=1}^{\upsilon(D)}
    \overline{D }_{s_j}\;.
\end{equation}
On sait qu'il existe   une
 fibration holomorphe en disques
 $\underline{{\pi}}_D$ localement triviale, d\'{e}finie sur un voisinage ouvert
 $\underline{\Omega} _D$ de
 $D$ et \`{a} valeur dans $D$, qui est l'identit\'{e} en restriction \`{a} $D$.
Pour chaque $j=1,\ldots ,\upsilon(D)$ on construit sans difficult\'{e}
une submersion surjective $\underline{{\pi}} _{D, s_j}$ de classe
$C^\infty$ d\'{e}finie sur $ K(s_j)$ et \`{a} valeurs sur $
 D\cap\left\{|x_{s_j}|<2\right\}$, qui est \'{e}gale \`{a}
${\underline{\pi}_D}$ en restriction \`{a}  $\underline{\Omega }_D\cap
K(s_j) \cap \{ |x_{s_j}|> \frac{3}{2}\}$ et \'{e}gale \`{a} la projection
$x_{s_j}$ en restriction \`{a} $K(s_j)\cap \{|x_{s_j}|<1\}$. Par des
recollements appropri\'{e}s on obtient une fibration en disques
\begin{equation}\label{fibr}
{\pi}_D : \Omega (D) \longrightarrow D
\end{equation}
de classe $C^\infty$ d\'{e}finie sur un voisinage ouvert $\Omega (D)$ de
$D$, \'{e}gale \`{a} l'identit\'{e} en restriction \`{a} $D$ et telle que pour
chaque $j=1,\ldots \upsilon(D)$ on a :
\begin{enumerate}
  \item  ${\pi} _D^{-1}(D_{s_j})
  =K(s_j)$; et en restriction \`{a} cet ouvert ${\pi} _D$ est \'{e}gal
  \`{a} la projection $x_{s_j}$,
  \item $
{\pi} _D^{-1}(\wt D^*)=
  \underline{{\pi}} _D^{-1}(\wt D^*)$ avec $
  \wt D^* := D\setminus \bigcup_{j=1}^{\upsilon(D)}
\{|x_{s_j}|\leq 2\}$; et en restriction
  \`{a} cet ouvert ${\pi} _D$ est \'{e}gal \`{a}
  $\underline{{\pi}}_D$.
\end{enumerate}

Maintenant appelons \textbf{bloc \'{e}l\'{e}mentaire de $\mc D$} les
adh\'{e}rences des composantes connexes du compl\'{e}mentaire dans $\mc D$
de l'union $\bigcup_{(D,\,s)\in \mf A}\partial D_s$ o\`{u}
$$\mathfrak{A}:= \big\{ (D, s) \in Comp(\mc D)\times Sing(\wF)\;|
 \; s\in D\big\}$$
et $Comp(\mc D)$ d\'{e}signant l'ensemble des composantes irr\'{e}ductibles
de $\mc D$. Notons
\begin{equation}\label{trans.equa.sep}
    F := f \circ E
\end{equation}
le compos\'{e} de l'\'{e}quation $f$ de la s\'{e}paratrice totale par
l'application de r\'{e}duction. Pour $\eta>0$ assez petit l'hypersurface
r\'{e}elle $$\mc H := \bigcup_{(D, s)\in \mf A} {\pi} _D^{-1}(\partial
D_s)$$ est transverse aux fibres $F^{-1}(z)$, $|z|\leq \eta$.
Consid\'{e}rons la famille $(\underline{{\mc T}}_\eta^j)_{j\in \mf J}$
des composantes connexes de $T_\eta \setminus \mc H $ o\`{u}, comme dans
l'introduction, ${\mc T}_\eta:= E^{-1}({T}_\eta)$. Chaque bloc
\'{e}l\'{e}mentaire de $\mc D$ est contenu dans l'adh\'{e}rence d'une et une
seule composante $\underline{{\mc T}}^j_\eta$.

\begin{defin}\label{bloc.mil}
Soit $A\subset \mc D$ nous appelons ici \textbf{bloc de Milnor de
hauteur $\eta>0$ associ\'{e} \`{a} $A$} l'union des adh\'{e}rences dans ${\mc
T}_\eta$ des {composantes}  $\underline{{\mc T}}_\eta^j$ de {$\mc
T_{\eta}\setminus\mc H$} qui intersectent $A$ :
$${\mc T}_\eta(A) :=
\bigcup_{j \in \mf J_A}\overline{\underline{{\mc T}}_\eta^j}\,,
\qquad \mf J_A := \{i\in \mf J \; | \; \underline{{\mc T}}_\eta^j
\cap A \neq \emptyset\}\,.$$ Nous notons aussi ${\mc
T}_\eta^*(A):={\mc T}_\eta(A)\setminus \mc D$.
\end{defin}
\noindent Le r\'{e}sultat suivant bien connu de sp\'{e}cialistes se montre
facilement \`{a} l'aide de la remarque (\ref{rVK}).
\begin{obs}\label{iso.bloc.miln} Il existe $\eta_1>0$ tel que pour tout
$A\subset \mc D$ et $0<\eta' \leq \eta\leq \eta_1$, les inclusions
${\mc T}_{\eta'}^*(A)\subset {\mc T}_\eta^*(A)$ et $\partial {\mc
T}_{\eta'}^*(A)\subset
\partial {\mc T}_\eta^*(A)$ induisent des isomorphismes au niveau du
groupe fondamental.
\end{obs}

\vspace{1em}

\subsection[*]{{Blocs feuillet\'{e}s adaptables \`{a} taille et rugosit\'{e} contr\^{o}l\'{e}es}}\label{preuveprinc}

%\subsection{Preuve du Th\'{e}or\`{e}me Principal}\label{preuveprinc}

Soit $\mc A$ l'ensemble dont les \'{e}l\'{e}ments sont :
\begin{enumerate}[(i)]
  \item\label{ptsing} les points singuliers de $\mc D$,
  \item \label{blocscomp} les composantes simples de $\mc E$,
  les composantes centrales des blocs agr\'{e}g\'{e}s de $\mc E$, les
transform\'{e}es strictes des s\'{e}paratrices.
\end{enumerate}
Pour chaque $\alpha \in \mc A$ consid\'{e}rons le sous-ensemble
$K_\alpha $ de $\mc D$ suivant

\begin{itemize}
  %\item { dans le  cas (i) avec ${\alpha  = s \notin Sing(\mc D)}$},  on pose $K_\alpha  := \overline{D_s}\,$, o\`{u} $D$ est la composante de $\mc D$ qui contient $s$,
  \item{ dans le cas (i) avec $\alpha  = s \in Sing(\mc D)$},  on pose $K_\alpha  :=
  \overline{D_s}\cup \overline{D'_s}\,$,
 o\`{u} $D$ et $D'$ sont les deux composantes de $\mc D$ qui
 s'intersectent
 en $s$,
  \item dans le cas (ii) : $K_{\alpha}=\overline{D}^{\sharp}\cup\bigcup_{j\in\mc J_{D}}\mf M_{j}$
o\`{u} $D=\alpha$ est la composante de $\mc D$ consid\'{e}r\'{e}, $\{\mf
M_{j}\}_{j\in\mc J_{D}}$ est la collection (peut-\^{e}tre vide) des
branches mortes adjacentes \`{a} $D$, $\{s_{j}\}=\mf M_{j}\cap D$ et
$D^{\sharp}= \overline{D}^*\cup\bigcup_{j\in \mc J_{D}}D_{s_{j}}$.
\end{itemize}
\noindent L'intersection de deux sous-ensembles distincts $K_\alpha
$ et $K_{\alpha '}$ est ou bien vide, ou bien diff\'{e}omorphe \`{a} un
cercle et celui-ci est une composante connexe commune de $\partial
K_{\alpha
'}$ et de $\partial K_\alpha$.\\

Pour chaque composante connexe $\mc C$ du bord d'un $K_\alpha $,
$\alpha \in \mc  A$, donnons-nous un champ de vecteur r\'{e}el r\'{e}gulier
$Z_{\mc C}$ qui d\'{e}finit le feuilletage $\wF_{|{\pi}_D^{-1}(\mc C)}$,
$D$ d\'{e}signant la composante de $\mc D$ qui contient $\mc C$. En
chaque point $P \in \mc C$ est d\'{e}finit un germe de fonction ``temps
de premier retour'' $ \tau _P : {\pi}_D^{-1}(P) \rightarrow \R_{>0}$
qui v\'{e}rifie, pour tout $m\in {\pi}_D^{-1}(P)$ suffisamment proche de
$P$ :
$$\Phi_t^{Z_{\mc C}}(m) \in
  {\pi} _D^{-1}(\mc C\setminus \{P\})\quad\mathrm{pour}\quad
  0< t < \tau (m)\quad \mathrm{et}\quad
\Phi_{\tau (m)}(m)\in {\pi} _D^{-1}(P)\,,$$ $\Phi_t^{Z_{\mc C}}$
d\'{e}signant le flot de $Z_{\mc C}$. Nous dirons qu'un ensemble
$V\subset {\pi} _D^{-1}(C)$ est de type \textbf{suspension  au
dessus de $\mc C$}, s'il existe un point $P\in \mc C$ et un disque
conforme ouvert $\underline{\Sigma} \subset {\pi} _D^{-1}(P)$
contenant $P$, tel qu'en notant $\Sigma := \underline{\Sigma
}\setminus \{P\}$ on a :
$$V= V_\Sigma  := \left\{ \Phi_t^{Z_{\mc C}}(m) \;/\;  m\in \Sigma
,\; 0\leq t\leq \tau _P(m) \right\}\,.$$ Appelons ici
\textbf{rugosit\'{e} de $V$ } et notons $\be _F(V)$, la $F$-rugosit\'{e} du
bord de $\overline{\Sigma }$. De m\^{e}me {nous} appelons \textbf{taille
de $U$} le r\'{e}el $\| U  \|_F := \| \Sigma  \|_F$.

La proposition suivante est vraisemblablement bien connue des
{sp\'{e}cia\-lis\-tes}.
\begin{prop}\label{sommetexecpt}
Si le transform\'{e} strict d'un germe de feuilletage ne poss\`{e}de aucun
point singulier de type selle-n{\oe}ud, alors il existe au plus une
composante du diviseur exceptionnel adjacente \`{a} au moins deux
branches mortes.
\end{prop}
Lorsqu'elle existe, nous appellerons cette composante la
\textbf{composante centrale} du diviseur exceptionnel. Remarquons
que toute autre composante, ou bien contient (au moins) deux points
singuliers de $\F$ qui ne sont pas des points d'attache de branches
mortes, ou bien est contenue dans une branche morte. Remarquons
aussi qu'une composante centrale de $\mc E$ est le centre d'un bloc
agr\'{e}g\'{e}. C'est donc un \'{e}l\'{e}ment de $\mc A$.
\begin{dem2}{de la proposition}
D'apr\`{e}s l'hypoth\`{e}se l'application $E$ de r\'{e}\-duction de $\F$ est
identique \`{a} l'application de r\'{e}duction de la s\'{e}paratrice totale
\cite{Cam.Lins.Sad}. Il suffit d\'{e}montrer la propri\'{e}t\'{e} suivante des
germes de courbes  $X\subset (\mb C^2, 0)$.
\begin{enumerate}[(a)]
  \item[] \it Le diviseur exceptionnel $\mc D_X$ de la
  r\'{e}duction de $X$ poss\`{e}de
  au plus une composante
  centrale et dans ce cas le
  diviseur cr\'{e}\'{e} par le premier \'{e}clate\-ment est la composante d'extr\'{e}mit\'{e} d'une
  branche morte.
\end{enumerate}
Raisonnons par r\'{e}curence sur le nombre  $N_X$ d'\'{e}clatement
n\'{e}cessaires \`{a} la r\'{e}duction de $X$. L'assertion est triviale lorsque
$N_X = 1$. D'autre part l'application de r\'{e}duction de $X$ est le
compos\'{e} de l'application d'\'{e}clatement de l'origine $\wt E :
\wt{\C^2}\rightarrow \C^2$ avec l'application $E_{X'} : M
\rightarrow \wt{\C^2}$ de r\'{e}duction de $X' := \wt E^{-1}(X)$. En
chaque point singulier $c_j$, $j=1,\ldots ,\varrho_X$ de $X'$ nous
pouvons appliquer l'hypoth\`{e}se de r\'{e}curence. Le diviseur $\mc D_X$
est l'union des diviseurs exceptionnels $\mc D_j = E_{X'}^{-1}(c_j)$
de la r\'{e}duction des germes de $X'$ en $c_j$ et du transform\'{e} strict
$D_0 $ de $ \wt E^{-1}(0)$ par $E_{X'}$. Chaque $\mc D_j$ est
adjacent \`{a} $D_0$. Il est \'{e}vident que $\mc D_X$ ne poss\`{e}de pas de
composante centrale si aucun des $\mc D_{X',\, c_j}$ n'en poss\`{e}de,
ou bien si $\varrho_X \geq 2$. Consid\'{e}rons le cas  :  $\varrho_X =
1$ et $\mc D_j $ poss\`{e}de deux {branches} mortes. On voit facilement
que, $ \wt E^{-1}(0)$ \'{e}tant lisse, $D_0 $ est n\'{e}cessairement
adjacent \`{a} une composante $D_1$ contenue dans l'une des branches
mortes $\mf M$. On consid\`{e}re alors s\'{e}par\'{e}ment les deux \'{e}ventualit\'{e}s
: $D_1$ est une composante d'extr\'{e}mit\'{e} de $\mf M$, ou non. Dans
chaque cas l'hypoth\`{e}se de r\'{e}curence est imm\'{e}diatement satisfaite.
\end{dem2}

D\'{e}signons toujours par $\eta_1$ la ``hauteur d'uniformit\'{e}'' des
blocs de Milnor donn\'{e}e par la remarque (\ref{iso.bloc.miln}).

\begin{teo}[\bf d'existence de blocs adaptables]\label{teoexistence}
Soit $\alpha \in \mc A$ et  $\varepsilon>0$. Il existe $\mc
B_\alpha\subset T_{\eta_1}^*(K_\alpha )$ tel que :
\begin{enumerate}
\item pour $ \eta >0$ assez petit, $\mc B_\alpha $ contient
$T_{\eta}^*(K_\alpha )$ et les applications d'inclusion induisent
des isomorphismes $$\pi _1(T_{\eta}^*(K_\alpha ),
\mathbf{\cdot})\iso \pi _1(\mc B_\alpha, \mathbf{\cdot})\quad
et\quad \pi _1(\partial T_{\eta}^*(K_\alpha ) ,\, p)\iso \pi
_1(\partial \mc B_\alpha  ,\, p)\,,$$ pour tout $p\in
\partial B_{\varepsilon_0,\,\eta}^\alpha$,
\item $\mc B_\alpha$
  est un bloc $\F$-adaptable,
\item les composantes connexes  $V_1,\ldots ,V_{r_\alpha }$ de $\partial\mc B_\alpha $
  sont des ensembles  de type suspension au dessus des
  composantes connexes
  de $\partial K_\alpha $,
  \item $\be_F(V_j)+ \|V_j\|_F
  \leq \varepsilon\,$, $j=1,\ldots ,r_\alpha $.
\end{enumerate}
Si nous supposons de plus que $\alpha $ n'est pas la composante
centrale de $\mc D$, alors il existe des constantes $C_{\alpha
,\,0}$, $C_{\alpha}$, $\xi_\alpha >0$ et une fonction $\mf d_\alpha
: \R_+\rightarrow \R_+$, $ \lim_{r\rightarrow 0} \mf d_\alpha (r) =
0$, telles que %:
%\begin{itemize}
  %\item
{pour} tout sous-ensemble $V$ de type suspension au dessus d'une
composante de $\partial K_\alpha $ v\'{e}rifiant $\be_F(V)\leq
{C_{\alpha,\,0}}$ et $\|V\|_F \leq {C_{\alpha,\,0}}$, il existe
$\mc B_\alpha\subset T_{\eta_1}^*(K_\alpha )$ qui satisfait les
propri\'{e}t\'{e}s (1), (2) et (3) pr\'{e}c\'{e}dentes, ainsi que :
\begin{enumerate}
  \item [(3')]  $V_1 \ucon V$,
  \item [(4')] $\be_F (V_j)\leq \{\!\!\{\be_F(V) +
\mf d_\alpha  (\| V\|_F )\}\!\!\}\,$ et $\| {V}_j\|_F \leq
  C_{\alpha} \| V\|_F^{\xi_\alpha }$, $j=1,\ldots ,r_\alpha $.
\end{enumerate}
%\end{itemize}
\end{teo}

Fixons maintenant un \'{e}l\'{e}ment $\alpha _0$ de $\mc A$, que nous
choisirons \'{e}gal \`{a} la composante centrale de $\mc D$ si celle-ci
existe. Consid\'{e}rons la filtration croissante du diviseur
$$
K_{\alpha _0}=: \mc K_0 \subset \mc K_1 \subset \cdots  \cdots
\subset \mc K_\kappa = \mc D\,,\qquad \mc K_{j+1} := \mc K_j \cup
\bigcup_{\alpha \in \mc A_j}K_\alpha\,,
$$
o\`{u} $\mc A_j$ est l'ensemble des $\alpha \in \mc A$ tels que
$K_\alpha \cap \mc K_j\neq \emptyset$. \\

Le th\'{e}or\`{e}me ci-dessous permet de construire par induction des
ensembles $\mc B_\alpha \subset {\mc T}_{\eta_1}^*$, $\alpha \in \mc
A$ de mani\`{e}re que l'union $U' := \bigcup_{\alpha \in \mc A}{\mc
B}_\alpha $ est un assemblage bord \`{a} bord feuillet\'{e} bien construit,
son adh\'{e}rence est un voisinage de $\mc D$ contenant ${\mc T}_{\eta}$
pour $\eta$ assez petit et $\overline{U'}\setminus \mc D = U'$.
L'image $U := E(\overline{U'})$ satisfait les conclusion du th\'{e}or\`{e}me
principal. En effet les assertions (2) et (3) sont une cons\'{e}quence
imm\'{e}diate du th\'{e}or\`{e}me (\ref{teoexistence}) ci-dessus, du th\'{e}or\`{e}me de
localisation (\ref{teo.localisation}) et de la remarque
(\ref{reassemblage}). Les assertion (1) et (2) du th\'{e}or\`{e}me
d'existence  montre que les composantes connexes du bord de chaque
${\mc T}_{\eta}^*(K_\alpha )$, sont incompressibles dans ${\mc
T}_{\eta}^*(K_\alpha )$. Ainsi ${\mc T}_{\eta}^*$ est un assemblage
bord \`{a} bord bien construit des ${\mc T}_{\eta}^*(K_\alpha )$. La
proposition (\ref{assemblages}) donne l'isomorphisme $\pi _1( {\mc
T}_{\eta}^*) \iso \pi _1( U')$  induit par l'inclusion. On en d\'{e}duit
l'isomorphisme  $\pi _1({T}_\eta\setminus S) \iso \pi _1( U\setminus
S)$ et l'assertion (1) du th\'{e}or\`{e}me principal de ``la propri\'{e}t\'{e}
d'uniformit\'{e}'' (\ref{iso.bloc.miln}).

\begin{obs}\label{pi1U}
En appliquant de fa\c{c}on r\'{e}currente le th\'{e}or\`{e}me de Seifert-Van Kampen
aux blocs $\{{\mc T}_{\eta}^*(K_{\alpha})\}_{\alpha\in\mc A}$ et la
pr\'{e}sentation (\ref{pi1B}) du groupe fondamental de $\pi_{1}({\mc
T}_{\eta}^{*}(K_{D}))$ qu'on d\'{e}montrera \`{a} la section
\ref{bloc-seifert}, on obtient une pr\'{e}sentation explicite de
$\pi_{1}(U\setminus S)$ ayant comme syst\`{e}me de g\'{e}n\'{e}rateurs un
ensemble de lacets $\{a_{D}\}$, index\'{e} par les composantes
irr\'{e}ductibles $D\subset\mathcal{D}$, v\'{e}rifiant
$\frac{1}{2i\pi}\int_{a_{D}}\frac{dF}{F}=\mathrm{ord}_{D}F$ et ayant
pour  relations
$$\prod_{D\subset \mc D}a_{D}^{(D,E)}=1,\qquad [a_{D},a_{E}]^{(D,E)}=1,
\qquad E\subset\mathcal{E}, \quad D\subset\mathcal{D}.$$
\end{obs}

Le reste de l'article est maintenant consacr\'{e} \`{a} la preuve de ce
th\'{e}or\`{e}me. Plus pr\'{e}cis\'{e}ment le chapitre \ref{bloc-seifert} est
consacr\'{e} au cas  o\`{u} $\alpha $ est une composante $D$ de $\mc D$. Le
bloc $\mc B_\alpha $ sera appel\'{e} \textbf{bloc de type Seifert}, car
il est munit par construction d'une (pseudo)-fibration de Seifert
$\mc B_\alpha \rightarrow D^{{\sharp}}$. L'existence de cette
fibration joue un role cl\'{e} dans la preuve de la $1$-$\wF$-connexit\'{e}
du bord de $\mc B_\alpha $. Le cas  o\`{u} $\alpha $ est un point
singulier de $\wF$ est trait\'{e} au chapitre \ref{decompcollier}.
Lorsque la singularit\'{e}  est lin\'{e}arisable, il suffira de reprendre
les exemples trait\'{e} en (\ref{subsect.exemples}). Il reste uniquement
le cas o\`{u} la singularit\'{e} est une selle r\'{e}sonante non-lin\'{e}arisable.
On construit alors un voisinage-collier et on introduit une
technique ``de rabotage'' qui permet de lever la difficult\'{e} \'{e}voqu\'{e}e
dans l'exemple 2 de (\ref{subsect.exemples}).

\vspace{2em}

\section{Blocs feuillet\'{e}s adaptables de type Seifert}
\label{bloc-seifert} Dans {ce chapitre} nous d\'{e}montrons le th\'{e}or\`{e}me
d'existence de blocs adaptables (\ref{teoexistence}) dans le cas o\`{u}
$\alpha\in\mc A$ est une composante irr\'{e}ductible $D$ du diviseur
$\mc D$ non-contenue dans une branche morte. Dans {le paragraphe
(\ref{3.1}) nous construisons le bloc $\mc B_{D}$ et nous prouvons
les propri\'{e}t\'{e}s (1),(3),(3'),(4) et (4') du th\'{e}or\`{e}me
(\ref{teoexistence}). Il restera seulement \`{a} prouver la propri\'{e}t\'{e}
(2) de $\F$-adaptabilit\'{e} du bloc $\mc B_{D}$. {Nous montrons}
comment la $1$-$\F$-connexit\'{e} du bord de $\mc B_{D}$ implique (2).
Avant de prouver cette propri\'{e}t\'{e} nous introduisons pr\'{e}alablement,
dans le paragraphe (\ref{3.2}), les constructions et les
descriptions auxiliaires n\'{e}c\'{e}ssaires. Finalement, dans (\ref{3.3})
nous prouvons la $1$-$\F$-connexit\'{e} du bord de $\mc B_{D}$, ce qui
ach\`{e}ve la d\'{e}monstration.}

\subsection{Construction  du bloc et r\'{e}duction de la preuve de (\ref{teoexistence})}\label{3.1}
{\`{A} partir d'un disque ouvert $\Sigma$ contenu dans la fibre de
${\pi}_{D}$ au dessus d'un point  fix\'{e} $P_{0}$ de $D^*$, nous allons
construire dans la premi\`{e}re \'{e}tape un voisinage ferm\'{e} $\overline{\mc
B}_{\Sigma}(D^*)$ de $D^*$ dans ${\pi}_{D}^{-1}(D^*)$. A la fin de
cette \'{e}tape (\ref{obs1}) nous indiquerons comment un bon choix de
$\Sigma$ et l'existence de ``bons'' voisinages  $ \overline{\mc
B}_{\Sigma}(\mf M_{j})$ des branches mortes $\mf M_{j}$ adjacentes \`{a}
$D$ permettent de d\'{e}finir, par une construction bord \`{a} bord
ad\'{e}quate, un bloc $\mc B_D$ qui satisfait le th\'{e}or\`{e}me
(\ref{teoexistence}). A l'\'{e}tape $2$ nous construirons les voisinages
$ \overline{\mc B}_{\Sigma}(\mf M_{j})$, puis \`{a} l'\'{e}tape $3$ nous
r\'{e}duirons la preuve du th\'{e}or\`{e}me (\ref{teoexistence}) \`{a} une forme
faible de
l'assertion (2) de son \'{e}nonc\'{e}. Enfin celle-ci est prouv\'{e}e \`{a} l'\'{e}tape $4$.}\\

\noindent{\textbf{Etape 1 : construction de $\mc B_{\Sigma}(D^*)$.}}
{Choisissons} la num\'{e}rotation des points singuliers $s_0,\ldots
,s_n$ de {$\wF$} sur $D$, avec $n:= \upsilon(D)-1$, pour que
$s_{k+1},\ldots,s_{n}$ ($0\le k\le n$) soient les points
d'attache des branches mortes $\mf M_{k+1},\ldots,\mf M_{n}$ adjacentes \`{a} $D$.
Consid\'{e}rons les sous-ensembles de $D$ suivants :
\begin{equation}\label{defDstarDbecarre}
    \ssD :=D\setminus\bigcup\limits_{j=0}^{k}{D}_{s_j},\quad D^*:=\ssD\setminus\bigcup_{j=k+1}^{n}\overline{D}_{s_j}.
\end{equation}
Donnons-nous des chemins simples r\'{e}guliers analytiques r\'{e}els
$$\sigma_{j}:[0,1]\to \ssD ,\quad j={0},\ldots,n$$
dont les images $|\sigma_{j}|$ v\'{e}rifient :
\begin{enumerate}[-]
\item $|\sigma_{j}|\cap|\sigma_{k}|=\emptyset$, pour $j\neq k$;
\item $|\sigma_{j}|\cap \overline{D}_{s_{k}}=\emptyset$, pour $k\neq 0,j$;
\item  $|\sigma_{j}|\cap \overline{D}_{s_j}=\{\sigma_{j}(0)\}\subset {\partial D_{s_{j}}}$ si $j=1,\ldots,k$ et $|\sigma_{j}|\cap \overline{D}_{s_j}$
est un chemin radial de $\sigma_{j}(0)=s_{j}$ \`{a} $\sigma_{j}(1/2)\in
\mc C_{j}$ si $j=k+1,\ldots,n$;
\item $|\sigma_{j}|\cap \overline{D}_{s_{0}}=\{\sigma_{j}(1)\}\subset {\partial D_{s_{0}}}$.
\end{enumerate}

\vspace{-1,5em}

\begin{center}
\begin{figure}[htp]
\includegraphics[width=8cm]{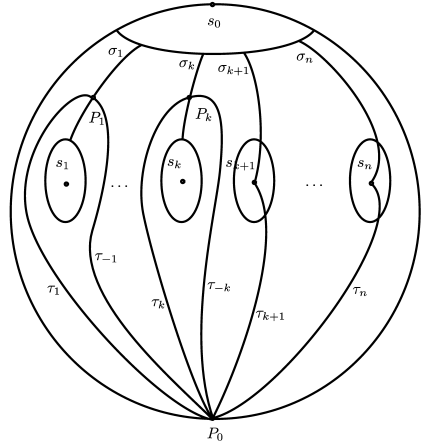}
\vspace{-3cm} {\caption{Num\'{e}rotation des points singuliers et
description des courbes auxiliaires.}}
\end{figure}
\end{center}
Consid\'{e}rons aussi le sous-ensemble de $D$ d\'{e}fini par
$$D^{\circ}=D^*\setminus\bigcup\limits_{j=1}^{n}|\sigma_{j}|\,.$$
{Fixons un point de base un  point base $P_0\in D^{\circ}$ et notons
encore $\pi_D : \mc T_{\eta_{1}}(D^*) \rightarrow D^*$ la
restriction au bloc de Milnor $\mc T_{\eta_{1}}(D^*)$, cf.
(\ref{bloc.mil}), de la fibration (\ref{fibr}) construite en
(\ref{decompdiv}).} Soit $\Sigma\subset{\pi}_{D}^{-1}(P_{0})$ un
disque ouvert dont le bord $\partial
\Sigma=\overline{\Sigma}\setminus\inte{\overline{\Sigma}}$ est lisse
par morceaux. La restriction $\F_{D^{\circ}} := \wF_{|
{\pi}_D^{-1}(D^\circ)}$ du feuilletage ${\wF}$ \`{a} ${\pi}
_D^{-1}(D^\circ)$ est triviale, car $D^{\circ}$ est simplement
connexe. Plus pr\'{e}cis\'{e}ment, il existe une constante $C>0$ %ne d\'{e}pendant que de $\F$ et du choix des courbes $\mc C_{j}$, $\sigma_{j}$
telle que si $|F(Q)|\le C$, la restriction de ${\pi}_{D}$ \`{a} la
feuille de $\F_{D^{\circ}}$ passant par $Q\in{\pi}_{D}^{-1}(P_{0})$
est un diff\'{e}omorphisme sur $D^{\circ}$. Supposons que l'on a :
$$\|\Sigma \|_F := \max\{|F(Q)|,\ Q\in\overline{\Sigma}
\}<C.$$ \noindent Consid\'{e}rons la r\'{e}union $\mc
B_{\Sigma}(D^{\circ})=\mathrm{Sat}(\Sigma,{\pi}_{D}^{-1}(D^{\circ}))$
des feuilles de $\F_{D^{\circ}}$ qui intersectent $\Sigma$.
Par rel\`{e}vement des chemins dans les feuilles suivant ${\pi}_{D}$, on
construit un biholomorphisme
$\varphi^{\circ}:D^{\circ}\times\Sigma\to \mc
B_{\Sigma}(D^{\circ})\subset {\pi}_{D}^{-1}(D^{\circ})$ qui conjugue
le feuilletage horizontal (i.e. dont les feuilles sont
$D^{\circ}\times \{Q\}$ avec $Q\in\Sigma$) \`{a} $\F_{D^{\circ}}$ et qui
v\'{e}rifie :
$${\pi}_{D}(\varphi^{\circ}(P,Q))=P\quad\textrm{et}
\quad \varphi^{\circ}(P_{0},Q)=Q\quad\textrm{pour tout }\quad  P\in
D^{\circ}\textrm{ et }Q\in\Sigma.$$
Comme le bord de $D^{\circ}$ est lisse par morceaux, il existe une
uniformization  $\psi^{\circ}:\mb D\iso \inte{D}$ de l'int\'{e}rieur de
$D^{\circ}$ dans $D$ qui s'\'{e}tend contin\^{u}ment au bord. Ainsi le
biholomorphisme
$\phi^{\circ}:=\varphi^{\circ}\circ(\psi^{\circ}\times\mathrm{id}_{\Sigma}):\mb
D\times\Sigma\to\mc B_{\Sigma}(\inte{D})$ s'\'{e}tend contin\^{u}ment \`{a}
$\phi:\overline{\mb D}\times\Sigma\to \overline{T}_{\eta}(D^*)$.
Soit $\overline{\mc B}_{\Sigma}(D^*)$ l'image de $\phi$ et
$$\mc B_{\Sigma}(D^*)=(\overline{\mc B}_{\Sigma}(D^*)\cap{\pi}_{D}^{-1}(D^*))\setminus  D.$$

Le type d'homotopie $\mc B_{\Sigma}(D^*)$ d\'{e}pend de la taille et de
la rugosit\'{e} de $\Sigma $. En effet, pour $P\in D^\circ$, la fibre
${\pi}_D^{-1}(P)\cap \mc B_{\Sigma}(D^*)$ est l'image de $\Sigma $
par l'application d'holonomie $h_\alpha : {\pi}_D^{-1}(P_0)
\rightarrow {\pi}_D^{-1}(P)$ du feuilletage le long d'un chemin
$\alpha $ trac\'{e} dans $D^\circ $ joignant $P_0$ \`{a} $P$. Par contre si
$P$ est un point d'une courbe $| \sigma _j |$, alors ${\pi}_D
^{-1}(P)\cap \mc B_{\Sigma}(D^*)$  est \'{e}gal \`{a} l'union $ h_\alpha
(\Sigma )\cup h_\beta (\Sigma )$ o\`{u} $\alpha $ et $\beta $ sont des
chemins d'origine $P_0$ contenu dans $D^\circ $, sauf leur extr\'{e}mit\'{e}
qui est commune, \'{e}gale \`{a} $P$ et atteinte par des cot\'{e}s diff\'{e}rents de
la courbe $|\sigma _j|$. A priori cette union pourrait ne pas \^{e}tre
simplement connexe. Une cons\'{e}quence directe de la proposition
(\ref{prop31}) est

\begin{obs}\label{control-rugosite}
{Suivant les notations de (\ref{doublecrocher}) et (\ref{diametre}),
il existe une constante $C_{D}>0$ telle que pour tout
$P\in\overline{D^*}$, la fibre $\Sigma_{P}:= {\pi}_D^{-1}(P)\cap
\overline{\mc B}_{\Sigma}(D^*)$ v\'{e}rifie les estimations suivantes
:\begin{equation}\label{CondRugGabarit} \mathbf{e}_F\,\left(
\Sigma_{P} \right )\; \leq \; \{\!\{ \,\mathbf{e}_F(\partial \Sigma
)+ C_{D} \| \Sigma
 \|_F \,\}\!\} \,\quad\textrm{et}\quad \|\Sigma_{P}\|_{F}\le C_{D}\|\Sigma\|_{F}.
\end{equation}}
\end{obs}
\noindent On en d\'{e}duit que si
$\max\{\mathbf{e}_{F}(\Sigma),\|\Sigma\|_{F}\}$ est assez petit
alors toutes les fibres $\Sigma_{P}$ ont une rugosit\'{e} finie et sont
donc \'{e}toil\'{e}es. Plus pr\'{e}cis\'{e}ment,

\begin{lema}\label{4.1.2.}
{Il existe une constante $C'>0$} %ne d\'{e}pendant que de $\F$ et du choix des courbes $\mc  C_{j}$, $\sigma_{j}$,
telle que si $\max\{\mathbf{e}_{F}(\Sigma),\|\Sigma\|_{F}\}\le C'$,
alors pour tout $0<\eta<\eta_{1}$ assez petit $\mc B_{\Sigma}(D^*)$
se r\'{e}tracte sur ${\mc T}_{\eta}^*(D^*)$. En particulier, l'inclusion
${\mc T}_{\eta}^*(D^*)\subset \mc B_{\Sigma}(D^*)$ induit un
isomorphisme au niveau des groupes fondamentaux et ceux-ci qui
admettent la pr\'{e}sentation suivante :
\begin{equation}\label{pi1D}
    \left< a_0, \ldots , a_{n}, c \;\left| \; a_0\cdots
a_{n} = c^{\nu},\;[a_{r},c]=1,\;r =0,\ldots,n \right.\right>
\;\simeq\; \mathbb{Z}^{*n}\times \mathbb{Z}\,,
\end{equation}
o\`{u} $\{a_{j}\}_{j=0}^{n}$  sont les relev\'{e}s par une section de la
fibration triviale ${\pi}_{D}:\mc B_{\Sigma}(D^*)\to D^*$ d'un
syst\`{e}me de g\'{e}n\'{e}rateurs g\'{e}om\'{e}triques de $\pi_{1}(D^*)$; $c$ est le
g\'{e}n\'{e}rateur du groupe fondamental de la fibre g\'{e}n\'{e}rique qui
s'identifie \`{a} $\mb D^*$ et $-\nu$ est \'{e}gal \`{a} l'auto-intersection
$(D,D)$. D'autre part, chaque composante connexe du bord de ${\mc
T}_{\eta}^*(D^*)$ est contenue dans une composante connexe du bord
de $\overline{\mc B}_{\Sigma}(D^*)$. Cette inclusion induit un
isomorphisme au niveau des groupes fondamentaux, ceux-ci admettant
pour pr\'{e}sentation : $\langle a_{j},c|\ [a_{j},c]=1\rangle$,
$j=0,\ldots,n$.
\end{lema}
\begin{proof}
{Supposons  $C'>0$ assez petit pour que les fibres $\Sigma_{P}$
soient toutes \'{e}toil\'{e}es. Pour d\'{e}finir la r\'{e}traction, il suffit
d'int\'{e}grer un rel\`{e}vement du champ radial $\partial/\partial z$ par
$F:\mc T_{\eta_{1}}\to\D_{\eta_{1}}$, qui est tangent aux fibres de
$\pi _D$.} {Remarquons que par d\'{e}finition les lacets
$\{\underline{\sigma}_{j}\vee \partial D_{s_{j}}\vee
\underline{\sigma}_{j}^{-1}\}_{j=0}^{n}$ {forment} un syst\`{e}me de
g\'{e}n\'{e}rateurs g\'{e}om\'{e}tri\-ques de $\pi_{1}(D^*)$, o\`{u}
$\underline{\sigma}_{j}=\sigma_{j}$ si $j=0,\ldots,k$ et
$\underline{\sigma}_{j}=\sigma_{j|[1/2,1]}$ si $j=k+1,\ldots,n$. Le
seul point de la pr\'{e}sentation (\ref{pi1D}) qui n'est pas \'{e}vident est
le fait que l'exposant $\nu$ de $c$ dans la relation $a_{0}\cdots
a_{n}=c^{\nu}$ soit {soit \'{e}gal \`{a} $-(D, D)$.} Ceci est facile a voir
dans le cas $\nu=1${; c' est} une cons\'{e}quence directe de la
description de $\mc O_{\mb P^{1}}(-\nu)$ comme le quotient $\mc
O_{\mb P^1}(-1)/\mb Z_{\nu}$ dans le cas g\'{e}n\'{e}ral, voir aussi
\cite{Mumford}.}
\end{proof}
\noindent{Notons que} chaque composante connexe
$${\partial_{j}\mc B_{\Sigma}(D^*):=\mc
B_{\Sigma}(D^*)\cap\pi_{D}^{-1}(\partial D_{s_{j}})}$$ du bord
$\partial\mc B_{\Sigma}(D^*)=\bigsqcup_{j=0}^{k}\partial_{j}\mc
B_{\Sigma}(D^*)$ est
 de type suspension si $j=1,\ldots,k$ et de type
multi-suspension si $j=0$, cf. (\ref{suspension}). En {effectuant}
le proc\'{e}d\'{e} du ``rabotage'' d\'{e}crit dans (\ref{rabote}) on
construit un sous-ensemble de type suspension $\partial_{0}\mc
B_{D}$ qui est $1$-$\F$-connexe dans $\partial_{0}\mc
B_{\Sigma}(D^*)$.

\begin{obs}\label{obs1} {A l'\'{e}tape suivante nous construisons} un voisinage $\mc B_{\Sigma}(\mf M_{j})\sqcup (\mf
M_{j}\cup D_{s_{j}})\subset {\mc T}_{\eta_{1}}(\mf M_{j})$ de chaque
branche morte $\mf M_{j}$, $j=k+1,\ldots,n$, adjacente \`{a} $D$ tel que
:
\begin{enumerate}[(i)]
\item $\mc B_{\Sigma}(\mf M_{j})$ est $\F$-satur\'{e} dans ${\mc T}_{\eta_{1}}^*(\mf M_{j})$,
\item le bord $\partial\mc B_{\Sigma}(\mf M_{j})$
est contenu dans ${\pi}_{D}^{-1}({\partial D_{s_j}})\cap
\overline{\mc B}_{\Sigma}(D^*)$ {et  est} invariant par l'holonomie
$h_{j}$ de ${\partial D_{s_j}}\subset D$,
\item $\mc B_{\Sigma}(\mf M_{j})$ est maximal pour la inclusion parmi les sous-ensembles de ${\mc T}_{\eta_{1}}^*(\mf M_{j})$ qui v\'{e}rifient les  propri\'{e}t\'{e}s (i) et (ii).
\end{enumerate}
{Cela nous permet de d\'{e}finir :}
$$\mc B_{\Sigma}(D):=\mc B_{\Sigma}(D^*)\cup\bigcup_{j=k+1}^{n}\mc B_{\Sigma}(\mf M_{j})$$
et
\begin{equation}\label{bloctotal}
\mc B_{D}:=(\mc B_{\Sigma}(D)\setminus\partial_{0}\mc
B_{\Sigma}(D))\cup\partial_{0}\mc B_{D},
\end{equation}
en choisissant le disque $\Sigma$ de la mani\`{e}re suivante :
\begin{itemize}
\item si $D$ n'est pas la composante centrale de $\mc D$, alors $k\ge
1$ et on num\'{e}rote les singularit\'{e}s pour avoir
${\pi_{D}(V)\subset\partial D_{s_{1}}}$;  on prend pour
$\sigma_{1}(0)\in{\partial D_{s_{1}}}$ le ``point de rupture'' de la
suspension $V$; on choisit $P_{0}$ sur ${\partial
D_{s_{1}}}\setminus\{\sigma_{1}(0)\}$ et on consid\`{e}re
$\Sigma:={\pi}_{D}^{-1}(P_{0})\cap V$;
\item si $D$ est la composante centrale de $\mc D$, on prend pour $P_0$ un point
quelconque de $D^{\circ}$ et $\Sigma$ assez
petit de rugosit\'{e} nulle.
\end{itemize}
\noindent Il est clair que $\partial \mc
B_{D}=\bigcup\limits_{i=0}^{k} V_{i}$ est de type suspension et
$V_{1}=V$ dans le cas o\`{u} $D$ n'est pas la composante centrale de
$\mc D$. Cela prouve les propri\'{e}t\'{e}s (3), (3') du th\'{e}or\`{e}me
(\ref{teoexistence}). Quant aux propri\'{e}t\'{e}s (4) et (4'), elles
r\'{e}sultent imm\'{e}diatement de la remarque (\ref{control-rugosite}).
\end{obs}

\noindent {\textbf{Etape 2 : construction et description de $\mc
B_{\Sigma}({\mf M_{j}})$.}} Il est bien connu que l'holonomie $h_{j}$
{de $\F$ le long du lacet} ${\partial D_{s_{j}}}\subset D$ est
p\'{e}riodique. {Notons} $p_{j}$ la p\'{e}riode {positive minimale} de
$h_{j}$. Il est aussi bien connu que $\F$ {poss\`{e}de au voisinage de
chaque branche morte $\mf M_j$, une int\'{e}grale premi\`{e}re holomorphe
$f_j$, avec $f_{j\;|\mf M_j}\equiv 0$. D'apr\`{e}s la proposition
(\ref{prop31}) et les propri\'{e}t\'{e}s (\ref{rugsom})} il existe une
constante $C''>0$
%ne d\'{e}pendant que de $\F$ et du choix
%des courbes $\mc C_{j}$, $\sigma_{j}$,
telle que si
$\max\{\mathbf{e}_{F}(\Sigma),\|\Sigma\|\}\le C''$ alors, pour
chaque $j=k+1,\ldots,n$, l'ensemble
$$\Sigma_{j}:=\Sigma_{\sigma_{j}(0)}\cap h_{j}(\Sigma_{\sigma_{j}(0)})\cap\cdots\cap h_{j}^{p_{j}-1}(\Sigma_{\sigma_{j}(0)})$$
est un disque conforme, invariant par l'holonomie $h_{j}$, dont le
bord est de $F$-rugosit\'{e} finie. {Consid\'{e}rons}
${C_{D,\,0}}:=\min\{C,C',C''\}>0$ et {d\'{e}finissons}
$$\mc B_{\Sigma}({\mf M_{j}}):=f^{-1}_{j}(\Delta^*_{j})\setminus {\pi}_{D}^{-1}(D\setminus \overline{D}_{s_{j}}),$$
o\`{u} $\Delta_{j}$ est l'image du  disque {ouvert} $\Sigma_{j}$ par une
int\'{e}grale premi\`{e}re holomorphe et primitive $f_{j}$ de $\F$ d\'{e}finie
au voisinage de $\mf M_{j}\cup
\overline{D}_{s_{j}}=f_{j}^{-1}(0)\setminus(D\setminus\overline{D}_{s_{j}})$.

Remarquons que ${\mc T}_{\eta}^*(\mf M_{j})$ est de la forme $\mc
B_{\Sigma}(\mf M_{j})$ pour le feuilletage donn\'{e} par $dF$ et {pour}
un disque $\Sigma$ tel que $F(\Sigma)=\D_{\eta}$.

\begin{prop}\label{branches-mortes} Si
$\max\{\mathbf{e}_{F}(\Sigma),\|\Sigma\|\}\le C''$, alors $\mc
B_{\Sigma}(\mf M_{j})$ v\'{e}rifie les propri\'{e}t\'{e}s suivantes :
\begin{enumerate}[(a)]
\item $\partial \overline{\mc B}_{\Sigma}(\mf M_{j})=\overline{\mc B}_{\Sigma}(\mf M_{j})\cap {\pi}_{D}^{-1}(\overline{D}^*)=\overline{\mc B}_{\Sigma}(\mf M_{j})\cap {\pi}_{D}^{-1}(\partial D_{s_{j}})$ est une $3$-vari\'{e}t\'{e} \`{a} bord qui est un voisinage de $\partial D_{s_{j}}$ dans ${\pi}_{D}^{-1}(\partial D_{s_{j}})$;
\item il existe une carte holomorphe $\tau_{j}=(x_{j},y_{j}):{\mc T}_{\eta_{1}}(s_{j})\to\DDD^*$
telle que $D_{s_{j}}=\tau_{j}^{-1}(\{y_{j}=0\})$,
$x_{j}\circ{\pi}_{D}=x_{j}$ et la restriction \`{a} ${\mc
T}_{\eta_{1}}(s_{j})$ de la int\'{e}grale premi\`{e}re holomorphe $f_{j}$
s'\'{e}crit sous la forme
$f_{j}(x_{j},y_{j})=x_{j}^{q_{j}}y_{j}^{p_{j}}$, o\`{u}
$-q_{j}/p_{j}\in\mb Q_{<0}$ est la fraction r\'{e}duite de l'indice de
Camacho-Sad du feuilletage $\wt \F$ par rapport \`{a} $D$ au point
singulier $s_{j}$;
\item il existe un $C^{1}$-diff\'{e}omorphisme $\Phi_{j}:
\overline{\mc B}_{\Sigma}^*(\mf M_{j})\to \overline{\mb
D}\times\overline{\mb D}^*$, {$\overline{\mc B}_{\Sigma}^*(\mf
M_{j}) :=\overline{\mc B}_{\Sigma}(\mf M_{j})\setminus \mc D$,} qui
conjugue la restriction de $\F$ au feuilletage horizontal et qui
envoie $\partial\mc B_{\Sigma}(\mf M_{j})$ sur
$\partial\overline{\mb D}\times\overline{\mb D}^*$. En particulier,
la restriction de $\F$ \`{a} $\partial\mc B_{\Sigma}(\mf M_{j})$ est un
feuilletage en cercles dont chaque feuille est le bord d'une feuille
de $\F_{|\mc B_{\Sigma}(\mf M_{j})}$;
\item les g\'{e}n\'{e}rateurs $a_{j}$ et $c$ du groupe fondamental de la composante
connexe {de} $\partial\overline{\mc B}_{\Sigma}(D^*)$ qui se
r\'{e}tracte sur $\partial\mc B_{\Sigma}(\mf M_{j})$ s'identifient
respectivement aux classes d'homotopie des lacets (positivement
orient\'{e}s) $\tau_{j}^{-1}(\{|x_{j}|=\varepsilon_{j},\
y_{j}=\varepsilon_{j}\})$ et
$\tau_{j}^{-1}(\{x_{j}=\varepsilon_{j},\
|y_{j}|=\varepsilon_{j}\})${, avec }$\varepsilon_{j}>0$ assez petit.
{De plus}, le groupe fondamental de $\mc B_{\Sigma}(\mf M_{j})$
admet la pr\'{e}sentation suivante :
\begin{equation}\label{pi1M}
\pi _1(\mc B_{\Sigma}(\mf M_{j}))\;\cong \; \left< a_{j},c
\;\left|\; [a_j, c] = 1,\; a_j^{p_j}=c^{q_{j}}\right. \right>\;=\;
\Z\langle a_{j}^{m_{j}}c^{-n_{j}}\rangle,
\end{equation}
o\`{u} $m_{j},n_{j}\in\mb N$, $n_{j}<p_{j}$, v\'{e}rifient {la relation}
$m_{j}p_{j}-n_{j}q_{j}=1$.
\item la restriction de ${\pi}_{D}$ \`{a}
$\partial\overline{\mc B}^*_{\Sigma}(\mf M_{j})$ se prolonge en une
$\overline{\mb D}^*$-fibration {de} Seifert ${\sigma}_{\mf
M_{j}}:\mc B_{\Sigma}(\mf M_{j})\to \overline{D}_{s_{j}}$ de classe
$C^{1}$ {et} transverse \`{a} $\F$, avec exactement une fibre
exceptionnelle ${\sigma}_{\mf M_{j}}^{-1}(s_{j})${, qui est} de type
$(p_{j},n_{j})$.
\end{enumerate}
\end{prop}

\begin{obs}\label{def-seifert}
La notion de $\overline{\D}^*$-fibration Seifert que nous
consid\'{e}rons ici est l'analogue \`{a}  la d\'{e}finition classique de
fibration Seifert en cercles. Plus pr\'{e}cis\'{e}ment, une application
diff\'{e}rentiable $\sigma:V\to S$ de une $4$-vari\'{e}t\'{e} \`{a} bord $V$ sur une
surface $S$ est appel\'{e}e \textbf{$\overline{\D}^*$-fibration de
Seifert} si tout point $P$ de $S$ admet un voisinage $U$ tel que la
restriction de $\sigma$ \`{a} $\sigma^{-1}(U)$ est \'{e}quivalente au mod\`{e}le
local $\sigma_{a,b}:\D\times \overline{\D}^*\to\D$ d\'{e}fini par
${\sigma_{a,b}(x,y)={x^{a}}{(y/|y|)^{-b}}}$,
%$\sigma_{a,b}(x,y)=\frac{x^{a}}{(y/|y|)^{b}}$
pour certains entiers
$0\le b<a$ premiers entre eux.
Remarquons que la restriction de $\sigma$ \`{a}
$\D\times\partial\overline{\D}$ (i.e. \`{a} $|y|=1$) est une
$\mathbb{S}^{1}$-fibrations de Seifert au sens classique, cf.
\cite{O72}.
On dit que la fibre correspondant \`{a} $\{x=0\}$ est {{\bf
exceptionnelle de type $(a,b)$}}. Toutes les autres fibres de
$\sigma^{-1}(U)$ admettent des voisinages tubulaires model\'{e}s par
$\sigma_{1,0}$, c'est \`{a} dire {sur lesquels} $\sigma$ est une
$\overline{\D}^*$-fibration localement triviale.
Comme dans le cas classique, pour toute $\overline{\D}^*$-fibration
Seifert $\sigma:V\to S$ il existe un rev\^{e}tement ramifi\'{e}
$\rho:\widehat{S}\to S$ tel que le pull-back
$\widehat{\sigma}:=\rho^*\sigma$ de $\widehat{V}:=\rho^*V$ en
$\widehat{S}$ est une $\overline{\D}^*$-fibration localement
triviale. En effet, il suffit de consid\'{e}rer le mod\`{e}le local
$\sigma_{a,b}$ et le rev\^{e}tement ramifi\'{e} $\rho_{a}:\D\to\D$ donn\'{e} par
$\rho_{a}(x)=x^{a}$.
\end{obs}

\begin{dem2}{de la proposition (\ref{branches-mortes})}
Les assertions (a) et (b) sont \'{e}videntes \`{a} partir de la d\'{e}finition
de $\mc B_{\Sigma}(\mf M_{j})$.
Pour montrer (c), remarquons d'abord que si $D'\subset\mf M_{j}$ est
la composante extr\'{e}male de valence $1$, alors il existe un
$C^{1}$-diff\'{e}omorphisme $\mc B_{\Sigma}(\mf M_{j})\cap {\mc
T}_{\eta_{1}}(D')\iso\D\times \D^*$ qui conjugue $\F$ au feuilletage
horizontal. Si $D'\subset\mf M_{j}$ est une composante de valence
$2$ ayant pour singularit\'{e}s du feuilletage $s',s''\in D'$ alors
$D'^* {:=} D'\setminus (D_{s'}'\cup D_{s''}')$ est une couronne {et}
poss\`{e}de un champ de vecteurs r\'{e}el dont chaque orbite relie un des
cercles de son bord $\partial D_{s'}'$ \`{a} l'autre $\partial
D_{s''}'$. En int\'{e}grant le relev\'{e} de ce champ aux feuilles de $\F$
via {$\pi _{D'}$} on obtient des $C^{1}$-diff\'{e}omorphismes
\begin{eqnarray*}
\mc B_{\Sigma}(\mf M_{j})\cap {\mc T}_{\eta_{1}}(D')&\iso& \left(\mc B_{\Sigma}(\mf M_{j})\cap\rho_{D'}^{-1}(\partial D'_{s'})\right)\times [0,1]\\
&\iso & \left(\mc B_{\Sigma}(\mf M_{j})\cap\rho_{D'}^{-1}(\partial
D'_{s''})\right)\times [0,1].
\end{eqnarray*}
D'autre part, pour {toute singularit\'{e} $s$ du feuilletage situ\'{e}e sur
$\mf M_{j}$, n\'{e}ces\-sai\-rement point d'intersection de} deux
composantes $D'$ et $D''$ de $\mc E$, la restriction de {$\wF$} est
d\'{e}crite par l'exemple {1.C) de la section \ref{subsect.exemples}}.
On obtient donc des diff\'{e}omorphismes
\begin{eqnarray*}
\mc B_{\Sigma}(\mf M_{j})\cap {\mc T}_{\eta_{1}}(s)&\iso& \left(\mc B_{\Sigma}(\mf M_{j})\cap\rho_{D'}^{-1}(\partial D'_{s})\right)\times [0,1]\\
&\iso & \left(\mc B_{\Sigma}(\mf M_{j})\cap\rho_{D''}^{-1}(\partial
D_{s}'')\right)\times [0,1].
\end{eqnarray*}
Il est clair que les feuilles de $\F_{|\mc B_{\Sigma}(\mf M_{j})}$
sont obtenues en collant bord \`{a} bord un disque avec une succession
des couronnes. En fait les diff\'{e}omorphismes pr\'{e}c\'{e}dents se recollent
pour donner un $C^{1}$-diff\'{e}omorphisme global
$$\Phi_{j}:\overline{\mc B}_{\Sigma}^*(\mf M_{j})\iso\DDD^*$$
satisfaisant les propri\'{e}t\'{e}s requises. Nous laissons les d\'{e}tails de
cette construction au lecteur.

La description des g\'{e}n\'{e}rateurs $a_{j}$ et $c$ dans l'assertion (d)
est \'{e}vidente d'apr\`{e}s le lemme (\ref{4.1.2.}). D'autre part, il est
clair que la restriction de $f_{j}$ \`{a} $\mc B_{\Sigma}(\mf M_{j})$
est une fibration localement triviale en disques (car  elle est
homotope \`{a} $\mathrm{pr}_{2}\circ\Phi_{j}$). On d\'{e}duit de la suite
exacte d'homotopie que
$\gamma\mapsto\frac{1}{2i\pi}\int_{\gamma}\frac{df_{j}}{f_{j}}$
induit un isomorphisme entre $\pi_{1}(\mc B_{\Sigma}(\mf M_{j}))$ et
$\mb Z$. D'apr\`{e}s la description locale des g\'{e}n\'{e}rateurs $a_{j}$ et
$c$, on a: $\frac{1}{2i\pi}\int_{a_{j}}\frac{df_{j}}{f_{j}}=q_{j}$
et $\frac{1}{2i\pi}\int_{c}\frac{df_{j}}{f_{j}}=p_{j}$. D'o\`{u}
$a_{j}^{p_{j}}=c^{q_{j}}$ dans $\pi_{1}(\mc B_{\Sigma}(\mf
M_{j}))\cong\mb Z$ et ce groupe admet $a_{j}^{m_{j}}c^{-n_{j}}$
comme g\'{e}n\'{e}rateur {d\`{e}s que} $m_{j}p_{j}-n_{j}q_{j}=1$.

Pour prouver l'assertion (e) consid\'{e}rons l'application $\phi_{j}:\mb
S^{1}\times\DD^*\to\mb S^{1}\times\DD^*$ donn\'{e}e par
$\phi_{j}(e^{i\theta_{1}},r
e^{i\theta_{2}})=\left(e^{i(m_{j}\theta_{1}+n_{j}\theta_{2})},r^{q_{j}}e^{i(q_{j}\theta_{1}+p_{j}\theta_{2})}\right)$,
o\`{u} $\mb S^{1}=\partial\DD$. C'est un diff\'{e}omorphisme dont l'inverse
est :
$$\phi_{j}^{-1}(e^{i\theta_{1}},re^{i\theta_{2}})=\left(e^{i(p_{j}\theta_{1}-n_{j}\theta_{2})},r^{\frac{1}{q_{j}}}e^{i(-q_{j}\theta_{1}+m_{j}\theta_{2})}\right).$$
La restriction de $\Phi_{j}$ \`{a} $\partial \mc B_{\Sigma}(\mf M_{j})$
co\"{\i}ncide avec $\phi_{j}\circ\tau_{j}:\partial B_{\Sigma}(\mf
M_{j})\to\mb S^{1}\times\DD^*$. L'application
${\varsigma}_{j}:=x_{j}\circ{\pi}_{D}\circ\Phi_{j}^{-1}:\mb
S^{1}\times\DD^*\to\mb S^{1}$, qui s'\'{e}crit aussi sous la forme
${\varsigma}_{j}(e^{i\theta_{1}},r
e^{i\theta_{2}})=e^{i(p_{j}\theta_{1}-n_{j}\theta_{2})}$, admet un
prolongement $\underline{{\varsigma}}_{j}:\DDD^*\times\DD$ obtenu en
posant
$\underline{{\varsigma}}_{j}(x, re^{i\theta_2}):= x^p e^{-in
\theta_2}=\sigma_{p,n}(x,y)$, $y=r
e^{i\theta_{2}}\in\overline{\mathbb{D}}^{*}(1)$. Celui-ci satisfait
visiblement les propri\'{e}t\'{e}s suivantes :
\begin{enumerate}
  \item[(a)]\label{pr1} la restriction de $\underline{{\varsigma}}_{j}$ \`{a}
  chaque feuille $\overline{\mathbb{D}} \times
\{re^{i\theta_2}\}$ du feuilletage horizontal $\mathcal{H}$ est un
rev\^{e}tement ramifi\'{e} avec l'origine comme unique point de ramification
(non-triviale,  sauf si  $p_{j}= 1$),
  \item[(b)]\label{pr2} chaque fibre de
  $\underline{{\varsigma}}_{j}$, y compris la fibre
exceptionnelle $\underline{{\varsigma}}_{j}^{-1}(0)$, est une courbe
lisse $C^\infty$ transverse \`{a} chaque feuille de $\mathcal{H}$.
\end{enumerate}
Posons ${\sigma}_\mathfrak{M_{j}} :=
\underline{{\varsigma}}_{j}\circ \Phi_{j}$. Pour achever la preuve
il suffit {de voir} que $p_{j}\neq 1$. Visiblement $p_{j}$ est la
p\'{e}riode de l'application d'holonomie $h_{j}$ du feuilletage $\wt \F$
le long du lacet $\partial D_{s_{j}}$. Or on sait que l'holonomie au
point de branchement d'une branche morte n'est jamais l'identit\'{e},
cf. \cite{Mattei-Salem} lemme (6.2.5).
\end{dem2}

\noindent\textbf{Etape 3 : r\'{e}duction de la preuve du th\'{e}or\`{e}me
(\ref{teoexistence}).}{  Les assertions (3), (3'), (4), (4') ont \'{e}t\'{e}
prouv\'{e}es en (\ref {obs1}). Montrons l'assertion (1).} D'apr\`{e}s la
pr\'{e}sentation (\ref{pi1M}), si
$\max\{\mathbf{e}_{F}(\Sigma),\|\Sigma\|\}\le C''$ et si
$0<\eta<\eta_{1}$ est assez petit, alors {les inclusions} ${\mc
T}_{\eta}^*(\mf M_{j})\subset\mc B_{\Sigma}(\mf M_{j})$ et $\partial
{\mc T}_{\eta}^*(\mf M_{j})\subset\partial \mc B_{\Sigma}(\mf
M_{j})$ induisent des isomorphismes \`{a} niveau des groupes
fondamentaux. Remarquons aussi que l'on a : $\partial {\mc
T}_{\eta}^*(\mf M_{j})={\mc T}_{\eta}^*(\mf M_{j})\cap {\mc
T}_{\eta}^*(D^*)$ et $\partial\mc B_{\Sigma}(\mf M_{j})={\mc
B}_{\Sigma}(\mf M_{j})\cap \overline{\mc B}_{\Sigma}(D^*)$. {Ainsi
l'assertion (1) de} (\ref{teoexistence}) d\'{e}coule {de la remarque}
(\ref{rVK}) et du lemme (\ref{4.1.2.}).

{Pour achever la preuve du th\'{e}or\`{e}me (\ref{teoexistence}) il ne reste
plus qu'\`{a} prouver l'assertion (2) de $\F$-adaptabilit\'{e} du bloc $\mc
B_{D}$, i.e.
\begin{itemize}
\item[(i)] le bord $\partial\mc B_{D}$ est incompressible,
\item[(ii)] $\F$ est transverse \`{a} $\partial \mc B_{D}$,
\item[(iii)] chaque feuille de $\F_{|\mc B_{D}}$ est incompressible dans $\mc
B_{D}$,
\item[(iv)] $\partial\mc B_{D}\uco{\F}\mc B_{D}$.
\end{itemize}
En appliquant le  th\'{e}or\`{e}me de Seifert-Van Kampen et {en utilisant}
les pr\'{e}sentations (\ref{pi1D}) et (\ref{pi1M}) on obtient la
pr\'{e}sentation explicite suivante du groupe {fondamental} $\pi_{1}(\mc
B_{D})$ :
\begin{multline}\label{pi1B}
\Big< a_0, \ldots , a_{n}, c \;\Big| \; a_0\cdots a_{n} =
c^{\nu},\;[a_{r},c]=1, \; a_{j}^{p_{j}}=c^{q_{j}},
\\
\;r =0,\ldots,n,\;j=k+1\,,\ldots, n \Big>.
\end{multline}
Visiblement le groupe fondamental  $\pi_{1}(V_{j})\cong\langle
a_{j},c\, |\, [a_{j},c]=1\rangle$ de chaque composante connexe
$V_{j}$ du bord de $\mc B_{D}$, $j=0,\ldots,k$, s'injecte dans
$\pi_{1}(\mc B_{D})$; D'o\`{u} (i).} {La propri\'{e}t\'{e} (ii) est \'{e}vidente par
construction. Nous allons voir maintenant que l'assertion (iii) est
une cons\'{e}quence de (iv).}

Soit $L$ une feuille de $\F_{|\mc B_{D}}$. Comme par construction de
$\mc B_{D}$, $L$ intersecte $\partial\mc B_{D}$, on peut consid\'{e}rer
un point $p\in\ L\cap\partial \mc B_{D}$. Soit $s$ la singularit\'{e} de
$D$ associ\'{e}e \`{a} la composante $\partial_{p}\mc B_{D}$ de $\partial\mc
B_{D}$ qui contient le point $p$. On distingue deux cas :
\begin{itemize}
\item si $\lambda_{s}\in\mathbb{Q}_{<0}$ et $s$ est lin\'{e}arisable, alors la description  finale de l'exemple 1.C. implique que $L\cap\partial_{p}\mc B_{D}$ est incompressible dans $\partial_{p}\mc B_{D}$;
\item sinon, ou bien $\lambda_{s}\notin\mathbb{Q}_{<0}$ auquel cas
$s$ est lin\'{e}arisable d'apr\`{e}s les hypoth\`{e}ses du th\'{e}or\`{e}me principal,
ou bien $s$ est une singularit\'{e} r\'{e}sonante.  Dans les deux cas,
$L\cap\partial_{p}\mc B_{D}\cong\R$ est simplement connexe, donc
incompressible dans $\partial_{p}\mc B_{D}$.
\end{itemize}
En combinant l'incompressibilit\'{e} de $L\cap\partial_{p}\mc B_{D}$
dans $\mc \partial_{p}\mc B_{D}$ avec l'assertion (iv) et la
transitivit\'{e} de la relation $\uco{\F}$ on obtient les relations
$$\{p\}\uco{\F}\partial_{p}\mc B_{D}\uco{\F}\mc B_{D},$$
c'est \`{a} dire l'incompressibilit\'{e} de $L$ dans $\mc B_{D}$.\\

Effectuons une derni\`{e}re r\'{e}duction du probl\`{e}me en utilisant la
relation $\partial_{0}\mc B_{D}\ucon\partial_{0}\mc B_{\Sigma}(D)$ :
par transitivit\'{e}{, la propri\'{e}t\'{e} (iv) est impliqu\'{e}e par %l'assertion
la proposition suivante :
\begin{prop}\label{v} Le bord de $\mc B_{\Sigma}(D)$ est $1$-$\F$-connexe dans $\mc B_{\Sigma}(D)$, i.e.
$$   \partial\mc
B_{\Sigma}(D)\ucon\mc B_{\Sigma}(D)\,.$$
\end{prop}}

%\begin{equation}\label{v}
 %   \partial\mc B_{\Sigma}(D)\ucon\mc B_{\Sigma}(D)\,.
%\end{equation}

\noindent {Finalement, pour achever la preuve du th\'{e}or\`{e}me
(\ref{teoexistence}) il suffit de prouver (\ref{v}). Le reste de ce
chapitre est consacr\'{e} \`{a} cette d\'{e}monstration, mais il nous faut
pr\'{e}alablement introduire quelques constructions auxiliaires.}
%dans la section suivante.

\subsection{Constructions et descriptions auxiliaires}\label{3.2}
{Consid\'{e}rons l'application l'application
$$\sigma :\mc B^\sharp :=  \mc B_\Sigma (D)\longrightarrow D^{\sharp}$$ qui est \'{e}gale \`{a} ${\sigma}_{\mf M_{j}}$
en restriction \`{a} chaque voisinage $\mc B_{\Sigma}(\mf M_{j})$,
$j=k+1,\ldots ,n$ et \`{a} $\pi _D$ ailleurs. Elle satisfait toutes les
conditions des fibrations de Seifert, sauf peut-\^{e}tre la locale
trivialit\'{e} aux points de recollement, c'est \`{a} dire aux points de
$\Xi := (\partial D^*\cup\bigcup_{j=1}^k|\sigma _j|)$. En fait,
comme chaque fibre de $\sigma _D$ est \'{e}toil\'{e}e, $\sigma $ est une
pseudo-fibration (de Seifert) singuli\`{e}re dans $\Xi$ dans le sens
suivant :
\begin{defin}
Nous dirons qu'une submersion surjective entre deux vari\'{e}t\'{e}s
diff\'{e}rentiables $\pi : E\rightarrow B$ est une
\textbf{pseudo-fibration (de Seiffert) singuli\`{e}re dans $S\subset
B$}, s'il existe $E'\subset E$ tel que la restriction $\pi _{|E'} :
E'\rightarrow B$ est une fibration de Seiffert dont aucune fibre
singuli\`{e}re n'est contenues dans $S$, s'il existe une r\'{e}traction par
d\'{e}formation de $E$ sur $E'$ donn\'{e}e par une homotopie qui commute
avec $\pi$ et si de plus $E'$ et $E$ co\"{\i}ncident au dessus du
compl\'{e}mentaire d'un voisinage ouvert de $S$ dans $B$, que l'on peut
choisir arbitrairement petit.
\end{defin}
} \noindent{\textbf{A. Rev\^{e}tement ramifi\'{e} adapt\'{e} \`{a} la
pseudo-fibration.}} {Nous allons} d\'{e}finir un rev\^{e}tement ramifi\'{e} qui
trivialisera $\sigma $ au voisinage de la fibre singuli\`{e}re. Il est
bien connu que la base $\ssD$ de $\sigma$ admet une structure
d'orbifold ayant comme points de ramification $s_{k+1},\ldots,
s_{n}$ d'ordres respectifs $p_{k+1},\ldots,p_{n}$. D'o\`{u}
{l'isomorphisme}
$$\pi_{1}^{\mathrm{orb}}(\ssD)=\langle a_{0},\ldots,a_{n}
|\prod\limits_{i=0}^{n}a_{i}=1,\ a_{j}^{p_{j}}=1,\
j=k+1,\ldots,n\rangle\cong\pi_{1}(D^*)/{\mc K},$$ o\`{u} ${\mc K}$ est
le sous-groupe normal de $\pi_{1}(D^*)=\langle
a_{0},\ldots,a_{n}|\prod\limits_{i=0}^{n}a_{i}=1\rangle$ engendr\'{e}
par les \'{e}l\'{e}ments $a_{j}^{p_{j}}$, $j=k+1,\ldots,n$. Soit $\rho^*:\wh
D^*\to D^*$ le rev\^{e}tement associ\'{e} \`{a} ${\mc K}$. Il existe une
extension $\rho:\wh D\to \ssD$ de $\rho^*$ qui est un rev\^{e}tement
ramifi\'{e} d'ordre $p_{j}$ au dessus de $s_{j}$, pour $j=k+1,\ldots,n$.
Remarquons que $\wh D$ est simplement connexe car  $\rho^*$ d\'{e}roule
tous les lacets de $D^{*}$ sauf {les lacets} $a_{j}^{p_{j}}$ et par
construction ces derniers bordent des disques dans $\wh D$.
Consid\'{e}rons le {pull-back  $\wh {\mc B} := \wh D \times_{(\rho ,
\sigma )} \mc B^\sharp$ de la pseudo-fibration de Seifert} par
$\rho:\wh D\to \ssD$, d\'{e}fini par le diagramme cart\'{e}sien :
$$\begin{CD}
\wh{\mc B} @>\wh \rho>> \mc B^\sharp \\
@V\wh \sigma VV @V\sigma VV\\
\wh D @>\rho>> \ssD
\end{CD}$$
Visiblement $\wh \sigma:\wh{\mc{B}}\to \wh D$ n'a pas de fibres
singuli\`{e}re. On verra dans le paragraphe suivant que $\wh D$ est
contractile et que par cons\'{e}quence $\wh \sigma$ est  une
pseudo-fibration triviale, dans un sens clair.
\begin{center}
\begin{figure}[htp]
\hspace{1cm}
\includegraphics[width=7cm]{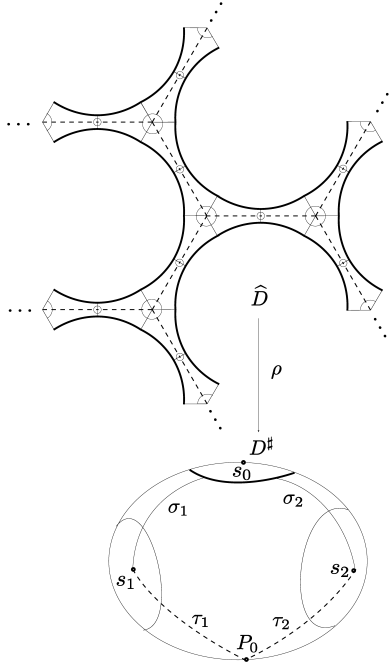}
\vspace{-0.5cm} \caption{Le rev\^{e}tement ramifi\'{e} $\rho:\wh D\to
D^{\sharp}$, pour $k=0$, $n=2$, $p_1 = 2$, $p_2 = 3$. {Les graphes $\mc G$ et $\wh{\mc{G}}$ sont en pointill\'{e}.}}
\end{figure}
\end{center}
\noindent{\textbf{B. Graphe adapt\'{e} au rev\^{e}tement ramifi\'{e} pr\'{e}c\'{e}dent.}}
 Afin d'unifier les notations nous \'{e}crirons $$\F^{\sharp}=\wt
\F_{|\mc B^{\sharp}}\quad \hbox{et}\quad \wh \F=\wh
\rho^*\F^{\sharp}\,.$$
Soient $P_{0}\in D^{\circ}$ et
$P_{j}\in|\sigma_{j}|\cap D^*$, $j=1,\ldots,k$. On notera aussi
$P_{j}=s_{j}$ si $j=k+1,\ldots,n$. Consid\'{e}rons des chemins r\'{e}guliers
simples (analytiques) $\tau_{j}:[0,1]\to \ssD$,
$j=-k,\ldots,-1,1,\ldots,n$ tels que, si on d\'{e}signe par $]\tau_{j}[$
l'image de $\tau_{j}$ priv\'{e}e de ses extr\'{e}mit\'{e}s, les propri\'{e}t\'{e}s
suivantes sont satifaites :
\begin{enumerate}[-]
\item $\tau_{j}(0)=P_{0}$ et $\tau_{j}(1)=P_{j}$;
\item $]\tau_{j}[\subset D^{\circ}$;
\item $]\tau_{j}[\cap ]\tau_{\ell}[=\emptyset$ si $j\neq \ell$;
\item $\mathrm{ind}_{s_{j}}(\tau_{\ell}\cdot\tau_{-\ell}^{-1})=\delta_{j\ell}$ si $1\le j,\ell\le k$.
\end{enumerate}
Consid\'{e}rons la structure de graphe combinatoire $\mc G$ d\'{e}finie sur
l'ensemble
$$ |\mc G|=\bigcup\limits_{j=-k}^{n}|\tau_{j|}
$$
en d\'{e}cr\'{e}tant que ses sommets sont
les points $P_{0},P_{1},\ldots, P_{n}$ et ses ar\^{e}tes les courbes $|\tau_{j}|$. %On peut d\'{e}finir une coloration sur $\mathcal G$ en d\'{e}cr\'{e}tant que le sommet $P_{0}$ est rouge et le reste de sommets sont bleus. De cette fa\c{c}on on obtient un graphe bi-color\'{e}e, c'est \`{a} dire, tel que deux sommets adjacentes sont toujours de couleurs diff\'{e}rentes.
{Nous laissons au lecteur le soin de prouver le lemme suivant, \`{a}
l'aide d' un champ de vecteur appropri\'{e} par exemple.
\begin{lema}\label{retract.par.def} Il existe une r\'{e}traction par d\'{e}formation $q:\ssD\to
|\mathcal{G}|$ telle que $q(|\sigma_{j}|)= \{P_{j}\}$, $j=1,\ldots
,n$ et telle que la restriction de $q$ \`{a} chaque composante connexe
de
  $$\partial D^\sharp \;\setminus\;\{\sigma _1(0),\ldots , \sigma _k(0),
  \sigma _1(1),\ldots , \sigma _n(1)\}$$ est injective.
\end{lema} \noindent
Relevons cette r\'{e}traction par $\rho $, c'est \`{a}
dire consid\'{e}rons l'unique r\'{e}traction par d\'{e}formation $\wh q : \wh D
\rightarrow |\wh{\mc G}|:=\rho^{-1}(|\mathcal G |)$ qui v\'{e}rifie la
relation de commutation $\rho \circ \wh q = \wh q \circ \rho $. Il
existe un unique graphe combinatoire $\wG$ dont $|\wG|$ est la
repr\'{e}sentation g\'{e}om\'{e}trique et telle que $\rho $ est la
repr\'{e}sentation g\'{e}om\'{e}trique d'un rev\^{e}tement de graphes $\wG
\rightarrow \mc G$.} Comme $\wh D$ est simplement connexe,
$\wh{\mc G}$ est un arbre et $\wh D$ est contractile.\\

\noindent {\textbf{C. Description des feuilles de $\wh \F$.} Appelons
\textbf{plaque de $\wh D$}} l'adh\'{e}rence {dans $\wh D$} de toute
composante connexe de l'image r\'{e}ciproque
$$\rho^{-1}\left(\ssD\setminus\bigcup_{j=1}^{n}|\sigma_{j}|\right)\subset \wh D.$$
Appelons \textbf{partie s\'{e}cable} d'une plaque $\mc M$ {de $\wh D$ }
toute composante connexe de l'ouvert $$\rho^{-1}\left(\bigcup_{j =
k+1}^n D_{s_{j}}\right) \cap \mc M.$$ Appelons \textbf{sous-plaque
{de la plaque $\mc M$},} {toute} diff\'{e}rence $\mc N := \mc M
\setminus \mc S$, ou $\mc S$ est une union non-vide de parties
s\'{e}cables de $\mc M$. La plaque $\mc M\supset \mc N$ sera dite
\textbf{plaque associ\'{e}e \`{a} la sous-plaque $\mc N$} et sera not\'{e}e $\mc
N^c$. Remarquons que l'on a :
\begin{equation}\label{aller.retour}
    \wh q (\mc N) = \wh q (\mc N^c)\,,\quad\wh q ^{-1}(\wh q (\mc N)) = \mc N^c\quad\hbox{et}\quad
    \mc N^c\cap \partial \wh D = \mc N \cap \partial \wh D\,.
\end{equation}
\noindent Appelons {maintenant} \textbf{plaque d'une feuille} $\wh
L$ de $\wh\F$, toute composante connexe $\mc M'$ de l'image
r\'{e}ciproque par $\wh \sigma$ d'une plaque $\mc M \subset \wh D$.
\begin{lema}\label{rho-diffeo-local}
La restriction de $\wh \sigma $ \`{a} une plaque de $\wh L$ est
injective et son image est une plaque ou bien une sous-plaque de
$\wh D$.
\end{lema}
\begin{proof}
L'holonomie autour du point d'attache d'une branche morte n'est
jamais l'identit\'{e}, i.e. $p_{j}>1$ pour tout $j=k+1,\ldots,n$, cf.
\cite{Mattei-Salem}.
\end{proof}
\noindent Si l'image par $\wh \sigma $ d'une plaque $ \mc M'$ de
$\wh L$ est une sous-plaque, nous dirons que $\mc M'$ est une
\textbf{petite plaque}.

Consid\'{e}rons une composante connexe $\Gamma $ de $\wh L \cap
(\rho\circ \wh \sigma )^{-1}(\partial D_{s_{j}})$, $j= k+1,\ldots
,n$. Une et une seule des deux \'{e}ventualit\'{e}s suivante est r\'{e}alis\'{e}e :
\begin{enumerate}
  \item $\Gamma $ est un cercle; alors $\Gamma \subset \inte{\wh L}$ borde un disque
  dans $\wh L$ et la restriction de $\wh \sigma $ \`{a} $\Gamma $ est
  injective,
  \item\label{ev2}  $\Gamma $ est un segment; alors $\Gamma $
  est contenu dans
  $\partial \wh L$, il existe des
  petites plaques $\mc N_1, \ldots \mc N_r$, $r< p_j$, deux \`{a} deux distinctes telles que :
  $\Gamma = \Gamma _1\cup \cdots\cup\Gamma _r$,
    $\Gamma _j \subset \partial\mc N_j$, $\sharp (\Gamma _j\cap \Gamma _{j+1}
    =1)$, $\Gamma _k \cap \Gamma _l = \emptyset$
    si $| k-l | \neq 1$, et  la restriction de $\wh \sigma $ \`{a} $\bigcup_{j=1}^r
    \mc N_j$ est injective.
\end{enumerate}
\noindent Dans l'evenualit\'{e} (\ref{ev2}), nous appellerons
\textbf{pi\`{e}ce associ\'{e}e \`{a} $\Gamma $} l'ensemble $$ \mc P(\Gamma ) :=
\overline{\bigcup_{j=1}^r \left(  \wh \sigma (\mc N_j)^c  \setminus
\wh \sigma (\mc N_j)\right)}
 $$
Remarquons que $\mc P(\Gamma )$ est un secteur du disque
$\rho^{-1}(D_{s_{j}})$ pour un $j\in \{k+1,\ldots , n\}$ appropri\'{e}.
Il est ainsi clair que $\wh \sigma (\Gamma )$ (qui est hom\'{e}omorphe \`{a}
$\Gamma $), est un r\'{e}tracte par d\'{e}formation de $ P(\Gamma )$.

Notons $(\Gamma _\alpha )_{\alpha \in \mathfrak{A}}$ la collection
des composantes $\Gamma $ du type (\ref{ev2}). Nous appelons
\textbf{compl\'{e}t\'{e} de la feuille $\wh L$} et {le quotient}
$$\wh L^c := \Big( \wh L \cup \bigsqcup_{\alpha \in \mathfrak{A}} \mc
P(\Gamma _\alpha )\Big)\Big/ \approx
$$
{o\`{u} $\bigsqcup$ est le symbole de l'union disjointe et $\approx$
est} la relation d'\'{e}quivalence suivante :
\begin{enumerate}
\item[(a)] Si $a\in \wh L$ et $b\in \mc P(\Gamma
_\alpha)$, alors $a\approx b$  ssi $a\in\Gamma _\alpha $ et $b = \wh
\rho (a)$
\item[(b)] $a\in \mc P(\Gamma
_\alpha )$ est \'{e}quivalent \`{a} $b\in \mc P(\Gamma _\beta  )$ ssi
$\alpha  = \beta $ et $a = b$.
\end{enumerate}

\begin{obs}
$\wh L$ se plonge proprement sur un ferm\'{e} de $\wh L^c$ qui est un
r\'{e}tract par d\'{e}formation de $\wh L^c$. Consid\'{e}rons donc d\'{e}sormais
$\wh L$ comme un sous-ensemble de $\wh L^c$. {La restriction {$\wh
\sigma_{\wh L}$} de} l'application $\wh \sigma$ {\`{a} $\wh L$} se
prolonge \`{a} $\wh L^c$ en une application, not\'{e}e {$\wh \sigma_{\wh
L}^c$}, dont l'image est l'union des plaques de $\wh D$ qui
intersectent $\wh \rho (\wh L)$. Ce prolongement est un
diff\'{e}omorphisme local : plus pr\'{e}cis\'{e}ment, sa restriction \`{a} chaque
composante connexe de l'image r\'{e}ciproque {$(\wh \sigma_{\wh
L}^{c}){}^{-1}(\mc M)$} d'une plaque $\mc M$ de $\wh D$ est un
diff\'{e}omorphisme sur $\mc M$. En effet, avec les notations de
{l'\'{e}ventualit\'{e}} (\ref{ev2}), on voit que la restriction de {$\wh
\sigma_{\wh L}^c$} \`{a} $\mc P(\Gamma )\cup \bigcup_{j=1}^r \mc N_j$
est un diff\'{e}omorphisme sur l'union des plaques $ \bigcup_{j=1}^r \wh
\sigma (\mc N_j)^c$.
\end{obs}

Nous appellerons \textbf{plaque de $\wh L^c$} toute composante
connexe de {$(\wh \sigma_{\wh L}^c){}^{-1}(\mc M)$}, o\`{u} $\mc M$ est
une plaque de $\wh D$.
\begin{center}
\begin{figure}[htp]
\includegraphics[width=6cm]{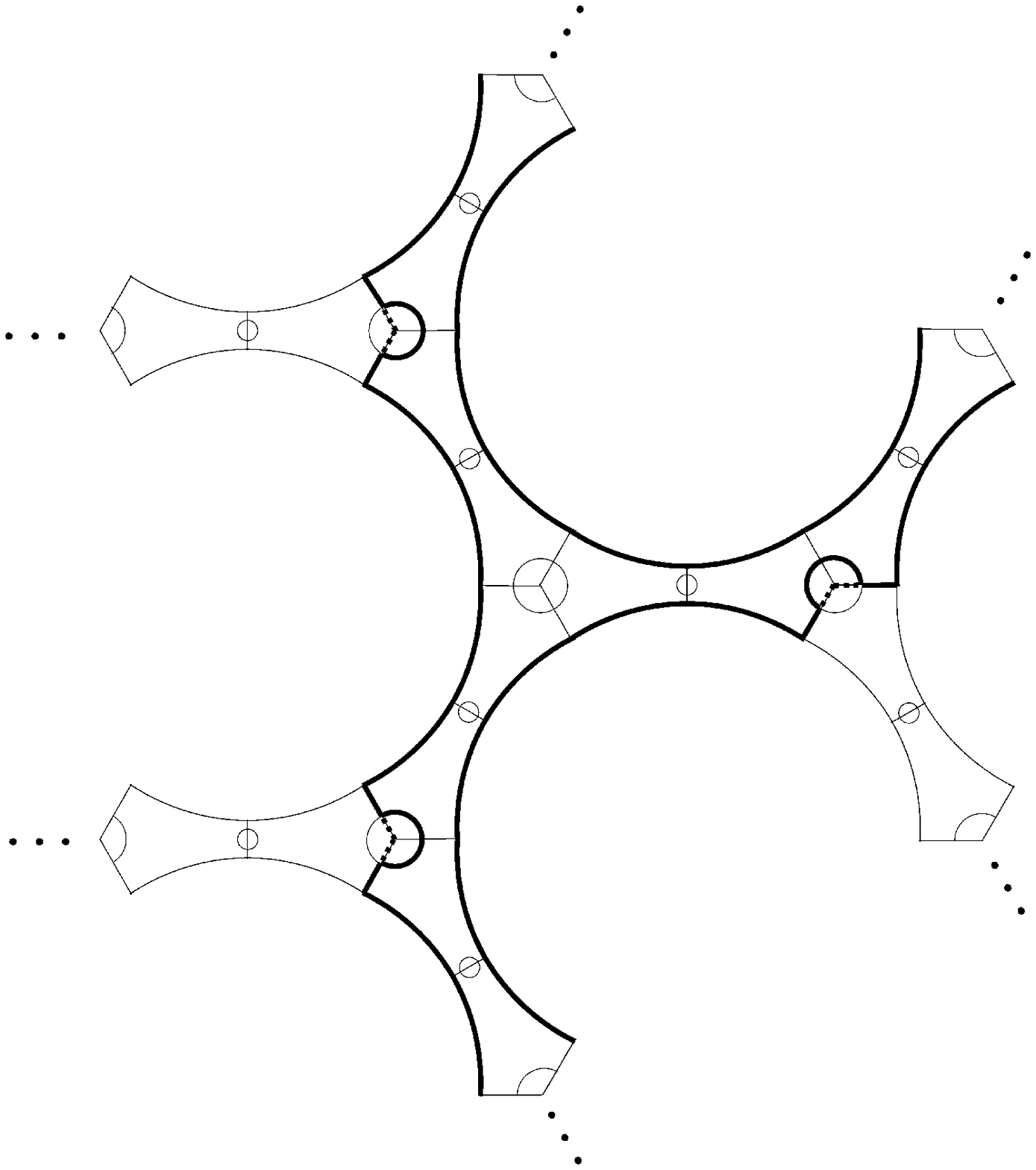}\includegraphics[width=5cm]{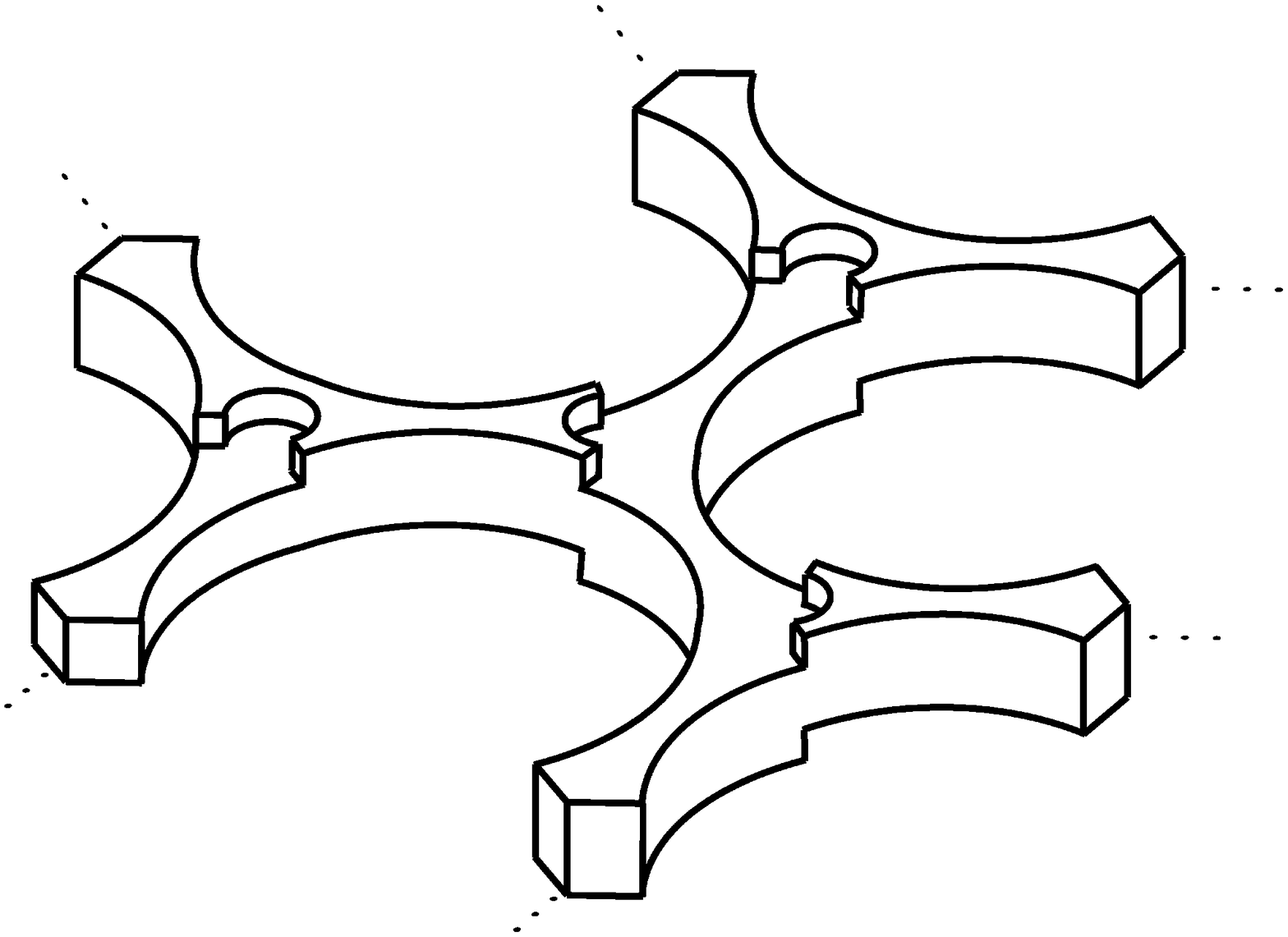}
\vspace{-3cm} {\caption{La feuille $\wh L$ et la feuille compl\'{e}t\'{e}e $\wh{L}^{c}$  vues comme
r\'{e}union de plaques. Image symbolique de $\wh{\mc{B}}^{\sharp}$.}}
\end{figure}
\end{center}
\vspace{-1em}
\begin{lema}\label{inj}
${\wh \sigma_{\wh L}^c} : \wh L^c \to \wh D$ est injective {et $\wh
L^c$ ainsi que $\wh L$ sont contractiles}.
\end{lema}
\begin{sublema}\label{injcomb}
Soit $h : \mc G\to \mc A$ un morphisme localement injectif d'un
graphe connexe $\mc G$ dans un arbre $\mc A$, i.e. la restriction de
$h$ \`{a} l' \'{e}toile $\rm{star}(s)$ de chaque sommet $s$ de $\mc G$ est
injective. Alors $\mc G$ est un arbre et $h$ est injective.
\end{sublema}
\begin{dem}
En consid\'{e}rant une exhaustion de $\mc A$ par des parties connexes
finies, on peut supposer que $\mc A$ est un arbre fini.  On raisonne
alors par r\'{e}currence sur le nombre de sommets de $\mc A$. Si $\mc A$
n'a qu'un seul sommet alors $\mc G$ aussi. Ainsi $\mc G$ est un
arbre et $h$ est injectif. Sinon, il existe un sommet d'extr\'{e}mit\'{e}
$s'$ de $\mc A$. Comme $h$ est localement injectif, tout $s\in
h^{-1}(s')$ est un sommet d'extr\'{e}mit\'{e} de $\mc G$. Si on enl\`{e}ve les
sommets d'extr\'{e}mit\'{e}s et leurs ar\^{e}tes adjacentes aux graphes $\mc A$
et $\mc G$, on obtient un nouveau morphisme localement injectif
$h':\mc G'\to\mc A'$, o\`{u} $\mc A'$ est un arbre avec un sommet de
moins que $\mc A$ et $\mc G'$ est toujours connexe.
\end{dem}
\begin{dem2}{du lemme (\ref{inj})}
Soit {$|\mc G_{\wh L}|$}la pr\'{e}-image par {$\wh \sigma\,_{\wh L}^c$}
de $|\wh{\mc G}|$. Comme {$\wh \sigma_{\wh L}^c$} est un
diff\'{e}o\-mor\-phis\-me local, cf. (\ref{rho-diffeo-local}), {il
existe un graphe combinatoire $\mc G_{\wh L}$ dont $|\mc G_{\wh L}|$
est la repr\'{e}sentation g\'{e}om\'{e}trique et un morphisme localement
injectif $ h: \mc G_{\wh L}\rightarrow \mc G$ dont la repr\'{e}sentation
g\'{e}om\'{e}trique $|h|$ est \'{e}gale \`{a} la restriction de $\wh \sigma_{\wh
L}^c$ \`{a} $|\mc G_{\wh L}|$.}  Or {$|\mc G_{\wh L}|$} est l'image par
la r\'{e}traction $\wh q$ de la feuille {compl\'{e}t\'{e}e} $\wh L^{c}$ qui
elle, est connexe. On en d\'{e}duit que {$\mc G_{\wh L}$} est un graphe
connexe. D'autre part, on a d\'{e}j\`{a} remarqu\'{e} que $\wh{\mc G}$ est un
arbre. Le sous-lemme (\ref{injcomb}) implique que $h$ est injectif
et { que $\mc G_{\wh L}$} est un arbre. {Ainsi $\wh L^c$, qui se
r\'{e}tracte sur $|\wh{\mc G}|$, est contractile et $|h|$ est injectif.
En utilisant de nouveau (\ref{rho-diffeo-local}), l'injectivit\'{e}
locale de $\wh \sigma_{\wh L}^c$ donne l'injectivit\'{e} de ${\wh
\sigma_{\wh L}^c}$. D'autre part la contractibilit\'{e} de $\wh L$
r\'{e}sulte directement de celle de $\wh L^c$, ce qui ach\`{e}ve la
d\'{e}monstration.}
\end{dem2}

\begin{lema}\label{retracte-bord}
La restriction de la r\'{e}traction $\wh q$ \`{a} toute composante connexe
$H$ de $\partial\wh D$ est injective et $H\cong \R$.
\end{lema}
\begin{dem}
{Il d\'{e}coule du lemme (\ref{retract.par.def})  que la restriction de
$\wh q$ \`{a} l'intersection de $H$ et d'une plaque de $\wh D$ est
injective. Ainsi il existe un graphe combinatoire $\mc H$ dont $H$
est la repr\'{e}sentation g\'{e}om\'{e}trique et un morphisme localement
injectif $\wh q_{\mc H} :\mc H\to \wh{\mc G}$ dont $\wh q$ est la
repr\'{e}sentation g\'{e}om\'{e}trique.} Le sous-lemme (\ref{injcomb}) implique
que {$\wh{q}_{\mc H}$} est injectif et $\mc H$ est un arbre.
{L'injectivit\'{e} de $\wh q$ en d\'{e}coule.} En fait $H$ est une copie de
$\R$ car $H$, \'{e}tant {une composante du} bord de la $2$-vari\'{e}t\'{e} $\wh
D$, est une $1$-vari\'{e}t\'{e} lisse.
\end{dem}

{
\begin{lema}\label{0-con}
Le bord de $\wh{\mc B}$ est $0$-$\wh\F$-connexe dans $\wh{\mc B}$, i.e.
$\partial\wh{\mc B}\zco{\wh\F}\wh{\mc B}$.
\end{lema}}
{
\begin{dem}
Consid\'{e}rons un chemin $\wh a:[0,1]\to \partial \wh{\mc B}$ dont les extr\'{e}mit\'{e}s sont dans une feuille $\wh L$ de $\wh\F$ et montrons qu'il existe un chemin $\wh c:[0,1]\to\partial \wh{\mc B}\cap\wh L$ de m\^{e}mes extr\'{e}mit\'{e}s que $\wh a$.
Par d\'{e}finition,
l'image  de {$|\wh a|$} par $\wh \sigma $ est contenue dans une
seule composante connexe de $\partial \wh D$. Celle-ci \'{e}tant une de
copie de $\R$ d'apr\`{e}s (\ref{retracte-bord}), consid\'{e}rons le chemin
g\'{e}od\'{e}sique $\check c $  joignant dans $\partial \wh D$ les
extr\'{e}mit\'{e}s de $\wh \sigma \circ a$. Le lemme (\ref{retracte-bord})
implique aussi que {$\wh q \circ\check c$ est un chemin g\'{e}od\'{e}sique}
de $|\wh{\mc G}|$. Ce chemin joint deux points de $\wh q\circ\wh
\sigma (\wh L)$. Cet ensemble est connexe  et, comme $\wh \sigma
(\wh L)$ est une union de plaques ou de sous- plaques, il est
d'apr\`{e}s (\ref{aller.retour}) la repr\'{e}sentation g\'{e}om\'{e}trique d'un
sous-arbre de $\wh{\mc G}$. Ainsi $|\wh q\circ \check{c}|$ est
n\'{e}cessairement contenu dans $\wh q\circ\wh \sigma (\wh L)$. Toujours
d'apr\`{e}s (\ref{aller.retour}) on a : $\wh q^{-1}(\wh q(\wh \sigma
(\wh L)))\cap \partial \wh D= \wh \sigma (\wh L)\cap \partial \wh
D$, cette \'{e}galit\'{e} \'{e}tant  vraie pour toute union de plaques et de
sous-plaques. Il vient l'inclusion : $\check c\subset \wh \sigma
(\wh L)$. D'apr\`{e}s le lemme (\ref{inj}), $\check c$ se rel\`{e}ve dans
$\wh L$ en un chemin $\wh c$ de m\^{e}mes extr\'{e}mit\'{e}s que $\wh a$.% et $\wh b$.
\end{dem}}

\vspace{1em}

\subsection{{Preuve de {la proposition} (\ref{v})}}\label{3.3} Pour
montrer la $1$-$\F^{\sharp}$-connexit\'{e} feuillet\'{e}e de $\partial \mc
B^{\sharp}$ dans $\mc B^{\sharp}$ {fixons une feuille $L^\sharp$ de
$\F^\sharp$} et consid\'{e}rons deux chemins {de m\^{e}mes extr\'{e}mit\'{e}s}
$$a:[0,1]\to\partial \mc B^{\sharp}\quad\textrm{et}\quad b:[0,1]\to L^{\sharp}\,,$$
{qui sont} homotopes dans $\mc B^{\sharp}$. Dans $L^{\sharp}$ on
peut homotoper $b$  \`{a} un chemin, qu'on continuera \`{a} noter $b$, ne
passant pas {par les disques conformes
$\sigma^{-1}(\overline{D}_{s_{j}})\cap \wh L^\sharp$},
$j=k+1,\ldots,n$,
$${|a|,|b|\subset B^*=B^{\sharp}\;\setminus\:\bigcup
\limits_{j=k+1}^{n}\sigma^{-1}(\overline{D}_{s_{j}})\,.}$$
{Comme la restriction $\wh \rho_{|\wh{\mc{B}}^*}  : \wh{\mc{B}}^* :=
\wh \rho ^{-1}(\mc B^*)\rightarrow \mc B^*$ est un rev\^{e}tement,} il
existe des $\wh \rho$-rel\`{e}vements $\wh a:[0,1]\to
\partial \wh{\mc B}$ et $\wh b:[0,1]\to \wh L$ {de $a$ et $b$,} de m\^{e}me origine. Le lacet $b^{-1}{\vee} a=\wh \rho(\wh{b}^{-1}{\vee}\wh a)$ est
homotope dans $\mc B^{\sharp}$ \`{a} un lacet constant  et sont image
par $\sigma$ est contenue dans {l'ensemble $D^*$, cf.
(\ref{defDstarDbecarre})}. D'apr\`{e}s (\ref{4.1.2.}) et (\ref{pi1B}) on
a les isomorphismes
$$\pi_{1}(\mc B^*)/\langle
c\rangle\cong\pi_{1}^{\mathrm{orb}}(\ssD)\cong\pi_{1}(D^*)/{\mc
K}\,.$$ {On dispose donc du diagramme commutatif}
$$\begin{CD}
\pi_{1}(\mc B^*)@>>>\pi_{1}(\mc B^{\sharp})\\
@VVV @VVV\\
\pi_{1}(D^*)@>>>\pi_{1}(D^*)/{\mc K}\end{CD}$$ {dont les fl\`{e}ches
horizontales sont le morphisme induit  par l'inclusion $\mc
B^*\subset \mc B^\sharp$ et le morphisme de passage au quotient, les
fl\`{e}ches verticales \'{e}tant induites par $\sigma $. On voit ainsi que
la classe d'homotopie de $\sigma(b^{-1}{\vee} a)$ appartient \`{a} ${\mc
K}$ et donc} $\wh\sigma\circ\wh a$ et $\wh\sigma\circ\wh b$
poss\`{e}dent les m\^{e}mes extr\'{e}mit\'{e}s.
Par construction (pull-back) la restriction de $\wh \rho$ \`{a} {toute}
fibre de $\wh \sigma$ est injective.
Il en r\'{e}sulte que $\wh a$ et $\wh b$ ont m\^{e}mes extr\'{e}mit\'{e}s.

{Il suffit maintenant d'appliquer le lemme \ref{0-con} pour trouver
un chemin $\wh c$ trac\'{e} dans $\wh L\cap\partial \wh B$ de m\^{e}mes
extr\'{e}mit\'{e}s que $\wh a$ et $\wh b$. En effet, comme d'apr\`{e}s le lemme
(\ref{inj}) $\wh L$ est contractile, les chemins $\wh c$ et $\wh b$
sont homotopes dans $\wh L$. Leurs images $c:= \wh \rho\circ \wh c$
et $b$ par $\wh \rho$ sont homotopes dans $L$. Ceci permet de
conclure, car l'existence d'une homotopie dans $\partial\mc
B^{\sharp}$ reliant $c$ \`{a} $a$, r\'{e}sulte de l'incompressibilit\'{e} de
$\partial\mc  B^{\sharp}$ dans $\mc B^{\sharp}$, d\'{e}j\`{a} d\'{e}montr\'{e}e \`{a}
l'\'{e}tape 3.}

\vspace{1em}

\section{{Blocs feuillet\'{e}s adaptables
attach\'{e}s \`{a} une singularit\'{e}}} \label{decompcollier}  Dans les
exemples (\ref{subsect.exemples}) nous avons donn\'{e} des mani\`{e}res de
construire des blocs feuillet\'{e}s adaptables en un point $s$ o\`{u} le
feuilletage est lin\'{e}arisable. Consid\'{e}rons maintenant le cas o\`{u} $\wF$
n'est pas lin\'{e}arisable en $s$. D'apr\`{e}s les hypoth\`{e}ses du th\'{e}or\`{e}me
principal, $\wF$ est un col r\'{e}sonnant. Pour simplifier l'\'{e}criture
nous notons $x:=x_s$ et $y:=y_s$ les coordonn\'{e}es d\'{e}finies en
(\ref{decompdiv}) et nous identifions le sous-ensemble $\{|x|<1,
\,|y|<1\}$ au polydisque $\D\times \D$ de $\C^2$, o\`{u} $\D$ d\'{e}signe le
disque unit\'{e} $\D(1) :=\{|z|<1\}$ de $\C$. Sur ce polydisque $\wF$
est d\'{e}fini par un champ de vecteur $X$ qui s'\'{e}crit :
$$X = x  \frac{\partial\phantom{x}}{\partial x} - y
\left(\lambda  + x y A(x,
y)\right)\frac{\partial\phantom{y}}{\partial y}\,,\quad \lambda =
\frac{p_0}{q_0}\in \mb Q_{>0}\,,\quad(p_0, q_0) = 1\,.$$ Nous allons
d'abord donner la construction de voisinages des axes  $\{xy=0\}$ de
type ``collier feuillet\'{e}s''. Ensuite, au paragraphe suivant, nous
montrerons le th\'{e}or\`{e}me (\ref{teoexistence}) d'existence de blocs
adapt\'{e}s pour $\alpha =s$, en admettant la proposition
(\ref{lema41}). Enfin les paragraphes
{(\ref{preuve.ex.reson})--(\ref{etapedulac})} seront consacr\'{e}s \`{a} la
preuve de cette proposition.

\subsection{Construction de colliers {et r\'{e}duction de la preuve de (\ref{teoexistence})}}\label{passgecols}
Pour $t \in \R i$, le flot $\phi_{t}^X $ de $X$ laisse invariant les
cylindres $\{|x| = c\} $, $c\in \R_{>0}$. De plus le diff\'{e}omorphisme
$\phi_{2\pi i}^X $ laisse invariant  chaque germe de droite $(\{x\}
\times \C, (x,0)) $ et, en restriction \`{a} ces droites, il est \'{e}gal au
germe d'holonomie de $\wF$ le long du lacet $\theta \mapsto e^{i
\theta} x\,$, $\theta \in [0, 2\pi ]$. Fixons $0 < \varepsilon_{3} <
1$ tel que $\phi_{t i}^{X}(x,y)$ appartient \`{a} $\overline{\D} \times
\overline{\D} $ pour tout $(x, y) \in \overline{\D} \times
\overline{\D}(\varepsilon_{3})$ et $ |t| \leq 2\pi$.\\

Dans tout le texte nous notons :
\begin{gather}\label{cercles}
   \quad C :=
   \partial  \D \times \{0\}\,,\quad \mc I_{\theta',\, \theta''}:=
   \{ (e^{i\theta}, 0) \;|\; \theta' \leq \theta
\leq \theta''\}\subset C\,,\\
\quad C' := \{0\}\times \partial \D, \quad \mc I'_{\theta',\,
\theta''}:=\{( 0, e^{i\theta})\;| \theta' \leq \theta \leq
\theta''\}\subset C'\,.\notag
\end{gather}

Soit $\Delta \subset \{ e^{i\theta_0}  \} \times
\D(\varepsilon_{3})$ et $ \theta_0<\theta_1\leq\theta_0+2\pi$.
L'ensemble :
$$U_\Delta  := \left\{\phi_{i t}^{X}(m)\;|\; m \in \Delta ,\; 0\leq t\leq
\theta_1 - \theta_0 \right\}\subset \mc
I_{\theta_0,\,\theta_1}\times \mathbb{D}\,,$$est appel\'{e}
\textbf{$\wF$-suspension de $\Delta$ au dessus de $\mc I_{\theta_0,
\theta_1}$ (suivant la projection $(x, y) \mapsto x $)}.

\begin{defin}\label{suspension}Nous dirons qu'un ensemble $U\subset
\partial \D\times \mathbb{D}$ est une $\wF$-\textbf{multi-suspension de
longueur $r$ au dessus de $\mc I_{\theta_0, \theta_r}$},
$\theta_r\leq \theta_0+2\pi $, s'il est \'{e}gal \`{a} une union
$\bigcup_{j=0}^{r-1} U_j$ de $\wF$-suspension $U_j$ au dessus de
$\mc I_{\theta_j, \theta_{j+1}}$, avec $\theta_0< \theta_1< \cdots
<\theta_r$
\end{defin}
\noindent Il est clair qu'une multi-suspension au dessus de $\mc
I_{\theta_0, \theta_r}$ est un un voisinage de $\mc I_{\theta_0,
\theta_r}$ dans $\mc I_{\theta_0, \theta_1} \times \D$, lorsque
$\Delta $ est un disque conforme ouvert contenant $( e^{i\theta_0},
0)$. D'autre part, pour $\varepsilon'_3<0$ assez petit, nous
d\'{e}finissons de la m\^{e}me mani\`{e}re la notion de \textbf{$\wF$-suspension
ou de $\wF$-multi-suspension au dessus de $\mc I_{\theta_0,
\theta_1}'$ suivant la projection $(x, y) \mapsto y $,} d'un
sous-ensemble $\Delta '$ de $\D(\varepsilon'_3)$.\\

\indent Consid\'{e}rons maintenant le champ
$$Y = x\left(\frac{p_0}{q_0} + x y A(x,
y)\right)^{-1}\frac{\partial\phantom{x}}{\partial x} - \, y
\frac{\partial\phantom{y}}{\partial y}\,.$$ Pour les temps $t$ r\'{e}el
son flot $\phi_{t}^{Y}$ laisse invariant les hyperplans r\'{e}els
$\{\arg y = c \}$, $c\in [0, 2\pi [$. On sait, cf. \cite{Loray} page
149, que si $\varepsilon_{4}
> 0$ est assez petit, il existe une (unique) fonction analytique
$\tau : \overline{\D}(\varepsilon_{4})^{*} \times \partial \D
\longrightarrow \R,$ telle que pour $0 < |x| \leq \varepsilon_{4}$
et pour $0 \leq \theta <2\pi$ on a :
\begin{enumerate}
\item[(i)] $\phi_{t}^{Y}(x,  e^{i \theta}) \in
\D^* \times \D^*$ pour tout $t \in ]0, \tau (x,  e^{i \theta})[\,$,
\item[(ii)] $\phi_{\tau(x,\;  \; e^{i
\theta})}^{Y}(x,  e^{i \theta}) \in \partial \D \times \D^\ast\,$,
\item[(iii)]\label{distzero} la distance entre
$\phi_{\tau(x, \, e^{i \theta})}^{Y} (x, \, e^{i \theta})$ et ${C}$
tend uniform\'{e}ment vers $0$ lorsque $x$ tend vers $0$.
\end{enumerate}

Fixons d\'{e}sormais $\varepsilon_{4} > 0 $ tel que ces propri\'{e}t\'{e}s
soient v\'{e}rifi\'{e}es. Donnons nous aussi est un disque conforme ferm\'{e}
$\Delta \subset \{ e^{i\theta_0}  \} \times \D(\varepsilon_{3})$ de
bord un lacet simple a. l. p. m. d'indice $1$ autour de l'origine.
Supposons que la taille  $\|\Delta\|_{y} $ d\'{e}finie en
(\ref{diametre}) est suffisamment petite pour que la suspension
$U_{\Delta}$ de $\Delta $ au dessus du cercle $C$ soit contenue dans
l'image de l'application suivante :
$$\Psi: \overline{\D}(\varepsilon_{4})^{\;\ast} \times
\partial \D\longrightarrow \partial \D \times
\overline{\D}^\ast\,,\quad \Psi(x,  e^{i \theta} ) := \phi_{\tau(x,
e^{i \theta})}^{Y}(x, e^{i \theta})\,.
$$

\begin{defin}\label{defin41} Le plongement $\Psi$ s'appellera
{\bf application passage du col dans le sens $x\rightharpoonup y$}
et l'ensemble : $$\mc{Col}(\Delta) := \{\phi_{t}^{Y}(m)\; /\; m \in
\overline{\D}(\varepsilon_{4})^\ast \times \partial \D,\; \Psi(m)
\in U_{\Delta},\; 0 \leq t \leq \tau(m)\}\,,$$ sera appel\'{e} \textbf{
$\wF$-collier de gabarit $\Delta\,$}.
\end{defin}

\noindent Il est clair que $\mc{Col}(\Delta)\cup Z$ est un voisinage
ferm\'{e} de $Z := \{xy = 0 \}\cap\mathbb{K}$ dans le polydisque ferm\'{e}
$\mathbb{K} := \overline{\D} \times \overline{\D}$. L'intersection
$\mc Col(\Delta )\cap
\partial \mathbb{K}$ est constitu\'{e} exactement des deux composantes connexes suivantes :
$$U_{\Delta}^{\ast}:= U_{\Delta} - C \quad
\mathrm{et}\quad V_{\Delta}^{\ast}:= \Psi^{-1}(U_{\Delta}^{\ast}).$$
Visiblement on a :
\begin{equation}\label{Vsigma}
    V_\Delta  := \overline{V_\Delta ^*} = V_\Delta ^*\cup
C'
\end{equation}
et cet ensemble est un voisinage de $C'$ dans
$\mathbb{D}(\varepsilon_4)\times \partial \D$. L'appellation
``collier'' est justifi\'{e}e par le fait que \textbf{l'application de
rectification}
\begin{equation}\label{retractcollier}
    R_\Delta {:} U_{\Delta}^{\ast} \times [0, 1] \longrightarrow
\mc{Col}(\Delta), \quad  R_\Delta (m, t):= \phi_{-t \cdot\tau \circ
\Psi^{-1}(m)}^{Y}(m)\,,
\end{equation}
\noindent est un diff\'{e}omorphisme analytique r\'{e}el qui conjugue le
feuilletage produit $\wF_{\mid U_{\Delta}^{\ast}} \times [0, 1]$ au
feuilletage $\wF_{\mid\mc{Col}(\Delta)} $. Visiblement $R_\Delta $
est \'{e}gale \`{a} l'identit\'{e} en restriction \`{a} $U_{\Delta}^{\ast} \times \{
0\}$ et \'{e}gale \`{a} $\Psi^{-1}$ en restriction \`{a} $U_{\Delta}^{\ast}
\times \{ 1\}$.
 Ainsi
l'on a trivialement les relations de connexit\'{e} feuillet\'{e}e :
\begin{equation}\label{equa02}
U_{\Delta}^{\ast} \; \ucon \; \mc{Col}(\Delta) \quad \mathrm{et}
\quad V_{\Delta}^{\ast} \; \ucon \; \mc{Col}(\Delta).
\end{equation}

\begin{prop}\label{lema41}
Il existe des constantes $c_0$, $c >0$ et une fonction $\mf d :
\R_{\geq 0}\rightarrow\R_{\geq 0}$, $ \lim_{r\rightarrow 0} \mf d(r)
= 0$, telles que, si $\|\Delta\|_{y}$ et $\mathbf{e}_{y}(\partial
\Delta)$ sont tous $\leq c_0$,  alors il existe un compact $\Delta'
\subset\D \times \{ 1\}$ v\'{e}rifiant les assertions suivantes :
\begin{enumerate}
\item le bord de $\Delta'$ est un lacet simple analytique lisse par
morceaux, d'indice $1$ par rapport \`{a} l'origine,
\item la suspension $U_{\Delta'}'$ de $\Delta'$ au dessus de $C'$ suivant la projection $(x,
y)\mapsto y$, est contenue dans $V_{\Delta}$,
\item $U_{\Delta'}^{'\ast} := ( U'_{\Delta'} -
C')$ est $1$-$\wF$-connexe dans $ V_{\Delta}^{\ast}$,
\item\label{propruglemme}  $\mathbf{e}_{x}(\Delta')\; \leq \;
\{\!\!\{ \mathbf{e}_{y}(\Delta) + \mf d(\|\Delta\|_{x}) \}\!\!\}$ et
$\|\Delta'\|_{x} \; \leq \; c \|\Delta\|_{x}^{\lambda }\,$.
\end{enumerate}
\end{prop}

%\subsection{Preuve du th\'{e}or\`{e}me  (\ref{teoexistence})}
\begin{proof}[{Preuve du th\'{e}or\`{e}me  (\ref{teoexistence}) dans le cas $\alpha=s\in Sing(\wF)$}]
Soit $\varepsilon>0$.
%Con\-si\-d\'{e}\-rons
{Trai\-tons} d'abord le cas %$\alpha = s\in Sing(\wF)$, avec
o\`{u} $\wF$ est
lin\'{e}arisable en $s$. Consid\'{e}rons les blocs  $\mc B_\alpha $ ou $\mc
B_{\eta,\kappa,\tau}$ que nous avons construit aux exemples
(\ref{subsect.exemples}). Dans le cas A  choisissons
\begin{equation}\label{bordexlin}
    \Sigma _1 = \{ x=x_0, \, |y|<\epsilon_1\} \quad \hbox{\rm et}
\quad\Sigma _2 =\{ y=y_0, \, |x|<\epsilon_2\}
\end{equation}
avec $\epsilon_1< \zeta/|x_0|$ et  $\epsilon_2 < \zeta/|y_0|$. Ces
blocs satisfont les assertions (1), (2) et (3) du th\'{e}or\`{e}me
(\ref{teoexistence}). Dans les trois cas A, B et C, leurs bords sont
des suspension d'ensembles $\Sigma _i$ du type (\ref{bordexlin})
ci-dessus. Pour prouver l'assertion (4) il suffit de voir que la
$F$-rugosit\'{e} et la taille des $\Sigma _i$ tend vers $0$ lorsque
$\zeta$ tend vers $0$. Cela d\'{e}coule de la proposition (\ref{prop31})
avec $g=F$, car $\be_y(\Sigma _1) = \be_x(\Sigma _2) =0$. Remarquons
aussi que les deux bords peuvent \^{e}tre {construits} de mani\`{e}re
compl\`{e}tement  ind\'{e}pendante; ceci rend trivial la preuve des
assertions (3') et (4').

Examinons maintenant le cas o\`{u} $\wF$ est r\'{e}sonnant non-lin\'{e}arisable
en $s$. Appliquons la proposition (\ref{lema41}) en choisissant
$\Sigma $ de $F$-rugosit\'{e} nulle. On obtient un bloc $\F$-adaptable
en posant $\mc B_\alpha := Col(\Sigma )\cup U_\Sigma ^*\cup
U'{}^*_{\Sigma '}$. L'assertion (1) du th\'{e}or\`{e}me (\ref{teoexistence})
r\'{e}sulte directement de la construction. Les proposition
(\ref{prop31}) et (\ref{cor.est.rug.}) permettent de traduire les
in\'{e}galit\'{e}s (4) de (\ref{lema41}) en terme de $F$-rugosit\'{e}. Il vient
:
$$
\be_F(\Sigma ')\leq K_1\| \Sigma \|_F\,,\quad \|\Sigma '\|\leq K_2
\|\Sigma \|_F^\lambda \,,
$$
pour des constantes $K_1$, $K_2$ appropri\'{e}es. Ainsi, avec les
notation du th\'{e}or\`{e}me (\ref{teoexistence}), $\be_F(U_\Sigma ^*)=0$ et
$\be_F(U'{}^*_{\Sigma '})$, $\|U_\Sigma ^*\|_F$, $\|U'{}^*_{\Sigma
'}\|_F$ tendent encore vers $0$ pour $\|\Sigma \|_F \rightarrow 0$;
ce qui prouve l'assertion (4). Les assertions (3') et (4') sont
encore une retranscription de la proposition (\ref{lema41}), en
posant $V_1 =V := U_\Sigma ^*$. Nous en laissons les d\'{e}tails au
lecteur.
\end{proof}

%\subsection{Description de $V_\Delta $}\label{preuve.ex.reson}

%\subsection{Op\'{e}ration de rabotage}\label{susect.rabotage}

%\subsection{Estimation de la rugosit\'{e}}\label{subs.estim.rug}

\subsection{Description de $V_\Delta$ {et du proc\'{e}d\'{e} de rabotage}}\label{preuve.ex.reson}
Notons $$h := \phi^X_{2\pi i}\mid_{1 \times \D(\varepsilon_{3})} :
\{ 1\} \times \D(\varepsilon_{3}) \longrightarrow \{ 1\} \times \C$$
l'application d'holonomie de $\F$ le long de $\partial \D \times \{
0\}$. Si la rugosit\'{e} de $\partial \Delta$ est assez petite, les
compacts $\Delta$ et $h(\Delta)$ sont \'{e}toil\'{e}s. Ils sont distincts
car $h$ ne poss\`{e}de pas de domaine invariant. La ``diff\'{e}rence
sym\'{e}trique''
$$\Delta \vartriangle h(\Delta ) :=
(\Delta - h(\overset{\scriptscriptstyle 0}{\Delta})) \cup (h(\Delta)
- \overset{\scriptscriptstyle 0}{\Delta})$$ est une union de
``lunules'' d\'{e}limit\'{e}es par des courbes simples. Plus pr\'{e}cis\'{e}ment
nous avons une subdivision de $[\theta_0, \theta_0+2\pi ]$ :
\begin{equation}\label{thetaj}
    \theta_{0} < \theta_1 < \cdots < \theta_{q}\; = \;\theta_{0}
+ 2\pi\,, \quad q \geq 1,
\end{equation}
\noindent et des courbes simples a.l.p.m. \`{a} valeur dans
$\partial\Delta\cup h(\partial \Delta)$,
$$\gamma_{j}(\theta) = \rho_{j}(\theta)\,e^{i\theta}, \quad
\widetilde{\gamma}_{j}(\theta)
=\widetilde{\rho}_{j}(\theta)\,e^{i\theta}, \quad 0 <
\rho_{j}(\theta) \leq \widetilde{\rho}_{j}(\theta),\quad \theta\in
I_{j} := [\theta_{j - 1}, \theta_{j}]\,,$$ avec $j = 1, \ldots, q$,
v\'{e}rifiant :
\begin{gather}\label{bordlunules}
\gamma _{j}(\theta_{j-1})= \wt \gamma _{j}(\theta_{j-1}) = \gamma
_{j-1}(\theta_{j})= \wt \gamma _{j-1}(\theta_{j })\,, \quad j = 1,
\ldots, q\,,\\ \gamma _{1}(\theta_0)= \wt \gamma _{1}(\theta_0) =
\gamma _q(\theta_q)= \wt \gamma _q(\theta_q)\,.,\notag
\end{gather}
et telles que $\Delta \vartriangle h(\Delta ) $ est l'union des
\textbf{lunules} :
\begin{equation}\label{lunules}
    \mc L_{j} = \{r\,e^{i\theta} \; / \; \theta \in I_{j}, \;
\rho_{j}(\theta) \leq r \leq \widetilde{\rho}_{j}(\theta)\}\,,
\qquad j=1,\ldots , q\,.
\end{equation}
\noindent Eventuellement $\mc L_j$ peut \^{e}tre r\'{e}duit \`{a} l'image
$\mathrm{Im}(\gamma _j)$ de $\gamma_{j}$, si $\gamma_{j} \equiv
\widetilde{\gamma}_{j}$, et dans ce cas \textbf{l'int\'{e}rieur de }
$\mc L_j$
$$\inte{\mc L}_j := \{r\,e^{i\theta} \; / \; \theta \in I_{j}, \;
\rho_{j}(\theta) < r < \widetilde{\rho}_{j}(\theta)\}$$ est vide.
Quitte \`{a} modifier la subdivision , nous supposons qu'il existe un
sous-ensemble d'indice $\mf K$ tel que :
$$\inte{\mc
L}_k \not= \emptyset \quad\hbox{pour}\quad k\in \mf K\quad
\hbox{et}\quad\gamma _j \equiv \wt \gamma _j\quad\hbox{pour}\quad
k\notin \mf K$$ Le bord de $\Delta$ est une concat\'{e}nation $\mu_{1}
\vee \cdots \vee \mu_q$ avec : $\mu_{j} = \gamma_{j}$ ou bien
$\mu_{j} = \widetilde{\gamma}_{j}$. Il en est de m\^{e}me du bord de
$h(\Delta)$. Comme l'holonomie ne poss\`{e}de pas de domaine invariant,
on ne peut pas avoir $\gamma_{j} \equiv \widetilde{\gamma}_{j} $
pour tout $j$ et il existe des lunules d'int\'{e}rieur non-vide.

\noindent L'holonomie $h$ v\'{e}rifie aussi $h(\Delta) \, \not\subset \,
\Delta $ et $\Delta \, \not\subset \, h(\Delta)$. Il en r\'{e}sulte
\begin{equation*}\label{nombrelunules}
\# \mf K \geq 2\,.
\end{equation*}
\noindent Nous pouvons d\'{e}composer le bord de $U_{\Delta}$ en l'union
disjointe  $$\partial U_{\Delta} = \mc L \cup
\partial^{\mathrm{sat}} U_{\Delta}\,$$ avec $\mc L := \bigcup_{k \in
\mf K} \mc L_{k}$ et
 $\partial^{\mathrm{sat}} U_{\Delta} \,$ d\'{e}signant l'ouvert
de $\partial U_{\Delta}$ constitu\'{e} des points $m$ satisfaisant la
\textbf{propri\'{e}t\'{e} de saturation} suivante :
\begin{enumerate}
\item[$(\star)$] \it Il existe un voisinage $W_{m}$ de $m$ et
$\varepsilon_{m} > 0$ tel que $\phi_{it}^{X}(m') \in \partial
U_{\Delta}$ pour tout $m' \; \in \;W_{m}\; \cap \;
 \partial U_{\Delta}$ et tout $|t| < \varepsilon_{m}\,$.
\end{enumerate}
\noindent Nous dirons que $m$ est un \textbf{point de saturation
g\'{e}n\'{e}rique de $U_{\Delta}$ pour $\mc F_{|U_\Delta }$}. Visiblement on
a:
$$\partial \mc L \,=\, \partial (\overline{\partial^{\mathrm{sat}}
 U_{\Delta}})\,=\, \bigcup_{k \in \mathfrak{K} }
 (\mathrm{Im}(\gamma _k) \cup \mathrm{Im}(\wt \gamma _k))
 \, ,$$
\noindent  Remarquons que les \textbf{lunules d\'{e}g\'{e}n\'{e}r\'{e}es} \'{e}point\'{e}es
de leurs extr\'{e}mit\'{e}s $$\overset{\diamond}{\mc L_j} := \{\gamma
_j(\theta)\,|\,\theta_{j-1}<\theta<\theta_j\}\,,\quad j\notin \mf
K\,,$$ sont contenues dans $\partial^{\mathrm{sat}} U_{\Delta}$.
D'autre part les points $m$ de $$\inte{\mc L}:= \bigcup_{j\in \mf
K}\inte{\mc L}_j\,$$ sont caract\'{e}ris\'{e}s par la propri\'{e}t\'{e} suivante :
\begin{enumerate}
\item[$(\star\star)$] \it L'une des assertions
suivantes est satisfaite pour $t > 0$ assez petit :
\begin{enumerate}
  \item $\phi_{it}^{X}(m) \; \in  \;
  \inte U_\Delta
\quad \mathrm{et} \quad \phi_{-it}^{X}(m) \; \not\in \; \inte
U_\Delta\,,$
  \item $\phi_{-it}^{X}(m) \; \in \;
\inte U_\Delta
 \quad
\mathrm{et} \quad \phi_{it}^{X}(m) \; \not\in \; \inte U_\Delta \,.$
\end{enumerate}
\end{enumerate}
\noindent Nous dirons que $m$ est un \textbf{point de non-saturation
g\'{e}n\'{e}rique de $U_{\Delta}$}.
\begin{obs}\label{une.eventualite}
Il est clair que pour chaque  $\inte{\mc L}_j$, $j\in \mf K $, une
seule des deux \'{e}ventualit\'{e}s $(a)$ ou $(b)$ est r\'{e}alis\'{e}e en tout
point.
\end{obs}

\indent Le diff\'{e}omorphisme  $\R$-analytique  $\Psi^{-1}$ d\'{e}finie en
(\ref{defin41})  conjugue les restrictions de $\mc F$ \`{a}
$\overline{\D}(\varepsilon_{4}) ^{\; \ast} \times \partial \D$ et \`{a}
$\partial \D \times \overline{\D} ^{\; \ast}$. Ainsi il  v\'{e}rifie
$$\Psi^{-1}(\partial
U_{\Delta}) = \partial V_{\Delta}\quad
\hbox{et}\quad\Psi^{-1}(\partial^{sat} U_{\Delta})=\partial^{sat}
V_{\Delta}\,,$$ o\`{u} $\partial^{\mathrm{sat}} V_{\Delta}$ est d\'{e}fini
par la propri\'{e}t\'{e} de saturation analogue \`{a} la propri\'{e}t\'{e} $(\star)$
mais relativement au champ $Y$ et \`{a} l'ensemble $\partial
V_{\Delta}$. Comme $\Psi$ laisse invariant les hyperplans $\{\arg y
= c\}$, $c\in \R_{>0}$, les  restriction de $\Psi^{-1}$ induisent
des diff\'{e}o\-morphismes
$$ U_\Delta^*  \,\cap\, (\partial \D\times\{ \arg y
\in I_k \})\overset{\sim}{\longrightarrow} V_\Delta^* \,\cap\, \mc T
(I_k)\,,$$ o\`{u} $\mc T (I_k) $ d\'{e}signe le ``3-tube'' $\mc T (I_k) :=
\overline{\D}(\varepsilon_{4}) ^{\; \ast} \times \{ e^{i\theta} \; /
\; \theta \in I_k\}$. Ainsi chaque image r\'{e}ciproque $ \Psi^{-1}(\mc
L_k)$, $k \in \mf K$, est contenue dans $\mc T (I_k)$. On a la
d\'{e}composition:
$$\partial V_{\Delta} =\mc L' \sqcup \partial^
{\mathrm{sat}} V_{\Delta}\,,\quad \hbox{avec}\quad \mc L' :=
\bigcup_{k \in \mf K} \mc L'_{k}\,,
 \quad \mc L'_{k} := \Psi^{-1}(\mc
L_{k})\,.$$ De m\^{e}me les points int\'{e}rieurs de $\inte{\mc L'} =
\Psi^{-1}(\inte{\mc L})$ sont caract\'{e}ris\'{e}s par l'analogue pour $Y$
et $\partial
V_{\Delta}$ de la propri\'{e}t\'{e} $(\star\star)$.\\

\indent Pour tout sous-ensemble $\Omega $ de $\C\times
\partial \D\,$, ou de $\partial \D\times \C \,$, pour  $
I\subset \R$, et pour $\theta \in \R$, nous adoptons maintenant les
notations suivantes :
\begin{equation}\label{tranches}
\Omega (I) := \Omega \cap \{\arg y \in I\}\,,\quad \Omega (\theta)
:= \Omega  \cap \{\arg y = \theta \}\,.
\end{equation}
\noindent Fixons $1\leq k\leq q$ et consid\'{e}rons maintenant les
restrictions de $\F$ aux ensembles $ V_\Delta(I_k)$. Visiblement on
a :
$$\Big( \partial \left( V_\Delta (I_k)\right)\Big)(\inte I_k) =
(\partial V_\Delta ) (\inte I_k) = \Big(\Psi^{-1} (\partial
U_\Delta) \Big)(\inte I_k) = \Psi^{-1} \Big((\partial U_\Delta)
(\inte I_k)\Big)$$ \indent Supposons d'abord $k\notin \mf K$. On a
l'\'{e}galit\'{e} $(\partial U_\Delta) (\inte I_k)=(\partial^{sat} U_\Delta)
(\inte I_k)$. Comme $\Psi$ pr\'{e}serve la propri\'{e}t\'{e} de saturation
g\'{e}n\'{e}rique, les \'{e}galit\'{e}s ci-dessus donnent :
$$(\partial V_\Delta) (\inte I_k)=(\partial^{sat} V_\Delta) (\inte
I_k)\,,\quad k\notin \mf K\,.$$ \noindent Ainsi dans ce cas, chaque
feuille est une courbe simple dont l'une des extr\'{e}\-mi\-t\'{e}s est
situ\'{e}e sur $V_{\Delta}(\theta_{j - 1})$, l'autre extr\'{e}mit\'{e} \'{e}tant
situ\'{e}e sur $V_{\Delta}(\theta_{j})$. D'o\`{u} :

\begin{obs}\label{susK} L'ensemble $V_{\Delta}(I_k)$ est une suspension
au dessus de $I_k$, pour $k\notin \mf K$.\end{obs}

\indent Supposons maintenant  $k \in \mf K$. Le m\^{e}me raisonnement
montre que les extr\'{e}mit\'{e}s des feuilles de $\mc
F|_{V_{\Delta}(I_{k})}$ sont situ\'{e}es sur $\mc L'_{j}\cup
V_{\Delta}(\theta_{k - 1})\cup V_{\Delta}(\theta_{k})$.

\begin{lema}\label{sublema41} Il n'existe pas de feuille de $\mc
F|_{V_{\Delta}(I_{k})} $ dont les deux extr\'{e}mit\'{e}s sont situ\'{e}es sur
${\inte{\mc L'}_k}$.
\end{lema}

\begin{dem} Raisonnons par contrapos\'{e}e. Une feuille $L' \in \mc
F|_{V_{\Delta}(I_{k})}$ telle que $\partial L' \subset {\inte{\mc
L'}_k}$ est l'image inverse par $\Psi$ d'une partie connexe d'une
feuille de $\F_{|U_{\Delta}}$, qui est param\'{e}tris\'{e}e par $t \mapsto
\phi_{it}^{X}(m)$, $t \in [0, 2\pi]$ et telle que $m$ et
$\phi_{i2\pi}(m)$ appartienne \`{a} $\mc L_{k}\cap\{\mathrm{arg} \,y\neq
\theta_k, \theta_{k + 1}\}$. Mais ceci est impossible d'apr\`{e}s
(\ref{une.eventualite}).
\end{dem}
\begin{obs}\label{extremites} Chaque feuille de
$\mc F|_{V_{\Delta}(I_{k})}, \; k \in \mf K$, dont l'une des
extr\'{e}\-mit\'{e}s $m$ est situ\'{e}e sur ${\inte{\mc L'}_k}$ a son autre
extr\'{e}\-mit\'{e}, que nous notons $l_{k}(m)$, situ\'{e}e sur $V(\theta_k)$ ou
sur $ V(\theta_{k-1})$. Par un argument de continuit\'{e} on voit qu'il
existe $\theta_{+}(k)\in \{\theta_{k - 1}, \theta_k\}$ tel que
$l_{k}(m)\in V_{\Delta}(\theta_{+}(k))$, pour tout $m\in {\inte{\mc
L'}_k}$. Par extension on obtient une application analytique r\'{e}elle
$l_{k}: \mc L'_{k} \longrightarrow V(\theta_{+}(k)), \quad m
\longmapsto l_{k}(m)$. Toutes feuilles de $\mc
F|_{V_{\Delta}(I_{k})}$ qui n'intersectent pas $\mc L'_k$ ont une de
leurs extr\'{e}mit\'{e}s situ\'{e}e sur $V(\theta_k)$ et l'autre sur $
V(\theta_{k-1})$. Finalement on obtient une application
$\R$-analytique `` de transport holonome'' que nous notons encore
$l_k$
$$
l_k : V_\Delta (I_k) \longrightarrow V_\Delta (\theta_+(k))\,,\quad
\{l_k(m)\}:= L_k(m)\cap V_\Delta (\theta_+(k))\,,
$$
avec $L_k(m)$ d\'{e}signant la feuille de $\wF_{|V_\Delta (I_k)}$ qui
contient le point $m$.
\end{obs}

%\subsection{Op\'{e}ration de rabotage}\label{susect.rabotage}
Consid\'{e}rons un  intervalle $I =[\theta' , \theta'' ]$, $0< \theta''
-\theta' \leq 2\pi$ et un sous-ensemble connexe ferm\'{e} $\Omega$ de $
\mathbb{D}\times \partial \D$ dont l'image par la projection $(x,
y)\mapsto y$ est l'arc de cercle $\mathcal{I}:=\{ (e^{i\theta}, 0)
\;|\;\theta\in I \}$. Nous supposons que $\F_{|\Omega }$ est un
feuilletage  en courbes lisses r\'{e}elles. Consid\'{e}rons le sous-ensemble
$\Omega ^0\subset \Omega $ des points $m$ qui sont int\'{e}rieurs \`{a}
$\Omega $ pour la topologie feuillet\'{e}e, i.e. il existe $\epsilon_m
>0$ tel que $\phi_{it}^Y(m)\in \Omega $ pour $- \epsilon_m \leq t \leq
\epsilon_m$. Notons $$\partial_\F\Omega := \Omega \setminus \Omega
^0\quad\hbox{et}\quad \Omega ^\vartriangle := \Omega \setminus
(\Omega(\theta')\cup \Omega(\theta''))= \Omega (\inte I)\,.$$
\begin{defin}\label{rabote}
Nous appelons \textbf{rabot\'{e} de $\Omega $} au dessus de $I$
l'ensemble :
$$
\mathrm{Rab}_I(\Omega ) := \overline{\Omega ^\vartriangle
\;\setminus\; \mathrm{Sat}\Big(\,
\partial_\F\Omega  \setminus (\Omega(\theta')\cup \Omega(\theta''))
\,,\,\; \Omega ^\vartriangle}\Big)\,.
$$
\end{defin}
\noindent Visiblement $\mathrm{Rab}_I(\Omega )$ est visiblement une
suspension (cf. {(\ref{suspension})}) au dessus de $\mathcal{I}$
pour la projection $(x, y)\mapsto y$, si $\theta'' -\theta' < 2\pi$.
C'est en fait le plus grand\footnote{Il est clair que l'union de
deux ensemble de type suspension au dessus de $\mc I$ est aussi de
type suspension au dessus de $\mc I$.} sous-ensemble de $\Omega $ de
type suspension au dessus de $\mathcal{I}$. Pour $\theta'' -\theta'
= 2\pi$, c'est une suspension au dessus du cercle
$C'$, point\'{e}e en $P_0:=(0, e^{i\theta'})$. \\

D'apr\`{e}s (\ref{extremites}) pour $k\in \mathfrak{K}$ on a l'\'{e}galit\'{e}
$$\mathrm{Rab}_{I_k}(V_{\Delta}(I_{k}))= \overline{V_{\Delta}(\inte I_{k}) -
\mathrm{Sat}_\F(\mc L'_{k}(\inte I_{k}) ;\, V_{\Delta}(\inte
I_{k}))}\,.$$ \noindent Cet ensemble est la suspension de $V_\Delta
(\theta_-(k))$ au dessus de $ \mathcal{I_k}:= \{0\}\times \{
e^{i\theta}\;/\;\theta\in I_k\}$, o\`{u} $\theta_-(k)$ est l'\'{e}l\'{e}ment de
$\{\theta_{k-1}, \theta_k\}$ autre que $\theta_+(k)$. D'autre part,
pour $k\notin \mathfrak{K}$, $V_\Delta (I_k)$ est d\'{e}j\`{a} une
suspension et $\mathrm{Rab}_{I_k}(V_{\Delta}(I_{k})) =
V_{\Delta}(I_{k})$.

%\vspace{1cm}

\begin{center}
\begin{figure}[htb!]
\includegraphics[width=12cm]{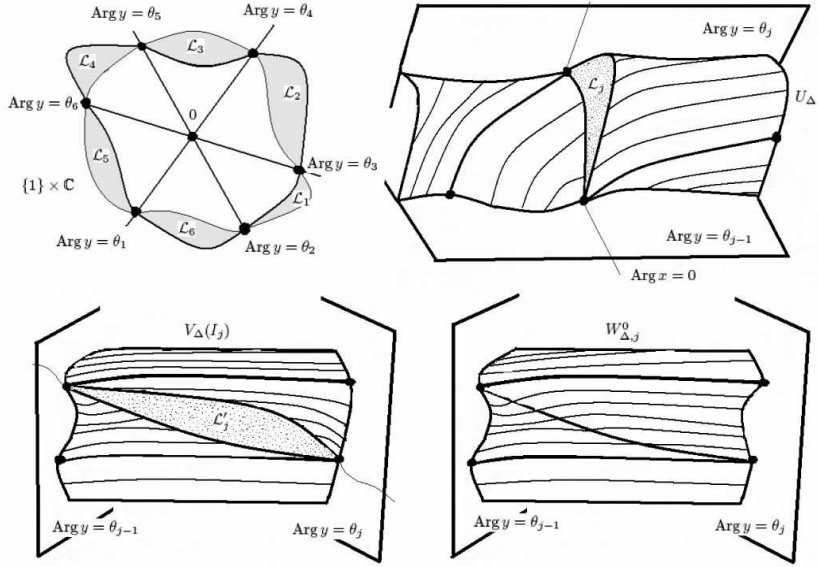}%, height=3cm,
\vspace{-0.75cm} {\caption{Description des lunules et du proc\'{e}d\'{e} de
rabotage.}}
\end{figure}
\end{center}

%\vspace{0.4cm}

%\begin{center}
%\begin{figure}%[]
%\includegraphics[width=5cm, height=3cm,]{jffigure2}
%\includegraphics[width=2cm, height=2cm,]{texte}
%\includegraphics[width=5cm, height=3cm,]{jffigure3}
%\includegraphics[width=8cm]{dessin_tranche_1}
%\vspace{-1cm} %\caption{Chirurgie d'homotopies.}
%\end{figure}
%\end{center}

%\vspace{1cm}

Ainsi l'ensemble
$$W_\Delta^{0} := \bigcup_{k=1}^q
W_{\Delta,\, k}^{0} \,,\quad \hbox{avec}\quad W_{\Delta,\, k}^{0} :=
\mathrm{Rab}_{I_k}(V_{\Delta}(I_{k}))$$ est une multi-suspension de
longueur $\leq q$ au dessus du cercle $C'$, cf. (\ref{suspension}).
Notons que d'apr\`{e}s la remarque (\ref{extremites}) on a :
\begin{equation}\label{paroies.rabotes}
    \left\{W_{\Delta ,\, k}^0(\theta_{k-1}),\, W_{\Delta ,\,
k}^0(\theta_{k})\right\} = \Big\{V_\Delta (\theta_-(k)),\; l_k(
V_\Delta (\theta_-(k)) ) \Big\}\,.
\end{equation}
Supposons $\varepsilon_4>0 $ assez petit pour que les
diff\'{e}omorphismes d'holonomies $\phi_{2\pi
 i}^X$ et $\phi_{2\pi
 i}^Y$ ne poss\`{e}dent pas  point p\'{e}riodique dans $U_\Delta ^*$ et
 $V_\Delta ^*$ respectivement. L'ensemble $W_\Delta^{0\,*} :=
 W_\Delta^{0}
\setminus C' $ est $0$-connexe dans $V_{\Delta}^*$, car chaque
feuille de la restriction de $\F$ \`{a} $\mathrm{Sat}_\F(\mc L'_{k};
V_{\Delta}(I_{k}))$ contient une extr\'{e}mit\'{e} d'une feuille de
$\F_{|V_\Delta^* }$. En  appliquant la proposition (\ref{prop3}) on
obtient : $\label{uconetape1} W_\Delta^{0\,*} \;\ucon\;
V_\Delta^*\,$.\\

\indent Nous allons maintenant construire une multi-suspension
$W_{\Delta}^{(1)\,*}\subset W_{\Delta}^{0\,*}$ de longueur $ \leq q
- 1$ au dessus ed $C'$, qui sera $1$-$\wF$-connexe dans
$W_{\Delta}^{0\,*}$. Consid\'{e}rons la restriction de $\mc F$ \`{a} $
W_{\Delta, \,q - 1}^{0} \cup W_{\Delta,\,q}^{0}$. Les feuilles sont
des courbes hom\'{e}omorphes \`{a} un intervalle ferm\'{e} dont les extr\'{e}mit\'{e}s
sont : ou bien toutes deux situ\'{e}es sur $W_{\Delta,\, q -
1}^0(\theta_{q-2})\cup  W_{\Delta,\, q}^0(\theta_{q})$, ou bien
l'une est situ\'{e}e sur $W_{\Delta,\, q - 1}^{0}(\theta_{q-2})\cup
W_{\Delta,\, q}^{0}(\theta_{q})$ et l'autre sur la diff\'{e}rence
sym\'{e}trique $ W_{\Delta,\, q - 1}^0(\theta_{q-1}) \vartriangle
W_{\Delta,\, q}^0(\theta_{q-1})$. Ainsi on a :
$$\mathrm{Rab}_{I_{q - 1} \cup
I_{q}}(W_{\Delta, q - 1}^0 \cup W_{\Delta,q}^0) =
\mathrm{Sat}\Big( W_{\Delta,\, q - 1}^0(\theta_{q-1})\cap
W_{\Delta, q}^0(\theta_{q-1}),  W_{\Delta, q - 1}^0 \cup
W_{\Delta,q}^0\Big).$$ Cet ensemble que nous notons $
W_{\Delta,q - 1}^1$ est une suspension au dessus de $ I_{q -
1}^{1} := I_{q - 1} \cup I_{q}$. L'ensemble
$$W_{\Delta}^1 := \overline{W_{\Delta}^0 - (W_{\Delta,\,q - 1}^0
\cup W_{\Delta,\,q}^0)} \cup W_{\Delta,\, q - 1}^{1}$$ est une
multi-suspension de longueur $\leq q - 1$. Par le m\^{e}me argument que
pr\'{e}c\'{e}demment la 0-$\wF$-connexit\'{e}e de $W_{\Delta}^{1\,*} :=
W_{\Delta}^1\setminus C'$ dans $W_{\Delta}^{0\,*}$ est \'{e}vidente et
il vient :  $W_\Delta^{1\,*} \;\ucon\; W_\Delta^{0\,*}$.

En it\'{e}rant ce proc\'{e}d\'{e} on obtient une succession de multi-suspension
$W_{\Delta}^j \; \subset \; W_{\Delta}^{j - 1}$ de longueur $\leq q
- j$ v\'{e}rifiant :
$$
W_{\Delta}^{q-1} \; \ucon \cdots \ucon\; W_{\Delta}^j\; \ucon \cdots
\ucon W_{\Delta}^0 \;\ucon\; V_\Delta ^*\,.
$$
Ainsi $W_{\Delta}^{q-1}$ est une suspension au dessus du cercle $C'
$ point\'{e}e en $P_0:=(0, e^{i\theta_0})$, d'un ensemble
$\widetilde{\Delta} '\subset \mathbb{D}\times \{P_0\}$ qui v\'{e}rifie
les propri\'{e}t\'{e}s $(1)$, $(2)$ et $(3)$ de la proposition
(\ref{lema41}). Si $\theta_0 = 0$, nous posons $\Delta
':=\widetilde{\Delta} '$ et $U'_{\Delta '} := W_\Delta ^{q-1}$. Si
$\theta_0$ n'est pas nul, pour obtenir une suspension point\'{e}e en
$(0, 1)$, nous effectuons une op\'{e}ration de rabotage suppl\'{e}mentaire :
nous posons $J:=[0, 2\pi ]$, l'ensemble :
$$
U'_{\Delta '} := \mathrm{Rab}_J\big(W_\Delta ^{q-1}\big)\;\ucon\;
W_{\Delta}^{q-1}\;\ucon\; V_\Delta ^*\,.
$$
convient.\\

Il reste maintenant \`{a} prouver l'assertion (4). Pour cela nous allons
estimer les ``pertes de
rugosit\'{e}'' \`{a} chaque \'{e}tape de l'induction.

{\subsection{Estimation de la rugosit\'{e}}\mbox{}\\
\noindent\textbf{Etape 1 : perte de rugosit\'{e} \corjf{par op\'{e}rations
de} rabotage.} Reprenons les notations (\ref{tranches}) du
paragraphe pr\'{e}c\'{e}dent. Pour $I=[\theta',\theta'']$ et pour $\Omega
\subset \D\times \{e^{i\theta},\ \theta'\le\theta\le\theta''\}$},
posons :
$$ \be_{I,x}(\Omega ) := \max\{ \be_x\left(\partial (\Omega
(\theta'))\right) , \, \be_x(\partial (\Omega (\theta''))\}\,,
$$
$$
\II\Omega \II_{I,x} := \max\{ \|\partial (\Omega (\theta'))\|_x , \,
\|
\partial (\Omega (\theta''))\|_x \}\,,
$$
en convenant que $\be_{I,x}(\partial (\Omega (\theta))) = +\infty$,
si $\partial (\Omega (\theta)))$ n'est pas un chemin a. l. p. m..\\

D'apr\`{e}s (\ref{paroies.rabotes}) on a pour tout $k =1,\ldots ,
q$\begin{gather} {\be}_{I_k, x} (W_{\Delta,\, k}^0) = \max\Big\{
\be_x ( V_\Delta (\theta_-(k))),\; \be_x \left( l_k( V_\Delta
(\theta_-(k)) \right) )\notag
\Big\} \,, \\
\II W_{\Delta,\, k}^0 \II_{I_k, x} = \max\Big\{ \II V_\Delta
(\theta_-(k))\II_x,\; \II ( l_k( V_\Delta (\theta_-(k)) )
\II_x\Big\}\,.\notag
\end{gather}
La proposition (\ref{prop31}) appliqu\'{e}e \`{a} la restriction du
transport holonome $l_k|_{V_\Delta (\theta_-(k))} : V_\Delta
(\theta_-(k)) \rightarrow V_\Delta (\theta_+(k))$ donne :
$$
\be _{I_k,x} ( W_{\Delta,\, k}^0 ) \leq  \{\!\!\{ \be_x ( V_\Delta
(\theta_-(k))) + c_1^{(0)}
 \|  V_\Delta (\theta_-(k)) \|_x \}\!\!\}
\,,$$
$$ \II  W_{\Delta,\, k}^0 \II_{I_k,x}  \leq c_2^{(0)} \|
V_\Delta (\theta_-(k))\|_x\,,
$$
pour des des constantes $c_1^{(0)}$, $c_2^{(0)}>0$ appropri\'{e}es. Il
vient :
\begin{equation}\label{estirug0}
     \be _{I_k,x} ( W_{\Delta,\, k}^0 )
 \leq  \{\!\!\{ \max_{j=1,\ldots , q} \be_x ( V_\Delta
(\theta(j))) +  c_1^{(0)}\max_{j=1,\ldots , q} \| V_\Delta
(\theta_j)\|_x \}\!\!\} \,,
\end{equation}
\begin{equation*}
\II  W_{\Delta,\, k}^0 \II_{I_k,x} \leq c_2^{(0)} \max_{j=1,\ldots ,
q}\| V_\Delta (\theta_j)\|_x\,.
\end{equation*}
\noindent Appliquons de m\^{e}me les in\'{e}galit\'{e}s (\ref{rugsom}) et la
proposition (\ref{prop31}) aux applications de transport holonome  :
$$W^0_{\Delta,\, q - 1}(\theta_{q-1})\cap W^0_{\Delta,\,
q}(\theta_{q-1})\rightarrow W^0_{\Delta, q-1 }(\theta_{q-2})\,,
$$
$$
W^0_{\Delta,\, q - 1}(\theta_{q-1})\cap W^0_{\Delta,\,
q}(\theta_{q-1})\rightarrow W^0_{\Delta, q }(\theta_{q})\,.$$ On
obtient ais\'{e}ment l'in\'{e}galit\'{e}
\begin{multline} \be_{I_{q-1}^{1},x}(
W^1_{\Delta,\, q - 1})\leq \{\!\!\{  \max\{ \be_x(W^0_{\Delta ,
q-1}(\theta_{q-1})), \, \be_x( W^0_{\Delta , q}(\theta_{q-1}))\} + \\
\wt c_1{}^{\!(1)} \max\{ \| W^0_{\Delta , q-1}(\theta_{q-1}) \|_x,
\, \| W^0_{\Delta , q}(\theta_{q-1}) \|_x \}\,\}\!\!\}\,.
\end{multline}
Ce qui donne grace \`{a} (\ref{estirug0})
$$
\be _{I_{q-1}^{1},x} ( W^1_{\Delta,\, k})
 \leq  \{\!\!\{ \max_{j=1,\ldots , q} \be_x ( V_\Delta
(\theta(j))) + c_1^{(1)}\max_{j=1,\ldots , q} \| V_\Delta
(\theta_j)\|_x \,\}\!\!\}\,,
$$
pour une constantes $c_1^{(1)}>0$ appropri\'{e}e. On montre de m\^{e}me une
in\'{e}galit\'{e} :
$$\II  W_{\Delta,\, k}^{}
\II_{I_{q-1}^{1},x}  \leq c_2^{(2)} \max_{j=1,\ldots , q}\| V_\Delta
(\theta_j)\|_x \,.$$ En it\'{e}rant ces majorations tout au long de
l'induction pr\'{e}c\'{e}dente, on aboutit \`{a} des in\'{e}galit\'{e}s similaires pour
$U'_{\Delta '}$ :
\begin{gather}\label{estirugq}
    \be _{I_k^{q-1}, x} (U'_{\Delta '} )
 \leq \{\!\!\{ \max_{j=1,\ldots , q}  \be_x ( V_\Delta
(\theta(j))) + c_1^{(q-1)} \max_{j=1,\ldots , q} \| V_\Delta
(\theta_j)\|_x \,\}\!\!\}\,,\notag\\
 \II U'_{\Delta '} \II_{I_k^{q-1}, x}  \leq c_2^{(q-1)}
\max_{j=1,\ldots , q} \| V_\Delta (\theta_j)\|_x \,.\notag
\end{gather}
\noindent Pour obtenir les majorations (\ref{propruglemme}) de la
proposition (\ref{lema41}) il suffit de prouver l'existence de
constantes $c_j$ et de fonctions $\mf d_j : \R_+\rightarrow\R_+$,
$\lim_{r\rightarrow 0}\mf d_j(r)= 0$ telles que :
\begin{equation}\label{inegderniere}
\be_x ( V_\Delta (\theta_j) )
  \leq  \{\!\!\{ \be_y( \Delta ) + \mf d_j (\|\Delta \|_y)\}\!\!\}\,,
\qquad  \|  V_\Delta (\theta_j)\|_x
  \leq   c_j \|\Delta \|^\lambda _y\,.
\end{equation}
\noindent Pour cela fixons nous allons utiliser les propri\'{e}t\'{e}s de
``l'application de Dulac'' du passage du col. Nous rappellerons
bri\`{e}vement cette notion et, pour plus de d\'{e}tails, nous renvoyons au
chapitre 7.2. du livre de F. Loray \cite{Loray}.\\

%\subsection{\'{E}tape 4: perte de rugosit\'{e} pour l'application de Dulac}\label{etapedulac}

\noindent{\textbf{Etape 2 : perte de rugosit\'{e} \corjf{par}
l'application de Dulac.}\label{etapedulac}} Fixons $j$ et d\'{e}signons
par $\zeta_j : V_\Delta \setminus V_\Delta(\theta_j)\rightarrow
\mathbb{D}\times \{  e^{i\theta_j} \} $ l'application de transport
holonome $\zeta_j(x, e^{i\theta}) := \phi_{i (\theta_j -
\theta)}^{Y}\,(x, \,e^{i\theta})$. D\'{e}signons aussi par $\wt
\D(\varepsilon_4)$  la surface de Riemann du disque ferm\'{e} \'{e}point\'{e}
$\{1\}\times \overline{\D}(\varepsilon_4)^\ast$. Nous notons
$r\underline{e}^{i\theta}$, $0<r\leq \varepsilon_4$, $\theta \in \R$
les points de $\wt \D(\varepsilon_4)$ et
$\chi(r\underline{e}^{i\theta}) := ( 1 , r{e}^{i\theta})$,
l'application de rev\^{e}tement. L'application holomorphe
$$\mf D_j := \zeta_j \circ \Psi^{-1}\circ \chi_{|\mc
S_j} : \mc S_j := \{r\underline{e}^{i\theta}\;|\; 0<r\leq
\varepsilon_4,\,\theta_j<\theta<\theta_j+2\pi \} \rightarrow
V_\Delta (\theta_j)\,,$$ avec $\Psi $  d\'{e}signant toujours
l'application de passage du col (\ref{defin41}), se prolonge au bord
du domaine. Plus pr\'{e}cis\'{e}ment, on voit facilement que $\mf
D_j(r\underline{e}^{i\theta})$ tend vers $\Psi^{-1}(re^{i\theta_j})$
pour $\theta\rightarrow\theta_j$ et  $\mf
D_j(r\underline{e}^{i\theta})$ {tend} vers $
h_j(\Psi^{-1}(re^{i\theta_j}))$ pour $\theta\rightarrow\theta_j+2\pi
$, o\`{u} $h_j := \phi_{2i\pi }^{Y}$ est l'holonomie de $\wF$ le long du
cercle $C'$. Ainsi en posant
$$
\wt{\mathfrak D}_j(1, r \underline{e}^{i(\theta +  2\pi{n} )}) :=
h_j^{\circ {n}}( \mf D_j(r \underline{e}^{i\theta}))\,,\quad {n}\in
\mathbb{Z}\,,\quad h_j := \phi_{2i\pi }^{Y}\,,
$$ on obtient un prolongement analytique de $\mf D_j$, not\'{e}
$\widetilde{\mathfrak{D}}_j$, qui peut \^{e}tre d\'{e}fini au voisinage d'un
ferm\'{e} $\wt{\mc S}_j\supset\mc S_j$ du type suivant :
$$\wt{\mc S}_j :=  \{\, r\underline{e}^{i\theta}\;/ \; 0 < r
\leq  \xi_j (\theta) \, \}\subset \wt \D(\varepsilon_4)\,,\quad
\hbox{\rm avec $ \xi_j : \R \rightarrow \R_{>0}\,$ continue}\,.$$
L'application $\wt{\mf D}_j : \wt{\mc S}_j \rightarrow V_\Delta
(\theta_j)$ s'appelle \textbf{application de Dulac de $\F$, r\'{e}alis\'{e}e
sur la transversale $\wt{\mc S}_j$ et \`{a} valeurs dans
$\mathbb{D}\times \{
 e^{i\theta_j} \}$}. Il est clair que pour tout $ r\underline{e}^{i\theta}\in
 \wt{\mc S}_j$ les points $(1, r e^{i\theta})$ et $\wt{ \mf D}_j
 (r\underline{e}^{i\theta})$ sont sur une m\^{e}me feuille du collier
 $\mc
Col(\wF)$.

\begin{lema}\label{unicite.dulac}Soit
$G: \Omega \rightarrow V_\Delta (\theta_j)$ une application continue
v\'{e}rifiant les propri\'{e}t\'{e}s : a) $\Omega$ est un sous-ensemble connexe
de $ \wt{\mc S}_j $, b) pour tout $m\in \Omega $ les points
$\chi(m)$ et $G(m)$ appartiennent \`{a} la m\^{e}me feuille de $\mc
Col(\wF)$, c) il existe un point $P_0\in \Omega $ tel que $G(P_0)=
  \wt {\mf D}_j(P_0)$.
Alors $G \equiv \wt {\mf D}_j{|}_{\Omega }$.
\end{lema}
\begin{dem}
On sait  que chaque feuille de  $\wF_{|\mc Col(\wF)}$ intersecte
$V_\Delta (\theta_j)$ suivant l'orbite de $h_j$ qui est un ensemble
discret. Ainsi, si l'on se donne la valeur de $G$ ou de $\wt {\mf
D}_j$ en un point, les valeurs aux points voisins sont enti\`{e}rement
d\'{e}termin\'{e}es par la propri\'{e}t\'{e} b). L'ensemble des points de $\Omega $
o\`{u} $G$, et $\wt {\mf D}_j$ co\"{\i}ncident est donc un ouvert. La
connexit\'{e} de $\Omega $ permet de conclure.
\end{dem}

D\'{e}signons par $\alpha_j(\theta) =: (1,
\varrho_j(\theta)e^{i\theta})$, $\theta\in [\theta_j, \theta_j+2\pi
] $ un chemin a. l. p. m. qui param\'{e}trise simplement $\partial
\Delta $ dans le sens direct et par $\wt \alpha_j : \R
\hookrightarrow \wt{\mathbb{D}}(\varepsilon_4)$,  son relev\'{e}  $\wt
\alpha_j(\theta):= \varrho_j(\theta)\underline{e}^{i\theta}$. Nous
allons admettre provisoirement les lemmes ci-dessous et prouver les
in\'{e}galit\'{e}s
(\ref{inegderniere}).\\

\begin{lema}\label{def1dulac} Pour $\|\Delta \|_x$ assez petit,
il existe  $\vartheta _j \in \mathbb{R}_{>0}$ tels que $\wt{\mf
D}_j\circ \wt \alpha _j ([\theta_j , \theta_j + \vartheta_j]) =
\partial V_\Delta (\theta_j)$ et $0<\vartheta_j \leq 2\pi (\lambda + 1)$.
\end{lema}
\begin{lema}\label{slemmelimdulac} Lorsque $y=r\underline{e}^{i\theta}$
tend vers $0$ avec $\theta$ variant dans un intervalle compact,
$\big(e^{i\theta_j}\wt {\mf D}_j(y)\big/ y^\lambda \big)$ tend vers
uniform\'{e}ment vers $1$ et  $\big({y\wt {\mf D}_j'(y)}\big/{\wt{\mf
D}_j(y)}\big)$ tend uniform\'{e}ment vers $\lambda $.
\end{lema}

Le bord de $V_\Delta (\theta_j)$ est l'image du chemin $\wt \alpha
_j$ par l'application $\wt{\mf D}_j$, d'apr\`{e}s le premier lemme
ci-dessus. La seconde in\'{e}galit\'{e} de (\ref{inegderniere}) r\'{e}sulte donc
directement du second lemme. La formule de composition analogue \`{a}
(\ref{com.rugos}) donne la premi\`{e}re in\'{e}galit\'{e} de
(\ref{inegderniere}) en prenant pour $\mf d_j(r)$ le maximum de la
valeur absolue de l'argument de $ \,{y\wt {\mf
D}_j'(y)}\big/{\wt{\mf D}_j(y)}$ pour $y=r\underline{e}^{i\theta}$,
avec $|y| \leq |r|$ et $\theta_j \leq \theta\leq \theta_j + 2\pi
(\lambda +1)$. D'apr\`{e}s (\ref{slemmelimdulac}) comme $\lambda $ est
r\'{e}el $>0$, $\mf d_j(r)$ tend vers $0$ lorsque $r$ tend vers $0$.

\begin{dem2}{du Lemme \ref{def1dulac}} Les feuilles de la restriction
de $\F$ \`{a} $\partial \D \times \mathbb D^*$ sont transverses aux
surfaces $T_\theta := \{ e^{i\theta}\}\times \C^*$, $\theta\in [0,
2\pi [$. Elles sont aussi transverses aux surfaces $T^\theta \subset
\partial \D \times \mathbb D^*$ d'\'{e}quation $\arg(y)= \theta$; pour
le voir on peut remarquer que $\F_{|\mathbb{D}^*\times
\partial \D}$ est transverse aux surfaces $T'_\theta \subset
\mathbb{D}^*\times
\partial \D$ d'\'{e}quations $\arg (y) = \theta$ et que
l'application $\Psi$ de passage du col \'{e}change ces deux feuilletages
et transforme $T'_\theta$ en $T^\theta$. Ainsi, \'{e}tant donn\'{e} un
chemin
\begin{equation}\label{chemfeuille}
    t\mapsto\phi^X_{it}(1, y_0)\,,\quad t\in [0, t_1]
\end{equation}
d'extr\'{e}mit\'{e}s  $(1, y_0) \in T_0$ et $m_1 = \phi^X_{it_1}(1, y_0)\in
T^{\theta_j}$, il existe une unique fonction analytique r\'{e}elle
$\varsigma$, d\'{e}finie sur un voisinage de $m_0$ dans $T^{\theta_j}$
et \`{a} valeurs dans $\R_{<0}$, telle que $\phi^X_{i\varsigma(m)}(m)\in
T_0$ et $\varsigma(m_1)= - t_1$. L'application $m\mapsto
\phi^X_{i\varsigma(m)}(m)$ est un diff\'{e}omorphisme local
$\R$-analytique de $(T_0,m_0)$ sur $(T^{\theta_j}, m_1)$.

D'autre part $\partial^{sat} U_\Delta $ est l'ensemble des points
$\phi^X_{it}(1, y)$ avec $(1, y)\in
\partial \Delta $ et $0<t<2\pi $. Ainsi tout point de
$\Gamma _j  := \partial^{sat} U_\Delta \cap T^{\theta_j}$ est reli\'{e}
dans $\partial \D\times \mb D^*$ \`{a} un point de $\partial \Delta $,
par un chemin du type (\ref{chemfeuille}) ci-dessus, avec
$0<t_1<2\pi $. On d\'{e}duit de ce qui pr\'{e}c\`{e}de que $\Gamma _j $ est
localement diff\'{e}omorphe \`{a} $\partial \Delta $. Plus pr\'{e}cis\'{e}ment
$\Gamma _j $ est un chemin analytique par morceaux, hom\'{e}omorphe \`{a} un
intervalle ouvert, trac\'{e} dans $T^{\theta_j}$.  Le transport holonome
induit une application ``localement injective'' $\beta _j:\Gamma
_j\rightarrow\partial \Delta $. Notons $P_j$ le point de $\partial
\Delta $ d'argument $\theta_j$. Visiblement $\{P_j\} =
\overline{\Gamma _j}\cap T_0 =\partial U_\Delta \cap T^{\theta_j}$.
Soit $\wh \Gamma _j$ la compactification de $\Gamma _j$ hom\'{e}omorphe
\`{a} un intervalle ferm\'{e}. Nous pouvons noter $\widehat{\Gamma}_j :=
\Gamma _j\sqcup\{P_j, P'_j\}$, avec l'application $\beta _j$ qui se
prolonge  continument en une application $\widehat{\beta} _j :
\widehat{\Gamma}_j \rightarrow
\partial \Delta $ telle que $\widehat{\beta} _j(P_j)=
\widehat{\beta} _j(P'_j) = P_j$. Celle-ci se factorise en une
injection $\wt \beta _j : \widehat{\Gamma}_j\hookrightarrow \wt{\mc
S}_j$, si $\Delta $ est de taille assez petite.

Remarquons que l'application
$$\mc D := (\Psi^{-1}\circ \wt\beta _j^{-1})\; :\; \Omega
\longrightarrow V_\Delta (\theta_j)\,,\quad \Omega := \wt \beta
_j(\widehat{\Gamma}_j)\,,$$ satisfait les hypoth\`{e}ses du lemme
(\ref{unicite.dulac}). En effet a) et b) sont {\'{e}videntes} par
construction et c)  r\'{e}sulte de la remarque suivante :
$\Psi^{-1}(re^{i\theta})= \wt{\mf D}_j (\underline{e}^{i\theta})$ si
$\theta = \theta_j$. Ainsi $\wt{\mf D}_j (\Omega )=
\Psi^{-1}(\overline{\Gamma _j})$. Mais par construction
$\Psi^{-1}(\overline{\Gamma _j}) = \partial V_\Delta (\theta_j)$.
Pour obtenir une \'{e}galit\'{e} $\wt{\mf D}_j\circ \wt \alpha _j ([\theta_j
, \theta_j + \vartheta_j]) =
\partial V_\Delta (\theta_j)$, il suffit de remarquer que
$\Omega$ est l'image d'un intervalle par $\wt \alpha _j$.

La longueur $\vartheta_j$ de cet intervalle est major\'{e} par $2\pi
d_j$, o\`{u} $d_j$ est le nombre maximum de {points} d'une fibre de
l'application $\beta _j$. Celui-ci se majore par la variation de
l'argument de $y$ dans un chemin $\gamma _{y_0}:
t\mapsto\phi^X_{it}(1, y_0)$, $t\in [0,\, 2\pi ]$, avec $y_0\in
\partial \Delta $. Ce qui donne : $$d_j \leq Re \left( \frac{1}{2i\pi
}\int_{\gamma _{y_0}}\frac{dy}{y}\right) = Re \left( \frac{1}{2i\pi
}\int_{\gamma _{y_0}}(\lambda
 + xy A(x, y))\frac{dx}{x}\right)\leq \lambda +  \epsilon_\Delta \,,$$
o\`{u} $\epsilon_\Delta $ est le maximum de $|y A(x, y))|$ pour $(x,
y)\in U_\Delta $. D'o\`{u} la conclusion.
\end{dem2}

\begin{dem2}{du Lemme (\ref{slemmelimdulac})}
Remarquons d'abord que si l'on compose $\wt {\mf D}_j$ par
l'application de transport holonome $\phi^Y_{-i\theta_j}$, on
obtient l'application de Dulac r\'{e}alis\'{e}e encore sur $\wt{\mc S}_j$
mais \`{a} valeur dans $\D\times \{1\}$. Toujours d'apr\`{e}s (\ref{prop31})
et (\ref{cor.est.rug.}) il revient au m\^{e}me de prouver le lemme pour
$\wt{\mf D}_j$ ou pour $\mf D$. Ainsi nous supposerons $\theta_j =
0$. L'application $\wt {\mf D}_j$ est alors l'application de Dulac
d\'{e}crite dans le chapitre 7 de \cite{Loray} et la premi\`{e}re
affirmation du lemme est \'{e}nonc\'{e} \`{a} la page 148 de cette r\'{e}f\'{e}rence.
\\

Consid\'{e}rons la forme normale formelle, du germe de $\F$ au point
singulier. Elle est donn\'{e}e par une $1$-forme qui s'\'{e}crit :
\begin{equation}\label{formenormale}
    \omega _{p_0/q_0,\,k,\, \alpha  } := q_0\, x\, (1 + \alpha\, x^{kp_0} y^{kq_0})\,dy
    + p_0\, y \,(1+ (\alpha -1)\, x^{k p_0} y^{kq_0})\, dx\,,
\end{equation}
avec $\alpha \in \C\,$, $ k \in \mathbb{N}^\ast$, d\'{e}finissant un
feuilletage not\'{e} $\FN$. D'apr\`{e}s \cite{MaRa} il existe des
applications $\Phi_l : U_l (r_0,\, \varepsilon) \rightarrow
\mc{Col}(\Delta )$, $l=0,\ldots , 2k-1$, d\'{e}finies et diff\'{e}rentiables
(au sens de Withney) sur les secteurs ferm\'{e}s
$$U_l (r_0,\, \varepsilon) := u^{-1}\left\{ re^{i\theta}\;/\;0
\leq r \leq r_0,\, \left|\theta - \frac{\pi }{2k} - l\frac{\pi }{k}
\right| \leq \pi /k \right\} \;\cap\; \D\times \D\,,
$$
holomorphes sur les secteurs ouvert $\inte{U}_l (r_0,\,
\varepsilon)$, tangentes \`{a} l'identit\'{e} en chaque point des axes et
qui conjuguent $\FN$ \`{a} $\F_{| \mc{Col}(\Delta )}$. Soit $m$, $n \in
\mathbb{N}$ des entiers qui satisfont : $mp_0 - n q_0 = 1$.
L'application multiforme
$$H := x^n y^m (x^{p_0} y^{q_0}) ^{\frac{\alpha  - m p_0}{p_0q_0}} \exp\left(\frac{-1}
{p_0q_0 u(x,y)^k}\right)$$ est, sur $\C^2$ priv\'{e} des axes, une
int\'{e}grale premi\`{e}re de $\F_{k, \alpha }$. Elle est uniforme sur $U_l
(r_0,\, \varepsilon)$ et s\'{e}pare les feuilles de la restriction de
$\FN$ \`{a} $U_l (r_0,\, \varepsilon)$, cf. \cite{MaRa} page 598. Ainsi,
si deux points $(1, y)$ et $(x,  e^{i\theta_0})$ appartiennent \`{a}
$U_l (r_0,\, \varepsilon)$ et satisfont l'\'{e}quation
\begin{equation}\label{equationdulac}
    H(1, y) = H(x,   e^{i\theta_0})\,,
\end{equation}
alors il existe un chemin dans une feuille de $\FN|_{U_l (r_0,\,
\varepsilon)}$ qui joint ces deux points. On en d\'{e}duit que toute
solution de (\ref{equationdulac}) est une d\'{e}termination de
l'application de Dulac de $\FN$. Il en r\'{e}sulte aussi que les
applications de conjugaison $\Phi_l$, $l=0,\ldots , 2k-1$, induisent
un syst\`{e}me complet de conjugaisons sectorielles (\`{a} la source et au
but) entre des d\'{e}terminations de l'application de Dulac de $\FN$  et
de $\F$. Ces conjugaisons sont holomorphes sur les secteurs ouverts
et poss\`{e}dent un d\'{e}veloppement assymptotique tangent \`{a} l'identit\'{e} \`{a}
l'origine. On en d\'{e}duit, gr\^{a}ce \`{a} (\ref{prop31}) et
(\ref{cor.est.rug.}), qu'il suffit de prouver le lemme
(\ref{slemmelimdulac}) pour le feuilletage $\FN$. \\

Supposons d\'{e}sormais $\F = \FN$. Pour simplifier les calculs,
remarquons d'abord que l'on peut se ramener au cas $p_0=q_0=k = 1$.
En effet $\FN$ est l'image r\'{e}ciproque de $\FNU$ par l'application
$R(x, y):=(x^{kp_0}, y^{kq_0})$. Celle-ci est un rev\^{e}tement en
restriction \`{a} $\C^*\times \C^*$ et ses restrictions \`{a} des secteurs
appropri\'{e}s contenus dans des transversales aux axes conjuguent les
applications de Dulac respectives. Les conjugantes $g$ sont du type
$\varphi(z) = c z^\beta $ avec $z=x$ ou $y$ et $\beta\in \mb
Q_{>0}$, qui ne modifient pas la rugosit\'{e} des courbes. Il en d\'{e}coule
que le lemme (\ref{slemmelimdulac}) est satisfait par $\FN$ d\`{e}s
qu'il l'est par
$\FNU$.\\

Supposons finalement $\F =\FNU$ et $\theta_j=0$. En d\'{e}rivant la
relation (\ref{equationdulac}) on voit que le graphe l'application
de Dulac $x = \wt{\mf D}_j(y)$ est une solution de l'\'{e}quation
diff\'{e}rentielle :
\begin{equation}\label{equadiffdulac}
\frac{y dx}{x dy} = \frac{x(1+\alpha  y)}{y(1 + (\alpha - 1) x)}\,.
\end{equation}
Comme $\wt{\mf D}_j(y)$ tend vers $0$ lorsque $y$ tend vers $0$ avec
un argument born\'{e}, on d\'{e}duit de cette \'{e}quation que $\big( {y\wt{\mf
D}_j'(y) }\big/{\wt{\mf D}_j(y)}\big)$ converge alors et poss\`{e}de la
m\^{e}me limite $1$, que $\big( \wt{\mf D}_j(y)\big/ y\big)$.
\end{dem2}

\end{document}